\documentclass[10pt,a4paper]{article}

\usepackage{amsmath,amssymb, bbm}
\usepackage{color}
\newcommand{\R}{\mathbb{R}}
\newcommand{\N}{\mathbb{N}}
\newcommand{\Z}{\mathbb{Z}}

\newcommand{\Co}{\mathbb{C}_0}
\newcommand{\Do}{\mathbb{D}_0}
\newcommand{\bP}{\mathbb{P}}
\def\cA{{\mathcal A}}
\def\cC{{\mathcal C}}
\def\cB{{\mathcal B}}

\def\cF{{\mathcal F}}
\def\cP{{\mathcal P}}
\def\cQ{{\mathcal Q}}
\newcommand{\ee}{\varepsilon}
\renewcommand{\aa}{\alpha}
\newcommand{\bb}{\beta}

\renewcommand{\div}{{\rm div}\,}

\newcommand{\Sum}{\displaystyle \sum}
\newcommand{\Hs}{\dot{H^s}}

\def\d{\partial}
\def\ddl{\dot \Delta_l}
\def\ddj{\dot \Delta_j}

\def\ddk{\dot \Delta_k^h}
\def\tilde{\widetilde}
\def\hat{\widehat}
\newcommand{\D}{\Delta}

\newcommand{\n}{\nabla}

\newcommand{\Fe}{F_{ext}}

\newcommand{\voe}{v_{0,\ee}}

\newcommand{\Ue}{U_\ee}
\newcommand{\ve}{v_\ee}
\newcommand{\Ve}{V_\ee}
\newcommand{\He}{H_\ee}
\newcommand{\Ven}{V_\ee^n}
\newcommand{\Hen}{H_\ee^n}

\newcommand{\Uoe}{U_{0,\ee}}
\newcommand{\UoeS}{U_{0,\ee,S}}
\newcommand{\Uoeosc}{U_{0,\ee, osc}}
\newcommand{\tU}{\tilde{U}}

\newcommand{\De}{D_\ee}
\newcommand{\Den}{D_\ee^n}
\newcommand{\Ee}{E_\ee}
\newcommand{\Een}{E_\ee^n}
\newcommand{\Fen}{F_\ee^n}

\newcommand{\Dosc}{D_{\ee, osc}}

\newcommand{\Phie}{\Phi_\ee}
\newcommand{\Thee}{\theta_\ee}
\newcommand{\tThe}{\tilde{\theta}}
\newcommand{\tThee}{\tilde{\theta}_\ee}
\newcommand{\Theeo}{\theta_{0,\ee}}
\newcommand{\tTheo}{\tilde{\theta}_0}
\newcommand{\tTheeo}{\tilde{\theta}_{0,\ee}}

\newcommand{\tv}{\tilde{v}}

\newcommand{\tvo}{\tilde{v}_{0}}

\newcommand{\tpi}{\tilde{\pi}}

\newcommand{\tge}{\tilde{g}_{\ee}}
\newcommand{\tG}{\tilde{G}}
\newcommand{\tZe}{\tilde{Z}_{\ee}}
\newcommand{\tKe}{\tilde{K}_{\ee}}

\newcommand{\qe}{q_{\ee}}

\newcommand{\We}{W_\ee}
\newcommand{\Wet}{W_\ee^T}
\newcommand{\de}{\delta_{\ee}}

\renewcommand{\Re}{R_\ee}
\newcommand{\re}{r_\ee}

\newcommand{\cPrR}{\cP_{\re, \Re}}
\newcommand{\cPrRb}{\cP_{\frac{\re}2, 2\Re}}

\newcommand{\exi}{(\ee,\xi)}

\newtheorem{thm}{Theorem}[section]
\newtheorem{lem}{Lemma}[section]

\newtheorem{prop}{Proposition}[section]
\newtheorem{defi}{Definition}[section]

\newtheorem{rem}{Remark}[section]
\title{New asymptotics for strong solutions of the strongly stratified Boussinesq system without rotation and for large ill-prepared initial data}

\author{Fr\'ed\'eric Charve\footnote{Univ Paris Est Cr\'eteil, CNRS, LAMA, F-94010 Cr\'eteil, 2 Univ Gustave Eiffel, LAMA, F-77447 Marne-la-Vall\'ee, France. E-mail: frederic.charve@u-pec.fr}}

\date{}
\begin{document}
%\tableofcontents

\maketitle

\begin{abstract}
In our previous work dedicated to the strongly stratified Boussinesq system, we obtained for the first time a limit system (when the froude number $\ee$ goes to zero) that depends on the thermal diffusivity $\nu'$ (other works obtained a limit system only depending on the visosity $\nu$). To reach those richer asymptotics we had to consider an unusual initial data which is the sum of a function depending on the full space variable and a function only depending on the vertical coordinate, and we studied the convergence of the weak Leray-type solutions. In the present article we extend these results to the strong Fujita-Kato-type solutions. In this setting, and compared to the case of weak solutions, we obtain far better convergence rates (in $\ee$) for ill-prepared initial data with very large oscillating part of size some negative power of the small parameter $\ee$. The main difficulties come from the anisotropy induced by the presence of $x_3$-depending functions.
\end{abstract}


\section{Introduction}

\subsection{Geophysical fluids: Strongly stratified Boussinesq system}

The strongly stratified Boussinesq system (without rotation) describes the motion of a geophysical fluid submitted to the influence of the gravity through the vertical stratification of the density. In the whole space, this model is written as follows:

\begin{equation}
\begin{cases}
\d_t \Ue +\ve\cdot \n \Ue -L \Ue +\frac{1}{\ee} \cB \Ue=\frac{1}{\ee} (-\n \Phie, 0),\\
\div \ve=0,\\
{\Ue}_{|t=0}=U_{0,\ee}.
\end{cases}
\label{Stratif}
\tag{$S_\ee$}
\end{equation}
The unknowns are $\Ue =(\ve, \Thee)=(\ve^1, \ve^2, \ve^3, \Thee)$, where $\ve$ denotes the velocity of the fluid and $\Thee$ the scalar potential temperature (linked to the density, temperature and salinity), and $\Phie$, which is still called the geopotential, and can  be decomposed as the sum of the pressure term and another penalized gradient term that could be seen as an analoguous of the centrifugal force (we refer to the introductions of \cite{FCPAA} and \cite{FCStratif1} for a more precise presentation of the model).

The diffusion operator $L$ takes into account two heat regularization effects  and is defined by
$$
L\Ue \overset{\mbox{def}}{=} (\nu \D \ve, \nu' \D \Thee),
$$
where $\nu, \nu'>0$ respectively denote the kinematic viscosity and thermal diffusivity (both will be called viscosities in the present article). The last term $\ee^{-1}\cB \Ue$ only takes into account stratification effects and $\cB$ is defined as the following skewsymmetric matrix:
$$
\cB \overset{\mbox{def}}{=}\left(
\begin{array}{llll}
0 & 0 & 0 & 0\\
0 & 0 & 0 & 0\\
0 & 0 & 0 & 1\\
0 & 0 & -1 & 0
\end{array}
\right).
$$
\begin{rem}
 \sl{System \eqref{Stratif} is obtained from the Primitive system only considering the Froude number (introduced by physicists to measure the importance of the stratification effect in the motion). As the rotating fluids and Primitive systems, this model belongs to the family of variations of the famous Navier-Stokes system showing better behaviour induced by the special structure brought by their respective penalized terms as $\ee$ goes to zero. We refer to \cite{FCStratif1} for more details about the geophysical fluids models and a survey on results about the rotating fluids system (\cite{CDGG, CDGG2, CDGGbook, Dutrifoy2}) the Primitive system (\cite{FC1, FC2, FCPAA, FCRF}) and System \eqref{Stratif} (see \cite{Wid, LT, T2, Scro3}).}
\end{rem}

\subsection{Notations}

For an $\R^3$ or $\R^4$-valued vector field, we will write $f^h=(f^1,f^2)$ and will define $f\cdot \n f=\sum_{i=1}^3 f^i \d_i f$. So that for instance we will indifferently write $\ve\cdot \n\Ue=\Ue \cdot \n\Ue$.

We will use the same notations as in \cite{FCPAA, FCcompl, FCStratif1}: for $s\in \R$ and $T>0$ we define the spaces:
$$
\begin{cases}
\vspace{0.2cm}
 \dot{E}_T^s=\mathcal{C}_T(\Hs (\R^3)) \cap L_T^2(\dot{H}^{s+1}(\R^3)),\\
 \dot{B}_T^s=\mathcal{C}_T(\dot{B}_{2,1}^s (\R)) \cap L_T^1(\dot{B}_{2,1}^{s+2}(\R)),
\end{cases}
$$
endowed with the following norms, where $\nu_0=\min(\nu,\nu')$):
$$
\begin{cases}
\vspace{0.2cm}
 \|f\|_{\dot{E}_T^s}^2 \overset{def}{=}\|f\|_{L_T^\infty \Hs }^2+\nu_0 \int_0^T \|f(\tau)\|_{\dot{H}^{s+1}}^2 d\tau,\\
 \|f\|_{\dot{B}_T^s} \overset{def}{=}\|f\|_{L_T^\infty \dot{B}_{2,1}^s}+\nu' \int_0^T \|f(\tau)\|_{\dot{B}_{2,1}^{s+2}} d\tau,
\end{cases}
$$
where $H^s(\R^3)$, $\dot{H}^s(\R^3)$ and $\dot{B}_{2,1}^s(\R)$ respectively denote the inhomogeneous and homogeneous Sobolev spaces of index $s\in \R$ and the homogeneous Besov space of indices $(s,2,1)$.\\
When $T=\infty$ we simply write $\dot{E}^s$ or $\dot{B}^s$ and the corresponding norms are understood as taken over $\R_+$ in time.

\subsection{The limit system}

The present article is the companion paper of \cite{FCStratif1} and focusses on the same question in the context of strong solutions. Let us recall that in \cite{FCStratif1}, we constructed and studied weak solutions to System \eqref{Stratif} which converge, as $\ee$ goes to zero, to a limit truly depending on $\nu'$ (this was not the case in previous papers, see for instance \cite{Scro3}). More precisely, we explained for the first time how we can formally obtain a limit described by the following two systems:
\begin{equation}
    \begin{cases}
  \d_t \tv^h +\tv^h \cdot \n_h \tv^h -\nu \D \tv^h & = -\n_h \tpi^0,\\
  \div_h \tv^h=0,\\
  \tv^h_{|t=0}= \tvo^h,
  \end{cases}
 \label{SNS3}
\end{equation}
and
\begin{equation}
 \begin{cases}
  \d_t \tThe-\nu' \d_3^2 \tThe = 0,\\
  \tThe_{|t=0}= \tTheo,
 \end{cases}
 \label{SNS4}
\end{equation}
where
\begin{equation}
 \tpi^0 =-\sum_{i,j=1}^2 \D_h^{-1} \d_i \d_j (\tv^i \tv^j).
\label{defq}
 \end{equation}
This suggested to consider initial data of the form:
$$
{\Ue}_{|t=0}(x) =\Uoe(x) +(0,0,0,\tTheeo (x_3)),
$$
connected to the previous systems according to:
\begin{equation}
 \begin{cases}
\bP_2 \Uoe^h (x) \underset{\ee \rightarrow 0}{\longrightarrow} \tvo^h (x), \quad \mbox{or equivalently } \bP_2 \Uoe (x) \underset{\ee \rightarrow 0}{\longrightarrow} (\tvo^h (x),0,0),\\
 \tTheeo (x_3) \underset{\ee \rightarrow 0}{\longrightarrow} \tTheo (x_3),
\end{cases}
\label{Condinit}
\end{equation}
where the projector $\bP_2$ is related to the structure of the limit system and is described in the following section.

\begin{rem}
 \sl{Had we only considered ${\Ue}_{|t=0}(x) =\Uoe(x)$ (it is a conventional initial data when $\tTheo (x_3)=0$), the limit would be $(\tv^h,0,0)$ (independant of $\nu'$). The same phenomenon occurs for the rotating fluids system: as explained in \cite{CDGG, CDGGbook} if the initial data is $u_0(x)$ instead of $u_0(x)+\bar{u}_0(x_1,x_2)$ we only obtain that the limit is zero (instead of the unique solution of the 2D-Navier-Stokes system with three components and initial data $\bar{u}_0(x_1,x_2)$).}
\end{rem}

\subsection{The Stratified/osc structure}

The structure of the formal limit system suggested us to introduce the following operators:
\begin{defi}(see \cite{FCStratif1})
 \sl{For a $R^4$-valued function, we introduce the following quantity, that we will call its vorticity:
 $$
 \omega(f)=\d_1 f^2 -\d_2 f^1.
 $$
 From this we define the \emph{stratified} and \emph{oscillating (or oscillatory) parts} of f, respectively denoted as $f_S$ and $f_{osc}$, according to:
 \begin{equation}
 f_S= \left(
 \begin{array}{c}
  \n_h^\perp \D_h^{-1} \omega(f)\\
  0\\
  0
 \end{array}\right)
=\left(
 \begin{array}{c}
  -\d_2 \D_h^{-1} \omega(f)\\
  \d_1 \D_h^{-1} \omega(f)\\
  0\\
  0
 \end{array}\right),
 \label{defS}
 \end{equation}
 and, denoting $\div_h f^h \overset{def}{=} \d_1 f^1 + \d_2 f^2$,
\begin{equation}
 f_{osc}=f-f_{S}= \left(
 \begin{array}{c}
  \n_h \D_h^{-1} \div_h f^h\\
  f^3\\
  f^4
 \end{array}\right)
 =\left(
 \begin{array}{c}
  \d_1 \D_h^{-1} \div_h f^h\\
  \d_2 \D_h^{-1} \div_h f^h\\
  f^3\\
  f^4
 \end{array}\right).
 \label{defosc}
 \end{equation}
}
\end{defi}
The following proposition gathers properties of the stratified/oscillating structure which is linked to the spectral properties of the linearized system.
\begin{prop}\cite{FCStratif1}
 \sl{With the notations from \eqref{defS} and \eqref{defosc}, there exist two pseudodifferential operators of order zero $\mathcal{P}$ and $\mathcal{Q}$ such that for any $f$,
$$
f_{S}=\mathcal{Q} f,\quad \mbox{and} \quad f_{osc}=\mathcal{P} f.
$$
These operators satisfy:
\begin{enumerate}
 \item $\mathcal{Q}=\bP_2$ and $\mathcal{P}=I_d-\bP_2$ (where the operators $\bP_k$ are the spectral projectors defined in Proposition \ref{estimvp}).
 \item For any $s\in \R$, we have $((I_d-\bP_2) f| \bP_2 f)_{H^s/\dot{H}^s} =0=(\cB f| \bP_2 f)_{H^s/\dot{H}^s}$ (when defined).
 \item $(I_d-\bP_2) f=f \Longleftrightarrow \bP_2 f =0 \Longleftrightarrow \omega(f)=0$.
 \item $(I_d-\bP_2) f=0 \Longleftrightarrow \bP_2 f =f \Longleftrightarrow f^3=f^4=0 \mbox{ and } \div_h f=0 \Longleftrightarrow$ there exists a scalar function $\phi$ such that $f=(-\d_2 \phi, \d_1 \phi,0,0) =(\n_h^\perp \phi, 0,0)$. Such a vector field is obviously divergence free (and horizontal divergence-free) and we will say that it is \emph{stratified}. It also satisfies $f=(f^h,0,0)$.
\item If $f$ is divergence-free, so is $(I_d-\bP_2) f$.
 \item $\cB \bP_2 f=0$ (in $\R^4$).
 \item $\bP_2 \bP=\bP \bP_2=\bP_2$ and $\bP_2 (I_d-\bP)=(I_d-\bP) \bP_2=0$ (in particular $\bP_2 (\n q,0)=0$).
 \item If f is a divergence-free vector field, then (we recall that we denote $f\cdot \n f=\sum_{i=1}^3 f^i \d_i f$)
 $$
 \omega(f\cdot \n f)=-\d_3 f^3\cdot \omega(f)+\d_1 f^3 \cdot \d_3 f^2 -\d_2 f^3 \cdot \d_3 f^1 +f\cdot \omega (f).
 $$
 \item If f is a stratified vector field, then $\omega(f\cdot \n f)=f\cdot \omega (f)$.
\end{enumerate}}
\label{PropSosc}
\end{prop}
\begin{rem}
\sl{
\begin{enumerate}
 \item As outlined in \cite{FCStratif1}, the previous decomposition is close to the nicer case $\nu=\nu'$ for the Primitive system, in the sense that we have $\cQ=\bP_2$, and $\bP_2$ is orthogonal to $\bP_3$ and $\bP_4$ in the general case (but it is only when $\nu=\nu'$ that $\bP_3$ and $\bP_4$ are also orthogonal projectors of norm $1$).
 \item For a $\R^2$-valued function $f=(f^1,f^2)=f^h$, we could introduce $\bP_2^h$ and $f_S= \bP_2^h f= \n_h^\perp \D_h^{-1} \omega(f)$ (but with a slight notational abuse, we may also denote $f_S=\bP_2 f$ and $f_{osc}=f-f_{S}=\n_h \D_h^{-1} \div_h f$.
 \item As outlined in \cite{FCStratif1}, the previous decomposition leads to the notion of \textbf{well-prepared} or \textbf{ill-prepared initial data}. We say that an initial data  is well-prepared if it is stratified in the sense of Point 4 from Proposition \ref{PropSosc} (which means it has a zero oscillating part) or has a small oscillating part. In other words the initial data \emph{already has the structure of the limit system or is close to it}. On the opposite, an ill-prepared initial data features a large oscillating part (and in the present article, it will be large of size a negative power of $\ee$). Historically for the Primitive or Rotating fluid systems, the first works consisted in studying the system for well-prepared initial data (see for instance \cite{DesGre, Dragos4}). Another case where we are forced to consider well-prepared initial data is when the Rossby and Froude number are equal and we cannot rely anymore on dispersive estimates (case $F=1$, see \cite {Chemin2, FCF1})
\end{enumerate}
}
\end{rem}
Now we can completely precise the initial data and limit system that we will consider in this article:
\begin{equation}
 {\Ue}_{|t=0}(x) =\Uoe(x)
 +\left(\begin{array}{c}
  0\\0\\0\\\tTheeo (x_3)
 \end{array}\right)
 =\UoeS(x) +\Uoeosc(x)
 +\left(\begin{array}{c}
  0\\0\\0\\\tTheeo (x_3)
 \end{array}\right).
 \label{Complete:data}
\end{equation}
And we will denote:
$$
\Uoe=\left(\begin{array}{c}
  \voe^1\\ \voe^2\\ \voe^3\\ \Theeo
 \end{array}\right)
 \quad \mbox{and} \quad
 \Ue=\left(\begin{array}{c}
  \ve^1\\ \ve^2\\ \ve^3\\ \Thee
 \end{array}\right)
 =\left(\begin{array}{c}
  \ve \\ \Thee
 \end{array}\right)
=\left(\begin{array}{c}
  \ve^h\\ \ve^3\\ \Thee
 \end{array}\right).
 $$

\subsection{Reformulation of the systems}

The first step in \cite{FCStratif1} was to rewrite Systems \eqref{SNS3} and \eqref{SNS4} and merge them into a more practical formulation: denoting as $\bP$ the orthogonal Leray projector onto divergence-free vectorfields, and setting $\tU\overset{def}{=}(\tv^h,0,\tThe)$, we obtained that $\tU$ satisfies:
\begin{equation}
\begin{cases}
 \d_t \tU +\tU \cdot \n \tU -L \tU + \frac1{\ee} \cB \tU = -\tG -\left(\begin{array}{c}\n \tge \\0\end{array}\right),\\
 \div \tv=0,\\
 \tU_{|t=0}=(\tvo^h,0, \tTheo).
\end{cases}
\label{SNS6}
%\tag{$\tilde{SNS_6}$}
\end{equation}
where $\tG$ is defined by ($\tpi^0$ has been introduced in \eqref{defq}):
\begin{equation}
 \tG= \bP \left(\begin{array}{c}
  \d_1 \tpi^0\\ \d_2 \tpi^0 \\0\\0 \end{array}\right)
  =\left(\begin{array}{c}
  \d_1 \d_3^2 \D^{-1} \D_h^{-1}\\\d_2 \d_3^2 \D^{-1} \D_h^{-1} \\-\d_3 \D^{-1}\\0 \end{array}\right)\sum_{i=1}^2 \d_i (\tv^h\cdot \n_h \tv^i).
  \label{deftG}
 \end{equation}
Moreover, we emphasize that
\begin{equation}
 \omega(\tG)=0=\div \tG \mbox{ and }\bP_2 \tG=0,\mbox{ and roughly }\tG \sim \n (\tv^h \otimes \tv^h) \sim \tv^h\cdot \n_h \tv^h.
 \label{ecritG}
\end{equation}
As precised in \cite{FCStratif1}, studying the system satisfied by $\Ue-\tU$ will not be possible, because of its initial data which prevent the use of classical results:
$$
(\Ue-\tU)_{|t=0}(x)=\Uoeosc(x)+ (\bP_2 \Uoe^h(x) -\tvo^h(x), 0, \tTheeo(x_3)-\tTheo(x_3)).
$$
In order to properly justify the construction of weak solutions with such initial data, we needed in \cite{FCStratif1} to rewrite System \eqref{Stratif} into a formulation where functions only depending on $x_3$ do not appear in the initial data anymore. Doing this moved these functions in the transport terms which required an adaptation of the proof of the classical Leray theorem.
\\

More precisely, in order to neutralize the $x_3$-only-dependent part, we simply defined the following function:
\begin{equation}
 \tZe=\left(\begin{array}{c}
  0\\0\\0\\ \tKe \end{array}\right), \quad \mbox{where } \tKe \mbox{ solves }
  \begin{cases}
   \d_t \tKe-\nu' \d_3^2 \tKe = 0,\\
  \tThe_{|t=0}= \tTheeo-\tTheo.
  \end{cases}
\end{equation}
and finally set:
\begin{equation}
 \De \overset{def}{=} \Ue-\tU-\tZe =\left(\begin{array}{c}
  \ve^h-\tv^h\\\ve^3\\ \Thee-(\tThe+\tKe) \end{array}\right) =\left(\begin{array}{c}
  \ve^h-\tv^h\\\ve^3\\ \Thee-\tThee \end{array}\right)=\Ue-\left(\begin{array}{c}\tv^h\\ 0\\ \tThee \end{array}\right),
  \label{defD}
\end{equation}
where the function $\tThee\overset{def}{=}\tThe+\tKe$ solves:
\begin{equation}
 \begin{cases}
   \d_t \tThee-\nu' \d_3^2 \tThee = 0,\\
  \tilde{\theta}_{\ee|t=0}= \tTheeo.
  \end{cases}
  \label{SysttThee}
\end{equation}
The results from \cite{FCStratif1} were obtained studying the following system satisfied by $\De=(\Ve,\He)$:
\begin{equation}
\begin{cases}
 \d_t \De -L \De + \frac1{\ee} \cB \De = -\left[\De\cdot \n \De
 +\left(\begin{array}{c}
  \De \cdot \n \tv^h \\0\\ \De^3\cdot \d_3 \tThee \end{array}\right) +\tv^h\cdot \n_h \De \right] +\tG -\left(\begin{array}{c}\n \qe \\0 \end{array}\right),\\
 \div \Ve=0,\\
 D_{\ee|t=0}=\Uoeosc +(\UoeS-(\tvo^h,0,0)) =\Uoeosc +(\UoeS^h-\tvo^h,0,0).
\end{cases}
\label{StratifD}
\end{equation}

Before presenting the results for the weak solutions, we will recall in the next section what we proved in \cite{FCStratif1} for the limit system.

\subsection{Study of the limit system}

Let us begin with System \eqref{SNS4}, which is only a one-dimensional heat equation (we refer for example to \cite{Dbook}, Section 3.4.1, Lemma 5.10 and Proposition 10.3, see also Definition \ref{deftilde}).

\begin{thm} \cite{FCStratif1}
 \sl{Let $s\in \R$. For any $\tTheo\in \dot{H}^s(\R)$ (respectively $\tTheo\in \dot{B}_{2,1}^s(\R)$) there exists a unique global solution $\tThe$ of \eqref{SNS4} and for all $t\geq 0$, we have:
 \begin{equation}
  \|\tThe\|_{\tilde{L}_t^\infty \dot{H}^s}^2 +\nu' \|\tThe\|_{L_t^2 \dot{H}^{s+1}}^2 \leq 2 \|\tTheo\|_{\dot{H}^s}^2.
 \end{equation}
\begin{equation}
(\mbox{respectively} \quad \|\tThe\|_{\tilde{L}_t^\infty \dot{B}_{2,1}^s} +\nu' \|\tThe\|_{L_t^1 \dot{B}_{2,1}^{s+2}} \leq \|\tTheo\|_{\dot{B}_{2,1}^s}.)
 \end{equation}
 More generally for $s\in \R$ and $p,r\in[1,\infty]$, there exists a constant $C>0$ such that if $\tTheo\in \dot{B}_{p,r}^s(\R)$ then for all $q\in [1,\infty]$
 \begin{equation}
  \|\tThe\|_{\tilde{L}_t^q \dot{B}_{p,r}^{s+\frac2{q}}}\leq \frac{C}{(\nu')^\frac1{q}} \|\tTheo\|_{\dot{B}_{p,r}^s}.
   \label{estimThetaBspr}
 \end{equation}
 }
 \label{ThHeat}
\end{thm}
\begin{rem}
 \sl{Thanks to this result, the previously defined $\tThe$, $\tKe$ and $\tThee$ are global and satisfy similar estimates.
}
\end{rem}
On the other hand, we observed in \cite{FCStratif1} that System \eqref{SNS3} is very close to the quasi-geostrophic system (see \cite{FC1, FC2}), and  we easily adapted Theorem 1 from \cite{FCPAA} and obtained the following theorem that generalizes the results from \cite{Scro3} as we need less initial regularity:

\begin{thm} \cite{FCStratif1}
 \sl{Let $\delta>0$ and $\tvo^h\in H^{\frac12+\delta}$ a $\R^2$-valued vectorfield such that $\div_h \tvo^h=0$. Then System \eqref{SNS3} has a unique global solution $\tv^h \in E^{\frac12+\delta}=\dot{E}^0 \cap \dot{E}^{\frac12+\delta}$ and there exists a constant $C=C_{\delta, \nu}>0$ such that for all $t\geq 0$, we have:
 \begin{multline}
  \|\tv^h\|_{L^\infty H^{\frac12+\delta}}^2 +\nu \|\n \tv^h\|_{L^2 H^{\frac12+\delta}}^2 \leq C_{\delta, \nu} \|\tvo^h\|_{H^{\frac12+\delta}}^2 \max (1,\|\tvo^h\|_{H^{\frac12+\delta}}^\frac1{\delta})\\
  \leq C_{\delta, \nu} \max (1,\|\tvo^h\|_{H^{\frac12+\delta}})^{2+\frac1{\delta}},
 \end{multline}
 Moreover, we can also bound the term $\tG$ introduced in \eqref{defD}: for all $s\in[0,\frac12+\delta]$,
 \begin{equation}
  \int_0^\infty \|\tG(\tau)\|_{\dot{H}^s} d\tau \leq C_{\delta, \nu} \max (1,\|\tvo^h\|_{H^{\frac12+\delta}})^{2+\frac1{\delta}}.
 \end{equation}
 }
 \label{ThSNS}
\end{thm}

\subsection{Existence and convergence results for the weak solutions}

We can now state the main results from \cite{FCStratif1}: first the analoguous of the Leray theorem for \eqref{StratifD} which provides global weak solutions for any $\ee>0$:
\begin{thm} (Existence of Leray weak solutions)
 \sl{For any $\delta>0$ and $\Co\geq 1$, let $\tv_0^h\in H^{\frac12+\delta}(\R^3)$ (with $\div_h \tv_0^h=0$), $\tTheeo \in \dot{B}_{2,1}^{-\frac12}(\R)$ (for all $\ee>0$) such that:
 $$
 \|\tv_0^h\|_{H^{\frac12+\delta}(\R^3)}\leq \Co \quad \mbox{and} \quad \sup_{\ee>0} \|\tTheeo\|_{\dot{B}_{2,1}^{-\frac12}(\R)} \leq \Co.
 $$
 Thanks to Theorems \ref{ThHeat} and \ref{ThSNS}, $\tv^h$ and $\tThee$ globally exist (for all $\ee>0$) and respectively belong to $\dot{E}^0 \cap \dot{E}^{\frac12+\delta}$ and $\dot{B}^{-\frac12}$.

 Moreover there exists a constant $C_{\delta,\nu,\nu'}>0$ such that for any fixed $\ee>0$, if $\Uoe\in L^2(\R^3)$, then there exists a weak global solution of \eqref{StratifD} $(\De,\qe)$ with $\De \in \dot{E}^0$ and $\qe \in \dot{E}^1+L^\frac43 (\R_+,L^2)$, satisfying for all $t\geq 0$,
 \begin{multline}
  \|\De(t)\|_{L^2}^2 +\nu_0 \int_0^t \|\n \De (\tau)\|_{L^2}^2 d\tau\\
\leq \left(\|\Uoeosc\|_{L^2}^2 +\|\UoeS^h-\tvo^h\|_{L^2}^2 +C_{\delta,\nu,\nu'} \Co^{2+\frac1{\delta}}\right) e^{C_{\delta,\nu,\nu'} \Co^{2+\frac1{\delta}}}.
\label{estimaprioriLeray}
 \end{multline}
 }
 \label{ThLeray}
\end{thm}

The second result in \cite{FCStratif1} was the convergence result, that rigourously validates as a limit what we formally obtained:
\begin{thm} (Convergence)
 \sl{For any $\delta>0$, $\Co\geq 1$, $\tTheo\in \dot{B}_{2,1}^{-\frac12}(\R)$ , $\tv_0^h\in H^{\frac12+\delta}(\R^3)^2$ (with $\div_h \tv_0^h=0$ or, equivalently, $ \tv_0^h=\bP_2\tv_0^h$) and any $\Uoe\in L^2$ (divergence-free) and $\tTheeo \in \dot{B}_{2,1}^{-\frac12}(\R)$ (for all $\ee>0$) with:
\begin{equation}
  \begin{cases}
  \vspace{0.1cm}
  \|\tv_0^h\|_{H^{\frac12+\delta}(\R^3)}\leq \Co,\\
  \vspace{0.1cm}
  \sup_{\ee>0} \|\Uoe\|_{L^2} \leq \Co,\\
  \|\UoeS^h-\tv_0^h\|_{L^2} \underset{\ee\rightarrow 0}{\longrightarrow} 0,
 \end{cases}
\quad \mbox{and} \quad
\begin{cases}
 \vspace{0.1cm}
 \|\tTheo\|_{\dot{B}_{2,1}^{-\frac12}}\leq \Co,\\
 \sup_{\ee>0} \|\tTheeo\|_{\dot{B}_{2,1}^{-\frac12}(\R)} \leq \Co,\\
 \|\tTheeo-\tTheo\|_{\dot{B}_{2,1}^{-\frac12}(\R)} \underset{\ee\rightarrow 0}{\longrightarrow} 0,
\end{cases}
\end{equation}
the global weak solution $\Ue$ (constructed in Theorem \ref{ThLeray}) converges to $(\tv^h,0,\tThe)$ (where $\tv^h$ and $\tThe$ are the global solutions of Systems \eqref{SNS3} and \eqref{SNS4}) in the following sense: if $\De=\Ue-(\tv^h,0, \tThee)$ (where $\tThee$ is the global solution of \eqref{SysttThee}), then
\begin{itemize}
 \item the stratified part $D_{\ee,S}=\bP_2 \De$ converges to zero: for all $q\in]2,6[$,
$$
\|D_{\ee,S}\|_{L_{loc}^2(\R_+, L_{loc}^q(\R^3)} \underset{\ee\rightarrow 0}{\longrightarrow} 0,
$$
\item the oscillating part $\Dosc =(I_d-\bP_2) \De$ converges to zero: for all $q\in]2,6[$, there exists $\ee_1=\ee_1(\nu,\nu',q)>0$ and, for all $t\geq 0$, a constant $\mathbb{D}_t=\mathbb{D}_{t,\delta,\nu,\nu',q,\Co}$ such that for all $\ee\in]0, \ee_1]$,

\begin{equation}
 \|\Dosc\|_{L_t^2 L^q} =\|\Dosc\|_{L^2([0,t], L^q(\R^3)} \leq \mathbb{D}_t \ee^{\frac{K(q)}{640}},\quad \mbox{with} \quad K(q)\overset{def}{=} \frac{\min(\frac6{q}-1,1-\frac2{q})^2}{ (\frac6{q}-1)}.
 \end{equation}
\end{itemize}
Moreover, when $\nu=\nu'$, the previous estimates can be upgraded into $\|\Dosc\|_{L_t^2 L^q} \leq \mathbb{D}_t \ee^{\frac{K(q)}{544}}$ (now valid for all $\ee>0$) and we can obtain global-in-time estimates with better convergence rate: there exists a constant $C=C_{\nu,\delta,\Co}>0$ such that, for any $\ee>0$,
$$
\|\Dosc\|_{\tilde{L}^\frac43 \dot{B}_{8, 2}^0 +\tilde{L}^1 \dot{B}_{8, 2}^0} \leq C \ee^\frac3{16}.
$$
 }
 \label{ThCV}
\end{thm}

\subsection{Existence and convergence results for the strong solutions}

Let us recall that our initial data is the one stated in \eqref{Complete:data}:
$$
{\Ue}_{|t=0}(x) =\Uoe(x)
 +\left(\begin{array}{c}
  0\\0\\0\\\tTheeo (x_3)
 \end{array}\right)
 =\UoeS(x) +\Uoeosc(x)
 +\left(\begin{array}{c}
  0\\0\\0\\\tTheeo (x_3)
 \end{array}\right).
$$
We are now able to state the main results of the present article. First, the general existence result which is the analoguous of the famous Fujita-Kato theorem (for $\ee>0$ fixed):
\begin{thm}(Existence of local Fujita-Kato strong solutions)
 \sl{Let $\ee>0$, $\delta\in]0,1]$, $\tv_0^h\in H^{\frac12+\delta}(\R^3)$ and $\tTheeo \in \dot{B}_{2,1}^{-\frac12}(\R)\cap \dot{B}_{2,1}^{-\frac12+\bb}(\R)$ (for some fixed $\bb>0$). For any $\Uoe=\UoeS+\Uoeosc \in H^\frac12$, there exists a unique local solution $\De$ of \eqref{StratifD} with lifespan $T_\ee^*>0$ such that for any $T<T_\ee^*$, $\De\in E_T^\frac12=\dot{E}_T^0 \cap \dot{E}_T^\frac12$. Moreover, the following properties are true:
 \begin{itemize}
  \item Regularity propagation: if in addition $\Uoe\in \dot{H}^s$ for some $s\in[0,\frac12+\delta]$ then  for any $T<T_\ee^*$, $\De\in \dot{E}_T^0 \cap \dot{E}_T^s$.
  \item Blow-up criterion: $\int_0^{T_\ee^*} \|\n \De (\tau)\|_{\dot{H}^\frac12}^2 d\tau <\infty \Longrightarrow T_\ee^*=\infty$.
 \end{itemize}
 }
 \label{TH0FK}
\end{thm}
\begin{rem}
 \sl{\begin{enumerate}
      \item The proof of this theorem is postponed to Section \ref{PreuveFK}.
      \item We emphasize that we only state a local existence result, with an unsusual low frequency assumption ($\Uoe\in L^2 \cap \dot{H}^\frac12$) which is needed to treat the additional term $\De^3\cdot \d_3 \Thee$. We will only need the previous blow-up criterion to prove global existence in the main results of the article.
      \item The usual domain for the propagation of regularity is $s\in]-\frac32, \frac32[$, in our case the constraint comes from the regularity of $\tG$.
     \end{enumerate}
 }
\end{rem}
Let us now state a simplified version of the main result of the present article.
\begin{thm}(Global existence and convergence)
 \sl{For all $\nu,\nu',\Co>0$, $\delta\in ]0,\frac18]$, $\tv_0^h\in H^{\frac12+\delta}(\R^3)$ and (for any $\ee>0$) $\Uoe=\UoeS+\Uoeosc \in H^\frac12$, $\tTheo,\tTheeo \in \dot{B}_{2,1}^{-\frac34}(\R)\cap \dot{B}_{2,1}^{-\frac14+\delta}(\R)$ such that for some $\aa_0>0$,
 \begin{equation}
  \begin{cases}
  \vspace{0.1cm}
  \|\tv_0^h\|_{H^{\frac12+\delta}(\R^3)}\leq \Co,\\
  \vspace{0.1cm}
    \|\UoeS^h-\tv_0^h\|_{H^{\frac12+\delta}} \leq \Co \ee^{\aa_0},
 \end{cases}
\quad \mbox{and} \quad
\begin{cases}
 \vspace{0.1cm}
 \|\tTheo\|_{\dot{B}_{2,1}^{-\frac34}(\R)\cap \dot{B}_{2,1}^{-\frac14+\delta}(\R)}\leq \Co,\\
 \|\tTheeo-\tTheo\|_{\dot{B}_{2,1}^{-\frac34}(\R)\cap \dot{B}_{2,1}^{-\frac14+\delta}(\R)} \underset{\ee\rightarrow 0}{\longrightarrow} 0,
\end{cases}
\label{HypTh1}
 \end{equation}
there exist $\ee_0,K,\gamma,c,\Do,q>0$ such that if
\begin{equation}
 \|\Uoeosc\|_{L^q} +\||D|^\frac12\Uoeosc\|_{L^q} +\|\Uoeosc\|_{\dot{H}^{\frac12-c\delta} \cap \dot{H}^{\frac12+\delta}} \leq \Co \ee^{-\gamma},
\end{equation}
then for any $\ee\in]0,\ee_0]$, there exists a unique global strong solution $\Ue$ of \eqref{Stratif} which satisfies $\Ue-(\tv,0, \tThee) \in \dot{E}^0 \cap \dot{E}^{\frac12+\frac{\delta}2}$ and
$$
\|\Ue-(\tv^h,0, \tThee)\|_{L^2 (\R_+, L^\infty(\R^3)} \leq \Do \ee^K.
$$
}
 \label{ThSimple}
\end{thm}

\subsection{Results for the classical Boussinesq system}

We recall that in the companion paper \cite{FCStratif1} (dedicated to the weak Leray solutions), we emphasized that System \eqref{Stratif} is related to the following well-known Boussinesq system:
\begin{equation}
 \begin{cases}
 \d_t v +v\cdot \n v -\nu \D v+\kappa^2 \rho e_3= -\n P,\\
 \d_t \rho + v\cdot \n \rho -\nu' \D \rho =0,\\
 \mbox{div }v=0.
\end{cases}
\label{Bo}
\end{equation}
Let us introduce the following explicit stationnary (and stably vertically stratified) solution of \eqref{Bo} (see other examples of explicit non-stationnary solutions in \cite{FCStratif1}):
\begin{equation}
  \bar{V}_\ee(x)=\left(\begin{array}{c}0 \\ \bar{\rho}_\ee(x)\end{array}\right) =\left(\begin{array}{c}0\\0\\0\\ \bar{\rho}_{0,\ee}-\frac{x_3}{\ee^2\kappa^2}\end{array}\right), \quad \bar{P}_\ee(x)=\bar{P}_{0,\ee}-\kappa^2 \bar{\rho}_{0,\ee} x_3+\frac{x_3^2}{2\ee^2}.
\label{Explicitsol}
  \end{equation}
Then $(V_\ee, P_\ee)$ solves \eqref{Bo} if, and only if, $(\Ue, \Phie)$ solves \eqref{Stratif}, where we have denoted:
\begin{equation}
 V_\ee=\left(\begin{array}{c}\ve \\ \rho_\ee\end{array}\right) =\left(\begin{array}{c}\ve\\ \bar{\rho}_\ee +\frac{\Thee}{\ee \kappa^2}\end{array}\right), \quad \Ue=\left(\begin{array}{c}\ve \\ \Thee\end{array}\right), \quad \frac1{\ee}\Phie= P_\ee-\bar{P}_\ee.
 \label{ChgtvarBouss}
\end{equation}
As in \cite{FCStratif1}, thanks to the change of variables from \eqref{ChgtvarBouss} the previous theorem can be rewritten and provide:
\begin{itemize}
 \item Global existence of strong solutions for the classical Boussinesq system \eqref{Bo} which are perturbations of the previous explicit solution $(\bar{V}_\ee, \bar{p}_\ee)$ (see \eqref{Explicitsol}) and corresponding to non-conventional vertically stratified initial data.
 \item Asymptotic expansion (in $\ee$) of these solutions.
\end{itemize}
More precisely the previous theorem can be reformulated as follows:
\begin{thm}
 \sl{With the assumptions and notations from Theorem \ref{ThSimple}, for any $\ee\in]0,\ee_0]$, there exists a unique strong global solution $V_\ee=(v_\ee,\rho_\ee)$ of \eqref{Bo} corresponding to the following initial data (the last term is $\Uoe$ with a scaling on its last component):
 %$$
 %V|_{t=0}=\Uoe(x) +\left(\begin{array}{c}0\\0\\0\\ -\frac{x_3}{\ee^2 \kappa^2} +\frac1{\ee \kappa^2}\tTheo(x_3)\end{array}\right).
 %$$
 $$
 V_\ee|_{t=0}=\left(\begin{array}{c}0\\0\\0\\ \bar{\rho}_{0,\ee}-\frac{x_3}{\ee^2 \kappa^2} \end{array}\right) +\left(\begin{array}{c}0\\0\\0\\ \frac{\tTheeo(x_3)}{\ee \kappa^2}\end{array}\right) +\left(\begin{array}{c} \UoeS^h(x) +v_{0,\ee,osc}^h(x)\\ v_{0,\ee,osc}^3(x) \vspace{0.2cm}\\ \frac{\theta_{0,\ee,osc}(x)}{\ee \kappa^2}\end{array}\right),
 $$
where $\Uoeosc$ is of size $\ee^{-\gamma}$. Moreover, we have an asymptotic expansion of this solution $V_\ee$ when $\ee$ goes to zero: there exist some $K>0$ and a four-component function $\De$ such that,
 $$
 \|\De\|_{L^2(\R_+, L^\infty(\R^3))}\leq \Do \ee^K,
 $$
and
$$
V_\ee(t,x)=\left(\begin{array}{c}\tv^h(t,x) +\De^h(t,x) \vspace{0.1cm}\\ \De^3(t,x)\\  \bar{\rho}_{0,\ee} -\frac{x_3}{\ee^2 \kappa^2} +\frac{\tThee(t,x_3) +\De^4(t,x)}{\ee \kappa^2}\end{array}\right).
$$
In other words, we have the following expansion:
$$
V_\ee(t,x) \underset{\ee \rightarrow 0}{=} \left(\begin{array}{c}0\\0\\0\\ \bar{\rho}_{0,\ee}-\frac{x_3}{\ee^2 \kappa^2} \end{array}\right) +\left(\begin{array}{c}0\\0\\0\\ \frac{\tThee(t,x_3)}{\ee \kappa^2}\end{array}\right) +\left(\begin{array}{c}\tv^h(t,x)\\ 0\\ 0\end{array}\right) +\left(\begin{array}{c} \mathcal{O}(\ee^K)\\ \mathcal{O}(\ee^K) \\ \mathcal{O}(\ee^K)\\ \mathcal{O}(\ee^{K-1})\end{array}\right).
$$
}
\label{ThBouss}
\end{thm}

\begin{rem}
 \sl{
 \begin{enumerate}
  \item As in \cite{FCStratif1}, we outline that the parameters $\bar{\rho}_{0,\ee}, \bar{P}_{0,\ee}$ and $\kappa$ are free, and we can be choosen depending on $\ee$ the way we wish, and choose for instance $\bar{\rho}_{0,\ee}=\bar{\rho}_0 \ee^{-2}$, $\kappa=\ee^{-1}$ or $\ee^{-\frac12}$.
  \item We refer to the recent article \cite{BBCD} about long-time asymptotics for solutions of the 2D inviscid Boussinesq system (in a periodic strip) near a stably stratified Couette flow.
 \end{enumerate}
 }
\end{rem}

\subsection{Precise statement of the main results}

As in \cite{FC2,FCPAA, FCcompl, FCRF} it is usual that we are not able to obtain convergence results without "removing" some waves. More precisely, due to the presence of the initial oscillating part and of $\tG$ as an independant of $\ee$ external force, any frontal approach with $\De$ is blocked as we could only obtain majorations by quantities independant of $\ee$. We first define the following waves $\We$ and $\Wet$, taylored to "eat" the blocking terms: if $\We$ is the global solution of the following system
\begin{equation}
 \begin{cases}
  \d_t \We -L \We +\frac{1}{\varepsilon} \mathbb{P} \mathcal{B} \We = \tG,\\
  {\We}_{|t=0}=\Uoeosc,
 \end{cases}
\label{We}
 \end{equation}
we also define its frequency truncation on the set $\mathcal{C}_{\re, \Re}$, denoted $\Wet =\cPrR \We$, where the general set $\cC_{r,R}$ is defined in \eqref{CrR}, $\re=\ee^m$ and $\Re=\ee^{-M}$ (the values of $m,M$ will be specified in the statements of the results) and the frequency truncation operator $\cPrR$ is defined in \eqref{PrR} so that $\Wet$ obviously satisfies:
\begin{equation}
 \begin{cases}
  \d_t \Wet -L \Wet +\frac{1}{\varepsilon} \mathbb{P} \mathcal{B} \Wet = \cPrR \tG,\\
  {\We}_{|t=0}=\cPrR \Uoeosc,
 \end{cases}
\label{WeT}
 \end{equation}
We are now able to give a more precise statement of the main results of this article. We emphasize that $\Uoe(x)$ \textbf{is only the conventional part} of the initial data (see \eqref{Complete:data}):
$$
{\Ue}_{|t=0}(x) =\Uoe(x)
 +\left(\begin{array}{c}
  0\\0\\0\\\tTheeo (x_3)
 \end{array}\right)
 =\UoeS(x) +\Uoeosc(x)
 +\left(\begin{array}{c}
  0\\0\\0\\\tTheeo (x_3)
 \end{array}\right).
$$
\begin{thm}(Global existence and convergence, general case)
 \sl{For all $\nu,\nu',\Co>0$, $\delta\in ]0,1]$ $\eta\in]0, \frac12]$ with $\eta \delta \leq \frac13$, $\tv_0^h\in H^{\frac12+\delta}(\R^3)$ and (for any $\ee>0$) $\Uoe=\UoeS+\Uoeosc \in H^\frac12$, $\tTheo,\tTheeo \in \dot{B}_{2,1}^{-\frac34}(\R)\cap \dot{B}_{2,1}^{-\frac14+\delta}(\R)$ satisfying \eqref{HypTh1} for some $\aa_0>0$, there exist $\ee_0,\Do>0$ (depending on $\nu,\nu',\Co,\delta,\eta$) such that for any $\ee\in]0,\ee_0]$, setting $\gamma\overset{def}{=} \frac{\delta}{2784}(1-\eta)$ and $q\overset{def}{=} \frac2{1+\delta}$, if we have
\begin{equation}
 \|\Uoeosc\|_{L^q} +\||D|^\frac12\Uoeosc\|_{L^q} +\|\Uoeosc\|_{\dot{H}^\frac12 \cap \dot{H}^{\frac12+\delta}} \leq \Co \ee^{-\gamma},
 \label{HypTh1osc}
\end{equation}
then there exists a unique global strong solution $\Ue$ of \eqref{Stratif}: the lifespan of $\De$ (given by Theorem \ref{TH0FK}) satisfies $T_\ee^*=+\infty$ and $\De\in \dot{E}^0 \cap \dot{E}^{\frac12+\eta \delta}$. Moreover, if we define $\de=\De-\Wet$ where $\Wet$ is defined in \eqref{WeT} for $(m,M)\overset{def}{=} (\frac1{259}, \frac1{1554})$, then
 $$
 \|\de\|_{\dot{E}^0\cap \dot{E}^{\frac12+\eta \delta}} \leq \Do \ee^{\min(\aa_0,\frac{\delta}{3108}(1-\eta), \frac1{9324})}.
 $$
 If in addition there exists $c>0$ such that
 $$
 \|\Uoeosc\|_{\dot{H}^{\frac12-c\delta} \cap \dot{H}^{\frac12+\delta}} \leq \Co \ee^{-\gamma},
 $$
 then we have:
 $$
 \|\De\|_{L^2 L^\infty} =\|\Ue-(\tv^h,0, \tThee)\|_{L^2 L^\infty} \leq \Do \ee^{\min(\aa_0,\frac{\delta}{3108}(1-\eta), \frac1{9324})}
 $$
 }
 \label{Th1}
\end{thm}
When $\nu=\nu'$ it is usual that some simplifications improve the results, as listed below.
\begin{thm}(Global existence and convergence, case $\nu=\nu'$)
 \sl{Let $\Co>0$, $\delta \in ]0,\frac18]$, $\tv_0^h\in H^{\frac12+\delta}(\R^3)$ and $\Uoe=\UoeS+\Uoeosc \in H^\frac12$, $\tTheo,\tTheeo \in \dot{B}_{2,1}^{-\frac34}(\R)\cap \dot{B}_{2,1}^{-\frac14+\delta}(\R)$ satisfying \eqref{HypTh1} for some $\aa_0>0$.
 \begin{enumerate}
  \item There exist $m_0,\ee_0>0$ such that if for some $c>0$ (as small as we want)
$$
  \|\Uoeosc\|_{\dot{H}^{\frac12-c\delta} \cap \dot{H}^{\frac12+\delta}} \leq m_0 \ee^{-\frac{\delta}2},
$$
then for any $\ee\in]0,\ee_0]$, there exists a global solution of \eqref{Stratif} and $\De\in \dot{E}^0\cap \dot{E}^{\frac12}$.
\item If there exists a function $m(\ee)\underset{\ee \rightarrow 0}{\longrightarrow}0$ such that for some $c>0$
$$
  \|\Uoeosc\|_{\dot{H}^{\frac12-c\delta} \cap \dot{H}^{\frac12+\delta}} \leq m(\ee) \ee^{-\frac{\delta}2},
$$
then if we define $\de=\De-\We$ (with $\We$ solving \eqref{We}), there exists $\Do=\Do(\nu,\Co,\delta)>0$ such that:
$$
\|\de\|_{\dot{E}^0\cap \dot{E}^{\frac12}} \leq \Do \max \big(\ee^{\aa_0}, \ee^{\frac{\delta}2}, m(\ee)\big)\underset{\ee \rightarrow 0}{\longrightarrow}0.
$$
\item Finally, if for some $c>0$ and $\gamma\in]0,\frac{\delta}2[$ we have
$$
  \|\Uoeosc\|_{\dot{H}^{\frac12-c\delta} \cap \dot{H}^{\frac12+\delta}} \leq \Co \ee^{-\gamma},
$$
then
$$
\|\de\|_{\dot{E}^0\cap \dot{E}^{\frac12+\frac{\delta}2-\gamma}} \leq \Do \ee^{\min(\aa_0, \frac{\delta}2-\gamma)},
$$
and for any $k\in]0,1[$ (as close to $1$ as we wish), there exists $\Do=\Do(\nu,\Co,\delta,k)>0$ such that:
$$
\|\De\|_{L^2 L^\infty} =\|\Ue-(\tv^h,0, \tThee)\|_{L^2 L^\infty}\leq \Do \ee^{\min\big(\aa_0, k(\frac{\delta}2-\gamma)\big)}.
$$
 \end{enumerate}
 }
 \label{Th2}
\end{thm}
The article is structured as follows: in the next section we prove Theorem \ref{Th1}, we first obtain apriori estimates then explain the bootstrap method. Section \ref{PrvTh2} is devoted to the proof of Theorem \ref{Th2}, which features better results as $\nu=\nu'$. We postponed to the appendix the proof of the anisotropic Strichartz estimates (which require a technical result from \cite{FCStratif1}) and of Theorem \ref{TH0FK} (which \emph{unusually} relies on a priori estimates in inhomogeneous Sobolev spaces, which are a particular case of the ones obtained Sections \ref{PrvTh1} and \ref{PrvTh2}).

\begin{rem}
 \sl{We emphasize that in this article we use the \textbf{isotropic Strichartz estimates} that we proved in the companion paper \cite{FCStratif1} (they are recalled in Section \ref{Section:iso}), and we prove in Section \ref{Section:aniso} the new \textbf{anisotropic estimates}, which are crucial to obtain convergence.}
\end{rem}

\section{Global existence and convergence of Strong solutions in the general case: proof of Theorem \ref{Th1}}
\label{PrvTh1}
The aim of this section is to prove the global existence of strong solutions when the Froude number $\ee$ is small enough, and the announced convergence rates in the general case (when we do not assume that $\nu=\nu'$).

\begin{rem}
 \sl{We emphasize that \emph{strong solutions} refers to strong solutions in the sense of Fujita and Kato, the existence and uniqueness of which is stated in Theorem \ref{TH0FK}. This theorem is proved in Section \ref{PreuveFK}.}
\end{rem}

\subsection{A priori estimates in the general case}
\label{Apriorigen}
Let us begin with the system satisfied by $\de=\De-\Wet$:
\begin{equation}
 \begin{cases}
  \d_t \de -L\de +\frac{1}{\ee} \mathbb{P} \mathcal{B} \de = \Sum_{i=1}^{11} F_i,\\
  {\de}_{|t=0}= (Id-\cPrR)\Uoeosc +(\UoeS^h-\tvo^h,0,0),
 \end{cases}
\label{de}
\end{equation}
where we define:
\begin{equation}
\begin{cases}
 F_1 \overset{def}{=}-\mathbb{P}(\de \cdot \n \de), \quad F_2 \overset{def}{=}-\mathbb{P}(\de \cdot \n \tv^h,0,0), \quad F_3 \overset{def}{=}-\mathbb{P}(\tv^h \cdot \n_h \de),\\
 F_4 \overset{def}{=}-\mathbb{P}(\de \cdot \nabla \Wet),\quad F_5 \overset{def}{=}-\mathbb{P}(\Wet \cdot \n \de),\quad F_6 \overset{def}{=}-\mathbb{P}(\tv^h \cdot \n_h \Wet),\\
 F_7 \overset{def}{=}-\mathbb{P}(\Wet \cdot \n \tv^h,0,0),\quad F_8 \overset{def}{=}-\mathbb{P}(\Wet \cdot \n \Wet),\\
 F_9 \overset{def}{=}-\mathbb{P}(0,0,0,\de^3 \cdot \d_3 \tThee),\quad F_{10} \overset{def}{=}-\mathbb{P}(0,0,0,\We^{T,3} \cdot \d_3 \tThee),\\
F_{11} \overset{def}{=} (Id-\cPrR) \tG.
\end{cases}
 \label{systde}
\end{equation}
Most of these terms will be estimated thanks to the following usual Sobolev product laws as in \cite{FCPAA, FCcompl, FCRF}.
\begin{prop}
 \sl{There exists a constant $C>0$ such that for any $s_1,s_2<\frac32$ with $s_1+s_2>0$ and any $u\in \dot{H}^{s_1}(\R^3)$, $v\in \dot{H}^{s_2}(\R^3)$, then $uv\in \dot{H}^{s_1+s_2-\frac32}(\R^3)$ and we have:
 $$
 \|uv\|_{\dot{H}^{s_1+s_2-\frac32}(\R^3)} \leq C \|u\|_{\dot{H}^{s_1}(\R^3)} \|v\|_{\dot{H}^{s_2}(\R^3)}.
 $$
 }
 \label{Prodlaws}
\end{prop}
As in \cite{FCStratif1} the terms involving a product with $\tThee$ will require special attention: we will need not only the following modified Sobolev product laws (that can be proved similarly as their bidimensional counterpart from \cite{CDGG} or \cite{IGTri} involving products with functions depending on $x_h$) but also anisotropic Strichartz estimates (that we prove in the appendix).
\begin{prop}
 \sl{There exists a constant $C>0$ such that for any $s_1,s_2<\frac12$ with $s_1+s_2>0$ and any $u\in \dot{H}^{s_1}(\R^3)$, $v\in \dot{H}^{s_2}(\R)$, then $uv\in \dot{H}^{s_1+s_2-\frac12}(\R^3)$ and we have:
 $$
 \|uv\|_{\dot{H}^{s_1+s_2-\frac12}(\R^3)} \leq C \|u\|_{\dot{H}^{s_1}(\R^3)} \|v\|_{\dot{H}^{s_2}(\R)}.
 $$
 }
 \label{prod1D3D}
\end{prop}

\begin{rem}
 \sl{As in \cite{FCPAA, FCcompl, FCRF}, we expect the previous external force terms to be small or absorbed:
 \begin{itemize}
  \item for $i\in\{1,2,3,9\}$, $F_i$ will be absorbed either by diffusion or thanks to Gronwall estimates,
  \item the first term in the initial data and $F_{11}$ are small thanks to frequency truncations,
  \item for other indices, $F_i$ will be split into a part absorbed thanks to diffusion, a part small thanks to dispersion, and in most of the cases a part absorbed thanks to Gronwall estimates.
 \end{itemize}
 }
\end{rem}

\subsubsection{Estimates in $\dot{H}^s$}

Computing the innerproduct in $\dot{H}^s$ ($s$ will be later fixed as $\frac12$ or $\frac12+\eta \delta$) of \eqref{de} with $\de$ we obtain (recall that $\nu_0=\min(\nu,\nu')$):
$$
\frac12 \frac{d}{dt} \|\de(t)\|_{\dot{H}^s}^2 +\nu_0 \|\n \de(t)\|_{\dot{H}^s}^2 \leq \Sum_{j=1}^{11} (F_j|\de)_{\dot{H}^s}.
$$
Using successively the classical Sobolev product laws for $(s_1,s_2)\in\{(\frac12,s),(s,\frac12),(1,s-\frac12)\}$, we obtain:
\begin{equation}
 |(F_1|\de)_{\dot{H}^s}|\leq \|\de \cdot \n \de\|_{\dot{H}^{s-1}} \|\de\|_{\dot{H}^{s+1}} \leq C \|\de\|_{\dot{H}^\frac12} \|\n\de\|_{\dot{H}^s}^2,
 \label{F1}
\end{equation}
\begin{multline}
 |(F_2|\de)_{\dot{H}^s}|\leq \|\de \cdot \n \tv^h\|_{\dot{H}^{s-1}} \|\de\|_{\dot{H}^{s+1}} \leq C \|\de\|_{\dot{H}^s} \|\n \tv^h\|_{\dot{H}^\frac12} \|\n\de\|_{\dot{H}^s}\\
 \leq \frac{\nu_0}{18}\|\n\de\|_{\dot{H}^s}^2+\frac{C}{\nu_0} \|\n \tv^h\|_{\dot{H}^\frac12}^2 \|\de\|_{\dot{H}^s}^2,
  \label{F2}
\end{multline}
\begin{multline}
 |(F_3|\de)_{\dot{H}^s}|\leq \|\tv^h \cdot \n_h \de\|_{\dot{H}^{s-1}} \|\de\|_{\dot{H}^{s+1}} \leq C \|\tv^h\|_{\dot{H}^1} \|\n_h \de\|_{\dot{H}^{s-\frac12}} \|\n\de\|_{\dot{H}^s}\\
 \leq C \|\tv^h\|_{\dot{H}^1} \|\n \de\|_{\dot{H}^{s-\frac12}} \|\n\de\|_{\dot{H}^s} \leq C (\|\tv^h\|_{\dot{H}^\frac12}^\frac12 \|\n \tv^h\|_{\dot{H}^\frac12}^\frac12)\|\de\|_{\dot{H}^{s+\frac12}} \|\n\de\|_{\dot{H}^s}\\
 \leq C \|\tv^h\|_{\dot{H}^\frac12}^\frac12 \|\n \tv^h\|_{\dot{H}^\frac12}^\frac12 \|\de\|_{\dot{H}^s}^\frac12 \|\n\de\|_{\dot{H}^s}^\frac32 \leq \frac{\nu_0}{18}\|\n\de\|_{\dot{H}^s}^2+\frac{C}{\nu_0^3} \|\tv^h\|_{\dot{H}^\frac12}^2 \|\n \tv^h\|_{\dot{H}^\frac12}^2 \|\de\|_{\dot{H}^s}^2.
  \label{F3}
\end{multline}
In the last estimate, we also used twice interpolation for Sobolev spaces, and the Young inequality for $(p,q)=(4,\frac43)$. We will estimate the next three terms (involving $F_4$, $F_5$ and $F_8$) reproducing what we did in \cite{FCRF} (see Section (2.2)), thanks to the Sobolev injections:
$$
  |(F_4|\de)_{\dot{H}^s}|\leq \|\de \cdot \nabla \Wet\|_{L^2} \|\de\|_{\dot{H}^{2s}} \leq C \|\de\|_{L^6} \|\nabla \Wet\|_{L^3} \|\de\|_{\dot{H}^{2s}} \leq C \|\de\|_{\dot{H}^1} \|\nabla \Wet\|_{L^3} \|\de\|_{\dot{H}^{2s}}.
$$
Thanks once more to interpolation ($1=(1-\aa)s+\aa (s+1)$ with $\aa=1-s$, and $2s=(1-\aa')s+\aa' (s+1)$ with $\aa'=s$, we will precise later that $s \in \{\frac12, \frac12+\eta \delta\}$), we obtain:
\begin{equation}
 |(F_4|\de)_{\dot{H}^s}|\leq C \|\de\|_{\dot{H}^s} \|\n \Wet\|_{L^3} \|\n\de\|_{\dot{H}^s}\\
 \leq \frac{\nu_0}{18}\|\n\de\|_{\dot{H}^s}^2+\frac{C}{\nu_0} \|\n \Wet\|_{L^3}^2 \|\de\|_{\dot{H}^s}^2.
  \label{F4}
\end{equation}
Similarly (but interpolating through $\frac32=(1-\aa)s+\aa (s+1)$ with $\aa=\frac32-s$ and $s \in \{\frac12, \frac12+\eta \delta\}$):
\begin{multline}
 |(F_5|\de)_{\dot{H}^s}| \leq \|\Wet \cdot \n \de\|_{L^2} \|\de\|_{\dot{H}^{2s}} \leq C \|\Wet\|_{L^6} \|\de\|_{\dot{H}^\frac32} \|\de\|_{\dot{H}^{2s}}\\
 \leq C \|\Wet\|_{L^6} \|\de\|_{\dot{H}^s}^\frac12 \|\n \de\|_{\dot{H}^s}^\frac32 \leq \frac{\nu_0}{18}\|\n\de\|_{\dot{H}^s}^2+\frac{C}{\nu_0^3} \|\Wet\|_{L^6}^4 \|\de\|_{\dot{H}^s}^2.
  \label{F5}
\end{multline}
Thanks to the Young inequality with $(\frac2{s},\frac2{1-s}, 2)$,
\begin{multline}
 |(F_8|\de)_{\dot{H}^s}| \leq \|\Wet \cdot \n \Wet\|_{L^2} \|\de\|_{\dot{H}^{2s}} \leq C \|\de\|_{\dot{H}^{s+1}}^s \left(\|\Wet\|_{L^6} \|\de\|_{\dot{H}^s}^{1-s}\right) \|\n\Wet\|_{L^3}\\
 \leq \frac{\nu_0}{18}\|\n\de\|_{\dot{H}^s}^2+\frac{C}{\nu_0^{\frac{s}{1-s}}} \|\Wet\|_{L^6}^\frac2{1-s} \|\de\|_{\dot{H}^s}^2 + C\|\n \Wet\|_{L^3}^2.
  \label{F8}
\end{multline}
The next terms can be estimated in an improved way compared to the methods from \cite{FCRF} as, on one hand, $\tv^h$ and $\tG$ are more
regular, and, on the other hand, we can "pay" on $\tThee$ what we need for $\Uoeosc$. More precisely, thanks to Theorem \ref{ThSNS}, by interpolation, we have that for all $s\in[0,\frac12+\delta]$ and for all $q\geq 2$,
$$
\|\tv^h\|_{L^q \dot{H}^{s+\frac2{q}}} \leq C_{\delta, \nu} \max (1,\|\tvo^h\|_{H^{\frac12+\delta}})^{1+\frac1{2\delta}},
$$
and as $\frac12+\delta+\frac2{q}=\frac32 \Leftrightarrow q=\frac2{1-\delta}$, we obtain that:
\begin{equation}
 \|\tv^h\|_{L^2 \dot{H}^{\frac32} \cap L^\frac2{1-\delta} \dot{H}^{\frac32}} \leq C_{\delta, \nu} \max (1,\|\tvo^h\|_{H^{\frac12+\delta}})^{1+\frac1{2\delta}}.
 \label{estimplage}
\end{equation}
Now, thanks to the Sobolev injection $\dot{H}^1(\R^3) \hookrightarrow L^6(\R^3)$, interpolation, and using once more the Young inequality with $(\frac2{s},\frac2{1-s}, 2)$:
\begin{multline}
 |(F_6|\de)_{\dot{H}^s}| \leq \|\tv^h \cdot \n_h \Wet\|_{L^2} \|\de\|_{\dot{H}^{2s}} \leq C \|\tv^h\|_{\dot{H}^\frac12}^\frac12 \|\n\tv^h\|_{\dot{H}^\frac12}^\frac12 \|\n\Wet\|_{L^3} \|\de\|_{\dot{H}^s}^{1-s} \|\de\|_{\dot{H}^{s+1}}^s\\
\leq C \|\de\|_{\dot{H}^{s+1}}^s \left(\|\n \tv^h\|_{\dot{H}^\frac12}^\frac12 \|\de\|_{\dot{H}^s}^{1-s}\right) \left(\|\tv^h\|_{\dot{H}^\frac12}^\frac12 \|\n\Wet\|_{L^3}\right)\\
 \leq \frac{\nu_0}{18}\|\n\de\|_{\dot{H}^s}^2+\frac{C}{\nu_0^{\frac{s}{1-s}}} \|\n \tv^h\|_{\dot{H}^\frac12}^\frac1{1-s} \|\de\|_{\dot{H}^s}^2 + C \|\tv^h\|_{\dot{H}\frac12} \|\n \Wet\|_{L^3}^2.
  \label{F6}
\end{multline}
\begin{rem}
 \sl{Observe that when $s=\frac12+\eta \delta$ with $\eta\in]0,\frac12]$, then $\frac1{1-s}\in[2,\frac2{1-\delta}]$ and we can use \eqref{estimplage}.}
\end{rem}
Similarly,
\begin{multline}
 |(F_7|\de)_{\dot{H}^s}| \leq \|\Wet \cdot \n \tv^h\|_{L^2} \|\de\|_{\dot{H}^{2s}} \leq C \|\Wet\|_{L^6} \|\n \tv^h\|_{\dot{H}^\frac12} \|\de\|_{\dot{H}^s}^{1-s} \|\de\|_{\dot{H}^{s+1}}^s\\
\leq C \|\de\|_{\dot{H}^{s+1}}^s \left(\|\n \tv^h\|_{\dot{H}^\frac12}^\frac12 \|\de\|_{\dot{H}^s}^{1-s}\right) \left(\|\n \tv^h\|_{\dot{H}^\frac12}^\frac12 \|\Wet\|_{L^6}\right)\\
 \leq \frac{\nu_0}{18}\|\n\de\|_{\dot{H}^s}^2+\frac{C}{\nu_0^{\frac{s}{1-s}}} \|\n \tv^h\|_{\dot{H}^\frac12}^\frac1{1-s} \|\de\|_{\dot{H}^s}^2 + C \|\n \tv^h\|_{\dot{H}^\frac12} \|\Wet\|_{L^6}^2.
  \label{F7}
\end{multline}
We easily obtain that:
\begin{multline}
 |(F_{11}|\de)_{\dot{H}^s}| \leq \|(Id-\cPrR) \tG\|_{\dot{H}^s} \|\de\|_{\dot{H}^s}\\
 \leq \frac12 \|(Id-\cPrR) \tG\|_{\dot{H}^s} +\frac12 \|(Id-\cPrR) \tG\|_{\dot{H}^s} \|\de\|_{\dot{H}^s}^2,
 \label{F11}
\end{multline}
and we are left with the new terms involving $\tThee(x_3)$. Let us begin with $F_{10}$: introducing the following anisotropic norms (with the classical adaptations for infinite exponents),
\begin{equation}
 \|f\|_{L_{v,h}^{p, q}} \overset{def}{=} \big\| \|f\|_{L^q(\R_h^2)} \big\|_{L^p(\R_v)} =\left(\int_\R \Big(\int_{\R^2} |f(x_h,x_3)|^q dx_h\Big)^\frac{p}{q} dx_3\right)^\frac1{p},
\label{Anisospace}
 \end{equation}
we have for some $\aa\in[0,1]$ (to be precised below):
\begin{multline}
 |(F_{10}|\de)_{\dot{H}^s}| \leq \|\We^{T,3} \cdot \d_3 \tThee\|_{L^2} \|\de\|_{\dot{H}^{2s}} \leq C \|\Wet\|_{L_{v,h}^{\infty, 2}} \|\d_3 \tThee\|_{L^2(\R)} \|\de\|_{\dot{H}^s}^{1-s} \|\de\|_{\dot{H}^{s+1}}^s\\
 \leq C \|\de\|_{\dot{H}^{s+1}}^s \left(\|\tThee\|_{\dot{H}^1(\R)}^{1-\aa} \|\de\|_{\dot{H}^s}^{1-s}\right) \left(\|\tThee\|_{\dot{H}^1(\R)}^{\aa} \|\Wet\|_{L_{v,h}^{\infty, 2}}\right)\\
 \leq \frac{\nu_0}{18}\|\n\de\|_{\dot{H}^s}^2+\frac{C}{\nu_0^{\frac{s}{1-s}}} \|\tThee\|_{\dot{H}^1(\R)}^{2\frac{1-\aa}{1-s}} \|\de\|_{\dot{H}^s}^2 + C \|\tThee\|_{\dot{H}^1(\R)}^{2\aa} \|\Wet\|_{L_{v,h}^{\infty, 2}}^2.
\end{multline}
\begin{rem}
 \sl{We emphasize that in \cite{CDGG}, the anisotropic norms were of the form $\|f\|_{L_{h,v}^{q, p}}$.}
 \label{Anisorem}
\end{rem}
As we aim for the best possible convergence rate, we will use below the Strichartz estimates for $\|\Wet\|_{L^8 L_{v,h}^{\infty, 2}}$ which will require us to estimate (thanks to \eqref{estimThetaBspr}):
\begin{equation}
\begin{cases}\vspace{0.2cm}
\|\tThee\|_{L^{2\frac{1-\aa}{1-s}} \dot{H}^1(\R)} \leq \|\tThee\|_{L^{2\frac{1-\aa}{1-s}} \dot{B}_{2,1}^1(\R)} \leq \frac{C}{(\nu')^{\frac{1-s}{2(1-\aa)}}} \|\tTheeo\|_{\dot{B}_{2,1}^\frac{s-\aa}{1-\aa}},\\
 \|\tThee\|_{L^{\frac{8\aa}3} \dot{H}^1(\R)} \leq \|\tThee\|_{L^{\frac{8\aa}3} \dot{B}_{2,1}^1(\R)} \leq \frac{C}{(\nu')^{\frac{3}{8\aa}}} \|\tTheeo\|_{\dot{B}_{2,1}^{1-\frac3{4\aa}}}.
\end{cases}
\label{HypThee}
\end{equation}
So, in addition to the assumption $\tTheeo\in \dot{B}_{2,1}^{-\frac12}$, we will need:
$$
\tTheeo\in \dot{B}_{2,1}^\frac{s-\aa}{1-\aa} \cap \dot{B}_{2,1}^{1-\frac3{4\aa}}.
$$
As we wish to make the least assumptions possible on $\tThee$, we will simply choose $\aa$ so that the previous Besov exponents are equal, that is $\aa=\frac3{7-4s}$. With this choice the additional assumptions on $\tTheeo$ reduce to $\tTheeo\in \dot{B}_{2,1}^{s-\frac34}$ as:
\begin{equation}
 \|\tThee\|_{L^{\frac8{7-4s}} \dot{H}^1(\R)} \leq \frac{C}{(\nu')^{\frac{7-4s}8}} \|\tTheeo\|_{\dot{B}_{2,1}^{s-\frac34}},
\end{equation}
and
\begin{equation}
 |(F_{10}|\de)_{\dot{H}^s}| \leq \frac{\nu_0}{18}\|\n\de\|_{\dot{H}^s}^2+\frac{C}{\nu_0^{\frac{s}{1-s}}} \|\tThee\|_{\dot{H}^1(\R)}^{\frac8{7-4s}} \|\de\|_{\dot{H}^s}^2 + C \|\tThee\|_{\dot{H}^1(\R)}^{\frac6{7-4s}} \|\Wet\|_{L_{v,h}^{\infty, 2}}^2.
 \label{F10}
\end{equation}
The last term is bounded thanks to the modified Sobolev product laws from Proposition \ref{prod1D3D}. Introducing, for some $\beta\in]0,\frac12[$ (which can be considered as small as we need, and will be precised in Remark \ref{Rembeta}) $(s_1,s_2)= (\frac12-\beta, s-1+\beta)$, and roughly bounding the following homogeneous Sobolev norm with an inhomogeneous one according to $\|\de\|_{\dot{H}^{\frac12-\beta}}\leq \|\de\|_{H^{\frac12-\beta}}\leq \|\de\|_{H^s}$, we obtain that when $s>\frac12$:
\begin{multline}
 |(F_{9}|\de)_{\dot{H}^s}| \leq \|\de^3 \cdot \d_3 \tThee\|_{\dot{H}^{s-1}} \|\de\|_{\dot{H}^{s+1}}\\
 \leq C \|\de\|_{\dot{H}^{\frac12-\beta}} \|\d_3 \tThee\|_{\dot{H}^{s-1+\beta}(\R)} \|\de\|_{\dot{H}^{s+1}} \leq C \|\de\|_{H^s} \|\d_3 \tThee\|_{\dot{H}^{s-1+\beta}(\R)} \|\de\|_{H^{s+1}}\\
 \leq \frac{\nu_0}4\|\n\de\|_{H^s}^2+\frac{C}{\nu_0} \|\tThee\|_{\dot{H}^{s+\beta}(\R)}^2 \|\de\|_{H^s}^2.
  \label{F9}
\end{multline}
When $s=\frac12$, we introduce $\bb_1,\bb_2>0$ with $\bb=\bb_1+\bb_2<1$ and use Proposition \ref{prod1D3D} with $(s_1,s_2)=(\frac12-\bb_2, -\frac12+\bb_1+\bb_2)$:
\begin{multline}
 |(F_{9}|\de)_{\dot{H}^\frac12}| \leq \|\de^3 \cdot \d_3 \tThee\|_{\dot{H}^{\frac12+\bb_1}} \|\de\|_{\dot{H}^{\frac32-\bb_1}}\\
 \leq C \|\de\|_{\dot{H}^{\frac12-\bb_2}} \|\d_3 \tThee\|_{\dot{H}^{-\frac12+\bb_1+\bb_2}(\R)} \|\de\|_{\dot{H}^{\frac32-\bb_1}} \leq C \|\de\|_{H^\frac12} \|\d_3 \tThee\|_{\dot{H}^{\frac12+\bb_1+\bb_2}(\R)} \|\de\|_{H^{\frac32}}\\
 \leq \frac{\nu_0}4\|\n\de\|_{H^\frac12}^2+\frac{C}{\nu_0} \|\tThee\|_{\dot{H}^{\frac12+\beta}(\R)}^2 \|\de\|_{H^\frac12}^2.
  \label{F9bis}
\end{multline}
\begin{rem}
 \sl{As in \cite{FCStratif1}, this term requires special attention and, in the present article, is dealt thanks to the low frequency assumptions.}
\end{rem}
Collecting \eqref{F1} to \eqref{F11}, with \eqref{F10} and \eqref{F9}, we obtain that for all $t<T_\ee^*$,
\begin{multline}
 \frac12 \frac{d}{dt} \|\de\|_{\dot{H}^s}^2 +\frac{\nu_0}2 \|\n \de\|_{\dot{H}^s}^2 \leq C \|\de\|_{\dot{H}^\frac12} \|\n\de\|_{\dot{H}^s}^2 +\frac{C}{\nu_0} \|\de\|_{H^s}^2 \|\tThee\|_{\dot{H}^{s+\beta}(\R)}^2 +\frac{\nu_0}4\|\n\de\|_{H^s}^2\\
 +\frac{C}{\nu_0} \|\de\|_{\dot{H}^s}^2 \Bigg( (1+\frac1{\nu_0^2} \|\tv^h\|_{\dot{H}^\frac12}^2) \|\n \tv^h\|_{\dot{H}^\frac12}^2 +\|\n \Wet\|_{L^3}^2 +\frac1{\nu_0^2} \|\Wet\|_{L^6}^4\\
 +\nu_0 \|\tG\|_{\dot{H}^s} +\frac1{\nu_0^{\frac{2s-1}{1-s}}} (\|\n \tv^h\|_{\dot{H}^\frac12}^\frac1{1-s} +\|\Wet\|_{L^6}^\frac2{1-s} +\|\tThee\|_{\dot{H}^1(\R)}^{\frac8{7-4s}})\Bigg)\\
 + C\Bigg[ (1+\|\tv^h\|_{\dot{H}\frac12}) \|\n \Wet\|_{L^3}^2 +\|\n \tv^h\|_{\dot{H}^\frac12} \|\Wet\|_{L^6}^2\\
 +\|\tThee\|_{\dot{H}^1(\R)}^{\frac6{7-4s}} \|\Wet\|_{L_{v,h}^{\infty, 2}}^2 +\frac12 \|(Id-\cPrR) \tG\|_{\dot{H}^s}\Bigg].
 \label{estimdHs}
\end{multline}

\subsubsection{Estimates in $L^2$}

As explained in the previous section, dealing with $F_9$ required additional low frequency assumptions, therefore we also need energy estimates in $L^2$. Computing the innerproduct in $L^2$ of \eqref{de} with $\de$, we similarly obtain that:
$$
\frac12 \frac{d}{dt} \|\de(t)\|_{L^2}^2 +\nu_0 \|\n \de(t)\|_{L^2}^2 \leq \Sum_{j=1}^{11} (F_j|\de)_{L^2}.
$$
The complete or horizontal divergence-free conditions imply that
$$
(F_1|\de)_{L^2} =(F_3|\de)_{L^2} =(F_5|\de)_{L^2} =0
$$
The next three terms are dealt with similar arguments as previously:
\begin{multline}
|(F_2|\de)_{L^2}|\leq \|\de \cdot \n \tv^h\|_{L^2} \|\de^h\|_{L^2} \leq
 C\|\de\|_{L^6} \|\n \tv^h\|_{L^3} \|\de\|_{L^2}\\
 \leq C\|\de\|_{\dot{H}^1} \|\n \tv^h\|_{\dot{H}^\frac12} \|\de\|_{L^2} \leq \frac{\nu_0}8 \|\n \de\|_{L^2}^2 +\frac{C}{\nu_0} \|\n \tv^h\|_{\dot{H}^\frac12}^2 \|\de^h\|_{L^2}^2,
\label{F2b}
\end{multline}
\begin{equation}
\begin{cases}
\vspace{0.2cm}
 |(F_4|\de)_{L^2}|\leq \|\de \cdot \n \Wet\|_{L^2} \|\de\|_{L^2} \leq \frac{\nu_0}8 \|\n \de\|_{L^2}^2 +\frac{C}{\nu_0} \|\n \Wet\|_{L^3}^2 \|\de\|_{L^2}^2.\\
 |(F_6|\de)_{L^2}|\leq \|\tv^h \cdot \n_h \Wet\|_{L^2} \|\de\|_{L^2} \leq C \|\n \Wet\|_{L^3}^2 +C \|\tv^h\|_{\dot{H}^1}^2 \|\de\|_{L^2}^2.
\end{cases}
\label{F46b}
\end{equation}
The next term is estimated differently (in order to minimize the assumptions on $\Uoeosc$), thanks to the Young inequality with indices $(4,4,2)$:
\begin{multline}
|(F_7|\de)_{L^2}| \leq \|\Wet\|_{L^6} \|\n \tv^h\|_{L^2} \|\de^h\|_{L^3} \leq C \|\Wet\|_{L^6} \|\n \tv^h\|_{L^2} \|\de\|_{\dot{H}\frac12}\\
\leq C \|\Wet\|_{L^6} \|\n \tv^h\|_{L^2} \|\de\|_{L^2}^\frac12 \|\n \de\|_{L^2}^\frac12 \leq C \|\n \de\|_{L^2}^\frac12 \big(\|\de\|_{L^2}^\frac12 \|\n \tv^h\|_{L^2}^\frac12 \big) \big(\|\n \tv^h\|_{L^2}^\frac12 \|\Wet\|_{L^6} \big)\\
\leq \frac{\nu_0}8 \|\n \de\|_{L^2}^2 +\frac{C}{\nu_0} \|\de\|_{L^2}^2 \|\n \tv^h\|_{L^2}^2 +\frac12 \|\n \tv^h\|_{L^2} \|\Wet\|_{L^6}^2.
\label{F7b}
\end{multline}
The following term also requires special attention, because if we use the same arguments as for $F_{2,4,6}$, we end-up with $\|\Wet\|_{L^2 L^6}$ which, in the case $\nu=\nu'$, would require additionnal assumptions on $\|\Uoeosc\|_{\dot{H}^\frac13}$. To avoid this, for $r_1,r_2>2$ such that $\frac1{r_1}+\frac1{r_2}=\frac12$, let us write:
\begin{equation}
|(F_8|\de)_{L^2}|\leq C \|\Wet\|_{L^{r_1}} \|\n \Wet\|_{L^{r_2}} \|\de\|_{L^2} \leq C \|\Wet\|_{L^{r_1}}^2 +C \|\n \Wet\|_{L^{r_2}}^2 \|\de\|_{L^2}^2.
\end{equation}
As we will see later, estimating $\|\Wet\|_{L^2 L^{r_1}}$ and $\|\n \Wet\|_{L^2 L^{r_2}}$ in the case $\nu=\nu'$ will make us deal with the norm of $\Uoeosc$ in the spaces $\dot{H}^{\sigma_1}$ and $\dot{H}^{\sigma_2}$ where:
$$
\sigma_1= \frac12-\frac3{r_1}+\theta_1(\frac12-\frac1{r_1})\quad \mbox{and} \quad \sigma_2= \frac32-\frac3{r_2}+\theta_2(\frac12-\frac1{r_2}),
$$
with $\theta_{1,2}\in]0,1]$. Using $\frac1{r_1}+\frac1{r_2}=\frac12$, we have $\sigma_2=\frac{3+\theta_2}{r_1}$, and the fewest assumptions are made when we choose:
$$
\sigma_1=\sigma_2=\frac{(3+\theta_2)(1+\theta_1)}{2(6+\theta_1+\theta_2)}.
$$
This function (of $(\theta_1,\theta_2)\in[0,1]^2$) reaches its maximum $\frac12$ when $\theta_1=\theta_2=1$ which corresponds to $(r_1,r_2)=(8, \frac83)$ so that we finally get the estimates:
\begin{equation}
|(F_8|\de)_{L^2}|\leq C \|\Wet\|_{L^8}^2 +C \|\n \Wet\|_{L^\frac83}^2 \|\de\|_{L^2}^2.
\label{F8b}
\end{equation}
Obviously
\begin{equation}
 |(F_{11}|\de)|_{L^2} \leq \frac12 \|(Id-\cPrR) \tG\|_{L^2} +\frac12 \|\tG\|_{L^2} \|\de\|_{L^2}^2,
 \label{F11b}
\end{equation}
and we are left with the last two terms, involving $\tThee$. The first one is bounded like in the proof of Proposition 3.1 in \cite{FCStratif1} using the Minkowski and Young estimates (twice for $(\frac43,4)$), the 1D-Sobolev injection $\dot{H}^\frac14 (\R) \hookrightarrow L^4(\R)$  and interpolation:
\begin{multline}
|(F_9|\de)_{L^2}|\leq \int_{\R^2} \left(\int_{\R} |\de(x_h,x_3)|^2 |\d_3 \tThee (x_3)| dx_3 \right) dx_h \leq C\|\d_3 \tThee\|_{L^2(\R)} \int_{\R^2} \|\De(x_h,\cdot)\|_{\dot{H}^\frac14(\R)}^2  dx_h\\
 \leq C\|\tThee\|_{\dot{H}^1(\R)} \int_{\R^2} \|\de(x_h,\cdot)\|_{L^2(\R)}^\frac32 \|\de(x_h,\cdot)\|_{\dot{H}^1(\R)}^\frac12 dx_h \leq C\|\tThee\|_{\dot{H}^1(\R)} \|\de\|_{L^2(\R^3)}^\frac32 \|\d_3 \de\|_{L^2(\R^3)}^\frac12\\
 \leq C\|\tThee\|_{\dot{H}^1(\R)} \|\de\|_{L^2}^\frac32 \|\n \de\|_{L^2}^\frac12 \leq \frac{\nu_0}8 \|\n \de\|_{L^2}^2 +\frac{C}{\nu_0^\frac13} \|\tThee\|_{\dot{H}^1(\R)}^\frac43 \|\de\|_{L^2}^2.
 \label{F9b}
\end{multline}
As for the $\dot{H}^s$-estimates, the last term will require adjustment in order to minimize the assumptions on $\tThee$. For some $\aa>0$ (to be precised later):
\begin{multline}
 |(F_{10}|\de)_{L^2}| \leq \|\We^{T,3} \cdot \d_3 \tThee\|_{L^2} \|\de\|_{L^2} \leq C \|\Wet\|_{L_{v,h}^{\infty, 2}} \|\d_3 \tThee\|_{L^2(\R)} \|\de\|_{L^2}\\
 \leq C \Big(\|\tThee\|_{\dot{H}^1(\R)}^{1-\aa} \|\de\|_{L^2} \Big) \Big(\|\tThee\|_{\dot{H}^1(\R)}^\aa \|\Wet\|_{L_{v,h}^{\infty, 2}} \Big)\\
 \leq C \|\tThee\|_{\dot{H}^1(\R)}^{2(1-\aa)} \|\de\|_{L^2}^2 +C \|\tThee\|_{\dot{H}^1(\R)}^{2\aa} \|\Wet\|_{L_{v,h}^{\infty, 2}}^2.
\end{multline}
As we aim for the best convergence rate, we will bound $\|\Wet\|_{L^8 L_{v,h}^{\infty, 2}}$ which, similarly to \eqref{HypThee}, will require a control on:
\begin{equation}
\begin{cases}\vspace{0.2cm}
\|\tThee\|_{L^{2(1-\aa)} \dot{H}^1(\R)} \leq \frac{C}{(\nu')^{\frac1{2(1-\aa)}}} \|\tTheeo\|_{\dot{B}_{2,1}^{-\frac{\aa}{1-\aa}}},\\
 \|\tThee\|_{L^{\frac{8\aa}3} \dot{H}^1(\R)} \leq \frac{C}{(\nu')^{\frac{3}{8\aa}}} \|\tTheeo\|_{\dot{B}_{2,1}^{1-\frac3{4\aa}}}.
\end{cases}
\label{HypThee2}
\end{equation}
The best choice is when both regularity indices are equal, that is when $\aa=\frac37$, so that we finally obtain:
\begin{equation}
 |(F_{10}|\de)_{L^2}| \leq C \|\tThee\|_{\dot{H}^1(\R)}^{\frac87} \|\de\|_{L^2}^2 +C \|\tThee\|_{\dot{H}^1(\R)}^{\frac67} \|\Wet\|_{L_{v,h}^{\infty, 2}}^2.
 \label{F10b}
\end{equation}
Gathering \eqref{F2b} to \eqref{F7b}, \eqref{F8b} to \eqref{F9b} and \eqref{F10b}, we obtain that for all $t<T_\ee^*$:
\begin{multline}
 \frac12 \frac{d}{dt} \|\de\|_{L^2}^2 +\frac{\nu_0}2 \|\n \de\|_{L^2}^2 \leq C \|\de^h\|_{L^2}^2 \Bigg( \frac1{\nu_0} (\|\n \tv^h\|_{\dot{H}^\frac12}^2 +\|\n \Wet\|_{L^3}^2)\\
 +(1+\frac1{\nu_0}) \|\n \tv^h\|_{L^2}^2 +\|\n \Wet\|_{L^\frac83}^2 +\|\tG\|_{L^2} +\frac1{\nu_0^\frac13} \|\tThee\|_{\dot{H}^1(\R)}^\frac43 +\|\tThee\|_{\dot{H}^1(\R)}^{\frac87} \Bigg)\\
 +C \Bigg[\|\n \Wet\|_{L^3}^2 +\|\n \tv^h\|_{L^2} \|\Wet\|_{L^6}^2 +\|\Wet\|_{L^8}^2 +\|\tThee\|_{\dot{H}^1(\R)}^{\frac67} \|\Wet\|_{L_{v,h}^{\infty, 2}}^2\\
  +\|(Id-\cPrR) \tG\|_{L^2} \Bigg]
 \label{estimdL2}
\end{multline}

\subsubsection{Estimates in $H^s$}

We recall that for any $s\in]0,1]$ and any $f\in H^s(\R^3)$,
\begin{equation}
 \frac12 (\|f\|_{L^2}^2+\|f\|_{\dot{H}^s}^2)\leq 2^{s-1}(\|f\|_{L^2}^2+\|f\|_{\dot{H}^s}^2)\leq \|f\|_{H^s}^2 \leq \|f\|_{L^2}^2+\|f\|_{\dot{H}^s}^2.
\label{Sobineq}
\end{equation}
Let us introduce
\begin{equation}
 T_{\ee,2}\overset{def}{=} \sup \Big\{t\in [0,T_\ee^*[/ \; \forall t'\in[0,t],\; \|\de(t')\|_{\dot{H}^\frac12} \leq \frac{\nu_0}{4C}.\Big\}
 \label{DefTe}
\end{equation}
As $\de(0)$ goes to zero when $\ee \rightarrow 0$ (we refer to \eqref{estiminit2} for details), we are sure that $T_{\ee,2}>0$ if $\ee>0$ is small enough so that, gathering \eqref{estimdHs} and \eqref{estimdL2} we obtain that for all $\beta\in]0,1[$ and all $t\leq T_{\ee,2}$,
\begin{equation}
 \frac{d}{dt} (\|\de\|_{L^2}^2 +\|\de\|_{\dot{H}^s}^2) +\frac{\nu_0}2 (\|\n \de\|_{L^2}^2 +\|\n \de\|_{\dot{H}^s}^2)\\
\leq C_{\nu,\nu',s} \left(\|\de\|_{H^s}^2 K(t) + J(t)\right),
\label{estimdHsfinale}
\end{equation}
where
\begin{multline}
 K(t)= (1+\|\tv^h\|_{\dot{H}^\frac12}^2) \|\n \tv^h\|_{\dot{H}^\frac12}^2 +\|\n \tv^h\|_{L^2}^2 +\|\n \tv^h\|_{\dot{H}^\frac12}^\frac1{1-s} +\|\tG\|_{L^2} +\|\tG\|_{\dot{H}^s} +\|\n \Wet\|_{L^3}^2\\
  +\|\Wet\|_{L^6}^4 +\|\Wet\|_{L^6}^\frac2{1-s} +\|\n \Wet\|_{L^\frac83}^2 +\|\tThee\|_{\dot{H}^1(\R)}^{\frac8{7-4s}} +\|\tThee\|_{\dot{H}^{s+\beta}(\R)}^2 +\|\tThee\|_{\dot{H}^1(\R)}^\frac43 +\|\tThee\|_{\dot{H}^1(\R)}^{\frac87}
\end{multline}
and
\begin{multline}
 J(t)= (1+\|\tv^h\|_{\dot{H}\frac12}) \|\n \Wet\|_{L^3}^2 +(\|\n \tv^h\|_{\dot{H}^\frac12} +\|\n \tv^h\|_{L^2}) \|\Wet\|_{L^6}^2 +\|\Wet\|_{L^8}^2\\
 +(\|\tThee\|_{\dot{H}^1(\R)}^{\frac6{7-4s}} +\|\tThee\|_{\dot{H}^1(\R)}^{\frac67}) \|\Wet\|_{L_{v,h}^{\infty, 2}}^2 +\|(Id-\cPrR) \tG\|_{L^2} +\|(Id-\cPrR) \tG\|_{\dot{H}^s}.
 \end{multline}
Thanks to the Gronwall lemma, and using once more \eqref{Sobineq}, we obtain that for all $t\leq T_{\ee,2}$,
\begin{multline}
 \|\de(t)\|_{H^s}^2+\frac{\nu_0}2 \int_0^t \|\n \de(t')\|_{H^s}^2 dt'\\
 \leq \left( \|\de(0)\|_{L^2}^2 +\|\de(0)\|_{\dot{H}^s}^2 + C_{\nu,\nu',s}\int_0^t J(t')dt'\right) e^{C_{\nu,\nu',s} \int_0^t K(\tau)d\tau}.
\end{multline}
Using Theorems \ref{ThHeat} and \ref{ThSNS} together with the assumptions on the initial data, we can bound $J$ and $K$ as follows:
\begin{multline}
 \int_0^t K(t')dt' \leq (1+\|\tv^h\|_{L_t^\infty \dot{H}^\frac12}^2) \|\n \tv^h\|_{L_t^2 \dot{H}^\frac12}^2 +\|\n \tv^h\|_{L_t^2 L^2}^2 +\|\n \tv^h\|_{L_t^\frac1{1-s}\dot{H}^\frac12}^\frac1{1-s}\\
 +\|\tG\|_{L_t^1 L^2} +\|\tG\|_{L_t^1 \dot{H}^s} +\|\n \Wet\|_{L_t^2 L^3}^2 +\|\Wet\|_{L_t^4 L^6}^4 +\|\Wet\|_{L_t^\frac2{1-s} L^6}^\frac2{1-s} +\|\n \Wet\|_{L_t^2 L^\frac83}^2\\
 +\|\tThee\|_{L_t^{\frac8{7-4s}} \dot{H}^1(\R)}^{\frac8{7-4s}} +\|\tThee\|_{L_t^2 \dot{H}^{s+\beta}(\R)}^2 +\|\tThee\|_{L_t^\frac43 \dot{H}^1(\R)}^\frac43 +\|\tThee\|_{L_t^\frac87 \dot{H}^1(\R)}^{\frac87}\\
 \leq (2+C_{\delta, \nu} \Co^{2+\frac1{\delta}})C_{\delta, \nu} \Co^{2+\frac1{\delta}} + (C_{\delta, \nu} \Co^{2+\frac1{\delta}})^\frac1{1-s}\\
 +C_{\nu',s}\left(\|\tTheeo\|_{\dot{B}_{2,1}^{s-\frac34}(\R)}^{\frac8{7-4s}} +\|\tTheeo\|_{\dot{B}_{2,1}^{s+\beta-1}(\R)}^2 +\|\tTheeo\|_{\dot{B}_{2,1}^{-\frac12}(\R)}^\frac43 +\|\tTheeo\|_{\dot{B}_{2,1}^{-\frac34}(\R)}^{\frac87}\right)\\
 +\|\n \Wet\|_{L_t^2 L^3}^2 +\|\Wet\|_{L_t^4 L^6}^4 +\|\Wet\|_{L_t^\frac2{1-s} L^6}^\frac2{1-s} +\|\n \Wet\|_{L_t^2 L^\frac83}^2\\
 \leq \Do +\|\n \Wet\|_{L_t^2 L^3}^2 +\|\Wet\|_{L_t^4 L^6}^4 +\|\Wet\|_{L_t^\frac2{1-s} L^6}^\frac2{1-s} +\|\n \Wet\|_{L_t^2 L^\frac83}^2,
\end{multline}
where $\Do=\Do(\nu,\nu',\Co, \delta,s)$. Similarly,
\begin{multline}
 \int_0^t J(t')dt'\leq (1+\|\tv^h\|_{L_t^\infty \dot{H}\frac12}) \|\n \Wet\|_{L_t^2 L^3}^2 +(\|\n \tv^h\|_{L_t^2 \dot{H}^\frac12} +\|\n \tv^h\|_{L_t^2 L^2}) \|\Wet\|_{L_t^4 L^6}^2\\
 +\|\Wet\|_{L_t^2 L^8}^2 +\|(Id-\cPrR) \tG\|_{L_t^1 L^2} +\|(Id-\cPrR) \tG\|_{L_t^1 \dot{H}^s}\\
 +(\|\tThee\|_{L_t^{\frac8{7-4s}} \dot{H}^1(\R)}^{\frac6{7-4s}} +\|\tThee\|_{L_t^{\frac87} \dot{H}^1(\R)}^{\frac67}) \|\Wet\|_{L_t^8 L_{v,h}^{\infty, 2}}^2\\
 \leq \|(Id-\cPrR) \tG\|_{L_t^1 L^2} +\|(Id-\cPrR) \tG\|_{L_t^1 \dot{H}^s}\\
 +\Do \left(\|\n \Wet\|_{L_t^2 L^3}^2 +\|\Wet\|_{L_t^4 L^6}^2 +\|\Wet\|_{L_t^2 L^8}^2 +\|\Wet\|_{L_t^8 L_{v,h}^{\infty, 2}}^2\right).
 \end{multline}
\begin{rem}
 \sl{Thanks to interpolation, the fact that $\|\tTheeo\|_{\dot{B}_{2,1}^{-\frac34}(\R) \cap \dot{B}_{2,1}^{-\frac14+\delta}(\R)} \leq \Co$ allowed us to properly bound every norm involving $\tThee$, including the norm $\|\tTheeo\|_{\dot{B}_{2,1}^{s+\beta-1}(\R)}^2$  for any $\beta>0$ as small as we need ($s+\beta-1\leq s-\frac34$ when $\beta\leq \frac14$): we simply choose $\beta\in]0,\frac14]$.
 }
 \label{Rembeta}
\end{rem}
This leads to the following estimates (we recall that we will choose $s=\frac12+\eta \delta$): there exists $\Do=\Do(\nu,\nu',\Co, \delta,s)$ such that for all $t\leq T_{\ee,2}$,
\begin{multline}
 \|\de(t)\|_{H^s}^2+\frac{\nu_0}2 \int_0^t \|\n \de(t')\|_{H^s}^2 dt' \leq \Bigg[\|(Id-\cPrR) \Uoeosc\|_{L^2 \cap \dot{H}^s}^2 +\|(Id-\cPrR) \tG\|_{L_t^1 (L^2 \cap \dot{H}^s)}\\
 +\|\UoeS^h-\tvo^h\|_{L^2 \cap \dot{H}^s}^2 +\Do \left(\|\n \Wet\|_{L_t^2 L^3}^2 +\|\Wet\|_{L_t^4 L^6}^2 +\|\Wet\|_{L_t^2 L^8}^2 +\|\Wet\|_{L_t^8 L_{v,h}^{\infty, 2}}^2\right)\Bigg]\\
 \times \exp \Bigg\{\Do \Big(1 +\|\n \Wet\|_{L_t^2 L^3}^2 +\|\Wet\|_{L_t^4 L^6}^4 +\|\Wet\|_{L_t^\frac2{1-s} L^6}^\frac2{1-s} +\|\n \Wet\|_{L_t^2 L^\frac83}^2\Big)\Bigg\}.
 \label{EstimEnergieHs}
\end{multline}

\begin{rem}
 \sl{In the following sections, we will show that the first two terms from the right-hand side are small thanks to frequency truncations (the third term being small thanks to the assumptions), and we will use the Strichartz estimates from the appendix for all the other terms.}
\end{rem}

\subsection{Estimates for the frenquency truncations}
The aim of this section is to prove the following result.
\begin{prop}
 \sl{There exists $C_{\delta,\eta,\nu,\Co}>$ such that (we recall that $s=\frac12+\eta \delta$):
\begin{equation}
\|(Id-\cPrR) \Uoeosc\|_{L^2 \cap \dot{H}^s}^2 \leq C_{\delta,\eta,\nu,\Co}\left(\ee^{2\big(M\delta(1-\eta)-\gamma\big)} +\ee^{2\big(m\delta -M \delta(\frac12+\eta) -\gamma\big)}\right).
\label{estiminit}
\end{equation}
and
 \begin{equation}
 \|(Id-\cPrR) \tG\|_{L_t^1 (L^2 \cap \dot{H}^s)} \leq C_{\delta,\eta,\nu,\Co} \left(\ee^{M\delta(1-\eta)} +\ee^{\frac{m}3 -M(\frac23+\delta)}\right).
  \label{estimG}
 \end{equation}
}
\end{prop}
\begin{rem}
 \sl{In particular, the initial data satisfies
\begin{equation}
\|\de(0)\|_{\dot{H}^\frac12}^2 \leq C_{\delta,\eta,\nu,\Co}\left(\ee^{2\aa_0} +\ee^{2\big(M\delta(1-\eta)-\gamma\big)} +\ee^{2\big(m\delta-M \delta(\frac12+\eta) -\gamma\big)}\right).
\label{estiminit2}
\end{equation}
 }
\end{rem}
\textbf{Proof:} we use here the methods from \cite{FCPAA} and as described in \cite{FCStratif1}, contrary to the QG/oscillating structure, here we have $\bP_2 \Uoeosc=0=\bP_2 \tG$ which simplifies a little the computations.
$$
 \|(Id-\cPrR) \Uoeosc\|_{\dot{H}^s}^2 \leq 2\left(\|\big(Id-\chi(\frac{|D|}{\Re})\big) \Uoeosc\|_{\dot{H}^s}^2 +\|\chi(\frac{|D|}{\Re}) \chi(\frac{|D_h|}{2\re}) \Uoeosc\|_{\dot{H}^s}^2 \right).
$$
The first part is easily bounded thanks to Plancherel and the Bienaym\'e-Tchebychev estimates (we recall that $s=\frac12+\eta \delta$ for some $\eta\in]0,1[$):
\begin{multline}
 \|\big(Id-\chi(\frac{|D|}{\Re})\big) \Uoeosc\|_{\dot{H}^s}^2 \leq C\int_{|\xi|\geq \frac{\Re}2} |\xi|^{2s} |\hat{\Uoeosc}(\xi)|^2 d\xi\\
 \leq C \int_{|\xi|\geq \frac{\Re}2} \left( \frac{2|\xi|}{\Re}\right)^{\delta(2-2\eta)}|\xi|^{1+2\eta \delta} |\hat{\Uoeosc}(\xi)|^2 d\xi \leq \frac{C_{\delta, \eta}}{\Re^{2\delta(1-\eta)}} \|\Uoeosc\|_{\dot{H}^{\frac12+\delta}}^2,
\end{multline}
while the second one is dealt thanks to Lemma \ref{lemaniso} introducing $q=\frac2{1+\delta}\in]1,2[$:
\begin{multline}
 \|\chi(\frac{|D|}{\Re}) \chi(\frac{|D_h|}{2\re}) \Uoeosc\|_{\dot{H}^s}= \||D|^{\frac12+\eta \delta} \chi(\frac{|D|}{\Re}) \chi(\frac{|D_h|}{2\re}) \Uoeosc\|_{L^2}\\
 \leq \Re^{\eta \delta} \|\chi(\frac{|D|}{\Re}) \chi(\frac{|D_h|}{2\re}) |D|^{\frac12} \Uoeosc\|_{L^2} \leq \Re^{\eta \delta} (\Re (2\re)^2)^{\frac1{q}-\frac12} \||D|^{\frac12} \Uoeosc\|_{L^q}\\
 \leq C_{\delta} \Re^{\delta(\frac12+\eta)} \re^\delta \||D|^{\frac12} \Uoeosc\|_{L^q}
\end{multline}
Similarly the $L^2$-norms are bounded according to:
\begin{equation}
 \begin{cases}
 \vspace{0.2cm}
  \displaystyle{\|\big(Id-\chi(\frac{|D|}{\Re})\big) \Uoeosc\|_{L^2}^2 \leq \frac{C_{\delta}}{\Re^{1+2\delta}} \|\Uoeosc\|_{\dot{H}^{\frac12+\delta}}^2},\\
  \displaystyle{\|\chi(\frac{|D|}{\Re}) \chi(\frac{|D_h|}{2\re}) \Uoeosc\|_{L^2}
  \leq C_{\delta} \Re^{\frac{\delta}2} \re^\delta \|\Uoeosc\|_{L^q},}
 \end{cases}
\end{equation}
so that we can finally write that:
\begin{multline}
 \|(Id-\cPrR) \Uoeosc\|_{L^2}^2 +\|(Id-\cPrR) \Uoeosc\|_{\dot{H}^s}^2\\
 \leq C_{\delta, \eta} \left(\frac1{\Re^{2\delta(1-\eta)}} \|\Uoeosc\|_{\dot{H}^{\frac12+\delta}}^2 +\Re^{\delta (1+2\eta)} \re^{2 \delta} (\||D|^{\frac12} \Uoeosc\|_{L^q}^2 +\|\Uoeosc\|_{L^q}^2) \right),
\end{multline}
which proves the first point. The truncated external force is also split into:
$$
 \|(Id-\cPrR) \tG\|_{L_t^1 \dot{H}^s} \leq \|\big(Id-\chi(\frac{|D|}{\Re})\big) \tG\|_{L_t^1 \dot{H}^s} +\|\chi(\frac{|D|}{\Re}) \chi(\frac{|D_h|}{2\re}) \tG\|_{L_t^1 \dot{H}^s},
$$
and the first term is estimated similarly as before:
\begin{equation}
   \|\big(Id-\chi(\frac{|D|}{\Re})\big) \tG\|_{L_t^1 \dot{H}^s} \leq \frac{C_{\delta, \eta}}{\Re^{\delta(1-\eta)}} \|\tG\|_{L_t^1 \dot{H}^{\frac12+\delta}} \leq \frac{C_{\delta, \eta, \nu, \Co}}{\Re^{\delta(1-\eta)}},
\end{equation}
The second term is also bounded as in \cite{FCPAA}, in a simpler way than the corresponding part in the initial data (we use specifically \eqref{ecritG} and Sobolev injections):
\begin{multline}
 \|\chi(\frac{|D|}{\Re}) \chi(\frac{|D_h|}{2\re}) \tG\|_{L_t^1 \dot{H}^s}= \|\chi(\frac{|D|}{\Re}) \chi(\frac{|D_h|}{2\re}) |D|^s \tG\|_{L_t^1 L^2} \leq \Re^s (\Re (2\re)^2)^{\frac1{\frac32}-\frac12} \|\tG\|_{L_t^1 L^\frac32}\\
 \leq C \Re^{\frac12+\eta \delta +\frac16} \re^{\frac13} \|\tv^h \cdot \n_h \tv^h\|_{L_t^1 L^\frac32} \leq C \Re^{\frac23+\delta}  \re^{\frac13} \int_0^t \|\tv^h(\tau)\|_{L^6} \|\n \tv^h(\tau)\|_{L^2} d\tau\\
 \leq C \Re^{\frac23+\delta}  \re^{\frac13} \|\n \tv^h\|_{L_t^2 L^2}^2 \leq C_{\delta, \nu, \Co} \Re^{\frac23+\delta}  \re^{\frac13}.
\end{multline}
Similarly, the $L^2$-norms are bounded as follows:
\begin{equation}
 \begin{cases}
 \vspace{0.2cm}
  \displaystyle{\|\big(Id-\chi(\frac{|D|}{\Re})\big) \tG\|_{L_t^1 L^2} \leq \frac{C_{\delta, \nu, \Co}}{\Re^{\frac12+\delta}},}\\
  \displaystyle{\|\chi(\frac{|D|}{\Re}) \chi(\frac{|D_h|}{2\re}) \tG\|_{L_t^1 L^2}\leq C_{\delta, \nu, \Co} \Re^{\frac16}  \re^{\frac13},}
 \end{cases}
\end{equation}
so that
\begin{equation}
  \|(Id-\cPrR) \tG\|_{L_t^1 \dot{H}^s} +\|(Id-\cPrR) \tG\|_{L_t^1 L^2} \leq C_{\delta, \eta, \nu, \Co} \left(\frac1{\Re^{\delta(1-\eta)}} +\Re^{\frac23+\delta}  \re^{\frac13} \right),
\end{equation}
which ends the proof of the second point. $\blacksquare$

\subsection{Strichartz estimates}

Thanks to the Strichartz estimates proved in the appendix, we are able to bound in \eqref{EstimEnergieHs} each term featuring $\Wet$, as collected in the following proposition.
\begin{prop}
 \sl{There exists a constant $\Do= \Do(\nu,\nu',\Co, \delta,\eta)>0$ such that for any $0<\ee\leq \ee_1$ ($\ee_1$ is defined in Proposition \ref{estimvp}) and any  $t\leq T_\ee^*$ (we recall that $s=\frac12+\eta \delta$):
 \begin{equation}
 \begin{cases}
 \vspace{0.15cm}
 \|\n \Wet\|_{L_t^2 L^3} \leq \Do \frac{\Re^6}{\re^7} \ee^{\frac1{24}-\gamma} =\Do \ee^{\frac1{24}-\gamma-(6M+7m)},\\
 \vspace{0.15cm}
 \|\Wet\|_{L_t^4 L^6} +\|\Wet\|_{L_t^{\frac2{1-s}} L^6} \leq \Do \ee^{\frac1{12}-\gamma-(6M+7m)},\\
 \vspace{0.15cm}
 \|\Wet\|_{L_t^2 L^8} \leq \frac{\Do}{(\re\Re)^\frac18} \ee^{\frac3{32}-\gamma-(6M+7m)},\\
 \vspace{0.15cm}
\|\n \Wet\|_{L_t^2 L^{\frac83}} \leq \Do \ee^{\frac1{32}-\gamma-(6M+7m)},\\
\|\Wet\|_{L_t^8 L_{v,h}^{\infty,2}} \leq \frac{\Do}{\re^\frac14} \ee^{\frac18-\gamma-(6M+7m)}
 \end{cases}
  \label{Estimdispgen1}
 \end{equation}
 }
 \end{prop}
\textbf{Proof:} the result is a consequence of Propositions \ref{EstimStri1} and \ref{EstimStri1aniso} (only for the last term). Choosing $(d,p,r,q)=(1,2,3,2)$ we can write that (thanks to Theorem \ref{ThSNS}, Lemma \ref{lemaniso}, Propositions \ref{injectionLr}, \ref{Propermut}, and \ref{EstimStri1}) there exists a constant $C=C_{\nu,\nu'}$ and a constant $C(\Co, \nu, \delta)>0$ such that:
\begin{multline}
 \|\n \Wet\|_{L_t^2 L^3} \leq \|\n \Wet\|_{L_t^2 \dot{B}_{3,2}^0}  \leq \|\n \Wet\|_{\tilde{L}_t^2 \dot{B}_{3,2}^0} \leq C \frac{\Re^4}{\re^{\frac{31}6}} \ee^{\frac1{24}} \left(\|\cPrR \Uoeosc\|_{\dot{H}^1} +\|\cPrR \tG\|_{L^1 \dot{H}^1}\right)\\
 \leq C \frac{\Re^4}{\re^{\frac{31}6}} \ee^{\frac1{24}} \Re^{\frac12-\delta}\left(\|\cPrR \Uoeosc\|_{\dot{H}^{\frac12+\delta}} +\|\cPrR \tG\|_{L^1 \dot{H}^{\frac12+\delta}}\right)\\
 \leq C \frac{\Re^{\frac92-\delta}}{\re^{\frac{31}6}} \ee^{\frac1{24}}\left(\Co \ee^{-\gamma} +C(\Co, \nu, \delta)\right).
\end{multline}
Choosing $(d,p,r,q)=(0,4,6,2)$ and using in addition the Bienaym\'e-Tchebychev estimates, we obtain:
\begin{multline}
 \|\Wet\|_{L_t^4 L^6} \leq \|\n \Wet\|_{\tilde{L}_t^4 \dot{B}_{6,2}^0} \leq C \frac{\Re^{\frac{11}2}}{\re^{\frac{35}6}} \ee^{\frac1{12}} \left(\|\cPrR \Uoeosc\|_{L^2} +\|\cPrR \tG\|_{L^1 L^2}\right)\\
 \leq C \frac{\Re^{\frac{11}2}}{\re^{\frac{35}6}} \ee^{\frac1{12}} \left(\re^{-\frac12}\|\cPrR \Uoeosc\|_{\dot{H}^\frac12} +\|\cPrR \tG\|_{L^1 L^2}\right)\\
 \leq C \frac{\Re^{\frac{11}2}}{\re^{\frac{19}3}} \ee^{\frac1{12}}\left(\Co \ee^{-\gamma} +C(\Co, \nu, \delta)\right).
\end{multline}
Taking $(d,p,r,q) \in \{(0,2,8,2), (1,2,\frac83,2), (0, \frac2{1-s},6,2)\}$, with the same arguments, we end-up with (we recall that $s=\frac12+\eta \delta$):
$$
\begin{cases}
\|\Wet\|_{L_t^2 L^8} \leq C_{\nu,\nu',\delta,\Co} \displaystyle{\frac{\Re^{\frac{47}8}}{\re^{\frac{57}8}}} \ee^{\frac3{32}-\gamma} \leq \Do \frac{\Re^6}{\re^7} \frac1{(\re\Re)^\frac18} \ee^{\frac3{32}-\gamma},\\
\|\n \Wet\|_{L_t^2 L^{\frac83}} \leq C_{\nu,\nu',\delta,\Co} \displaystyle{\frac{\Re^{\frac{33}8-\delta}}{\re^{\frac{39}8}}} \ee^{\frac1{32}-\gamma} \leq \Do \frac{\Re^6}{\re^7} \ee^{\frac1{32}-\gamma},\\
\|\Wet\|_{L_t^{\frac2{1-s}} L^6} \leq C_{\nu,\nu',\delta,\Co} \displaystyle{\frac{\Re^{\frac{11}2}}{\re^{\frac{19}3-\eta \delta}}} \ee^{\frac1{12}-\gamma} \leq \Do \frac{\Re^6}{\re^7} \ee^{\frac1{12}-\gamma}.
\end{cases}
$$
\begin{rem}
 \sl{In all the previous estimates, the condition $p\leq \frac8{1-\frac2{r}}$ is obvious, except for the last term, wich requires that $\eta \delta \leq \frac13$ (we recall that we already ask $\eta\leq \frac12$).
 }
 \label{Condetadelta}
\end{rem}
The anisotropic term is dealt with the same arguments but using Proposition \ref{EstimStri1aniso} and for $m=\infty$, we obtain that:
$$
\|\Wet\|_{L_t^8 L_{v,h}^{\infty,2}} \leq C_{\nu,\nu',\delta,\Co} \frac{\Re^6}{\re^{\frac{29}4}} \ee^{\frac18-\gamma}\leq \Do \frac{\Re^6}{\re^7} \frac1{\re^\frac14} \ee^{\frac18-\gamma}. \quad \blacksquare
$$

\subsection{Bootstrap and convergence}

We are now able to finish the bootstrap argument. Into \eqref{EstimEnergieHs} we inject, on one hand \eqref{estiminit}, \eqref{estimG} (to deal with the first two terms from the right-hand side of \eqref{estiminit}), and on the other hand \eqref{Estimdispgen1} (to show that the terms involving $\Wet$ are small thanks to Strichartz estimates), so that, uniformly denoting from line to line as $\Do$ a constant depending on $(\nu,\nu',\Co, \delta,\eta)$, we obtain that for all $0<\ee\leq \ee_1$ (defined in Proposition \ref{estimvp}) and all $t\leq T_{\ee,2}$ (defined in \eqref{DefTe}),
\begin{multline}
 \|\de(t)\|_{H^s}^2+\frac{\nu_0}2 \int_0^t \|\n \de(t')\|_{H^s}^2 dt' \leq \exp{\Bigg\{\Do \Big(1 +(1+\ee^{\frac1{48}}) E_\ee^2 +\big(\ee^{\frac5{96}} E_\ee \big)^4 +\big(\ee^{\frac5{96}} E_\ee \big)^{\frac4{1-2\eta \delta}}\Big)\Bigg\}}\\
 \Bigg[\ee^{2\aa_0} +\ee^{2(M\delta(1-\eta)-\gamma)} +\ee^{2(m\delta -M \delta(\frac12+\eta) -\gamma)} +\ee^{M\delta(1-\eta)} +\ee^{\frac{m}3 -M(\frac23+\delta)}\\
 +\Do \Big(\ee^{\frac1{48}} +\ee^{\frac5{48}} +\frac{\ee^{\frac18}}{(\re\Re)^\frac14} +\frac{\ee^{\frac3{16}}}{\re^\frac12}\Big) E_\ee^2\Bigg],
\end{multline}
where we have introduced the small quantity $E_\ee =\ee^{\frac1{32}-\gamma-(6M+7m)}$. Thanks to Remark \ref{Condetadelta} we have $4 \leq \frac4{1-2\eta \delta} \leq 12$, and if we ask that:
\begin{equation}
\frac1{32}-\gamma-(6M+7m)\geq 0,
 \label{CondmMgamma1}
\end{equation}
then we are sure that if $\ee>0$ is small enough, $E_\ee \leq 1$ and  $1 +(1+\ee^{\frac1{48}}) E_\ee^2 +\big(\ee^{\frac5{96}} E_\ee \big)^4 +\big(\ee^{\frac5{96}} E_\ee \big)^{\frac4{1-2\eta \delta}} \leq 5$, which implies that for all $t\leq T_{\ee,2}$:
\begin{multline}
 \|\de(t)\|_{H^s}^2+\frac{\nu_0}2 \int_0^t \|\n \de(t')\|_{H^s}^2 dt' \leq \Do e^{5\Do}
 \Bigg[\ee^{2\aa_0} +\ee^{2\big(M\delta(1-\eta)-\gamma\big)} +\ee^{2\big(m\delta -M \delta(\frac12+\eta) -\gamma\big)}\\
 +\ee^{M\delta(1-\eta)} +\ee^{\frac{m}3 -M(\frac23+\delta)} +\ee^{\frac1{48}} +\ee^{\frac14 \big(\frac12+M-m\big)} +\ee^{\frac3{16}-\frac{m}2}\Bigg].
\end{multline}
If we observe that
$$
\begin{cases}
 \frac14 \big(\frac12+M-m\big) \geq \frac1{48},\\
 \frac3{16}-\frac{m}2 \geq \frac1{48},
\end{cases}
\Longleftrightarrow \quad m\leq \frac13 \mbox{ and } m-M\leq \frac5{12} \quad \Longleftarrow \eqref{CondmMgamma1},
$$
then for all $t\leq T_{\ee,2}$,
$$
 \|\de(t)\|_{H^s}^2+\frac{\nu_0}2 \int_0^t \|\n \de(t')\|_{H^s}^2 dt' \leq \Do \ee^{2 N(\ee)},
 $$
 where
 \begin{multline}
 N(\ee) \overset{def}{=} \min \left(\aa_0, M\delta(1-\eta)-\gamma, m\delta -M \delta(\frac12+\eta) -\gamma, M\frac{\delta}2 (1-\eta), \frac{m}6 -M(\frac13+\frac{\delta}2), \frac1{96}\right).
\end{multline}
Using that $\delta\in]0,1]$, $\eta \leq \frac12$, and asking that $\gamma \leq M\frac{\delta}2 (1-\eta)$, we obtain that:
$$
N(\ee) \geq \min \left(\aa_0, M\frac{\delta}2 (1-\eta), (m -\frac32 M) \delta, \frac{m-5M}6, \frac1{96}\right).
$$
Choosing $M=\frac{m}6$ the previous estimates turns into:
$$
N(\ee) \geq \min \left(\aa_0, \frac{m\delta}{12} (1-\eta), \frac34 m\delta, \frac{m}{36}, \frac1{96}\right) \geq \min \left(\aa_0, \frac{m\delta}{12} (1-\eta), \frac{m}{36}\right).
$$
With these choices for $M$ and $\gamma$, we also have
$$
m\leq \frac1{259} \Longrightarrow 6M+ 7m+\gamma \leq \frac{97}{12} m \leq \frac1{32} \Longrightarrow \eqref{CondmMgamma1},
$$
So that, choosing $(m,M,\gamma)=(\frac1{259}, \frac1{1554},\frac{\delta}{3108}(1-\eta))$, if $\eta\delta \leq \frac13$ and $\eta\leq \frac12$, we finally obtain that for all $t\leq T_{\ee,2}$,
$$
 \|\de(t)\|_{H^s}^2+\frac{\nu_0}2 \int_0^t \|\n \de(t')\|_{H^s}^2 dt' \leq \Do \ee^{2 \min \left(\aa_0, \frac{\delta}{3108} (1-\eta), \frac1{9324}\right)}.
$$
Assuming that $T_{\ee,2} <T_{\ee}^*$, if $\ee>0$ is so small that $\Do \ee^{2 \min \left(\aa_0, \frac{\delta}{3108} (1-\eta), \frac1{9324}\right)} \leq (\frac{\nu_0}{8C})^2$ then the previous estimates implies that in particular $\|\de(T_{\ee,2})\|_{H^s} \leq \frac{\nu_0}{8C}$ which contradicts the definition of $T_{\ee,2}$ (see \eqref{DefTe}). We have proved by contradiction that $T_{\ee,2} =T_{\ee}^*$ and from the previous estimates,
$$
\int_0^{T_{\ee}^*} \|\n \de(t')\|_{H^s}^2 dt' <\infty,
$$
which entails, by the usual blow-up criterion, that $T_{\ee}^*=\infty$. Moreover we have obtained that for all $t\geq 0$,
$$
  \|\de(t)\|_{H^{\frac12+\eta \delta}}^2+\frac{\nu_0}2 \int_0^t \|\n \de(t')\|_{H^{\frac12+\eta \delta}}^2 dt' \leq \Do \ee^{2 \min \left(\aa_0, \frac{\delta}{3108} (1-\eta), \frac1{9324}\right)},
$$
which implies (thanks to Proposition \ref{injectionLr} and Lemma \ref{estimBsHs}) that
\begin{multline}
\|\de\|_{L_t^2 L^\infty} \leq C \|\de\|_{L_t^2 \dot{B}_{2,1}^\frac32} \leq C\|\n \de\|_{L_t^2 \dot{H}^{\frac12-\delta}}^\frac12 \|\n\de\|_{L_t^2 \dot{H}^{\frac12+\delta}}^\frac12\\
\leq C \|\n\de\|_{L_t^2 H^{\frac12+\delta}} \leq \Do \ee^{\min \left(\aa_0, \frac{\delta}{3108} (1-\eta), \frac1{9324}\right)}.
  \label{estimdHsfinale2}
\end{multline}
To finish the proof we use once more the Strichartz estimates from Proposition \ref{EstimStri1} with $(d,p,r,q)=(0,2,\infty,1)$,
\begin{multline}
 \|\Wet\|_{L_t^2 L^\infty} \leq \|\Wet\|_{L_t^2 \dot{B}_{\infty,1}^0} \leq C \frac{\Re^7}{\re^{\frac{15}2}} \ee^{\frac18} \left(\|\cPrR \Uoeosc\|_{\dot{B}_{2,1}^0} +\|\cPrR \tG\|_{L^1 \dot{B}_{2,1}^0}\right)\\
 \leq C \frac{\Re^7}{\re^{\frac{15}2}} \ee^{\frac18} \re^{-\frac12}\left(\|\Uoeosc\|_{\dot{B}_{2,1}^\frac12} +\|\cPrR \tG\|_{L^1 \dot{B}_{2,1}^\frac12}\right)\\
 \leq C \frac{\Re^7}{\re^8} \ee^{\frac18} \left(\|\Uoeosc\|_{\dot{H}^{\frac12-c\delta}}^\frac12 \|\Uoeosc\|_{\dot{H}^{\frac12+c\delta}}^\frac12 +\|\tG\|_{L^1 \dot{H}^{\frac12-c\delta}}^\frac12 \|\tG\|_{L^1 \dot{H}^{\frac12+c\delta}}^\frac12\right)\\
 \leq \Do \frac{\Re^7}{\re^8} \ee^{\frac18} (\ee^{-\gamma}+1) \leq \Do \ee^{\frac18-\gamma-(7M+8m)}.
 \label{estimBsHs2}
\end{multline}
With the previous choices for $(m,M,\gamma)$,
$$
\frac18-\gamma-(7M+8m) \geq \frac18 -(\frac1{12}+8+\frac76)\frac1{259} =\frac{555}{6216},
$$
so that we end-up with
\begin{multline}
\|\Ue-(\tv,0, \tThee)\|_{L_t^2 L^\infty} =\|\De\|_{L_t^2 L^\infty} =\|\de -\Wet\|_{L_t^2 L^\infty}\\
\leq \Do (\ee^{\min \left(\aa_0, \frac{\delta}{3108} (1-\eta), \frac1{9324}\right)} +\ee^{\frac{555}{6216}}) \leq 2\Do \ee^{\min \left(\aa_0, \frac{\delta}{3108} (1-\eta), \frac1{9324}\right)},
\end{multline}
and the proof of Theorem \ref{Th1} is complete. $\blacksquare$

\section{Proof of Theorem \ref{Th2}}
As usual, in the particular case $\nu=\nu'$, we can take advantage of simplifications: the computation of the eigenvalues for the linearized system does not require anymore truncations in frequency, and the projectors $\bP_3$ and $\bP_4$ become orthogonal.
\label{PrvTh2}
\subsection{A priori estimates}

Let us consider $\de=\De-\We$ where $\We$ solves \eqref{We}:

\begin{equation}
 \begin{cases}
  \d_t \de -L\de +\frac{1}{\ee} \mathbb{P} \mathcal{B} \de = \Sum_{i=1}^{10} G_i,\\
  {\de}_{|t=0}= (\UoeS^h-\tv^h,0,0),
 \end{cases}
\label{denunu}
\end{equation}
where:
\begin{equation}
\begin{cases}
 G_1 \overset{def}{=}-\mathbb{P}(\de \cdot \n \de), \quad G_2 \overset{def}{=}-\mathbb{P}(\de \cdot \n \tv^h,0,0), \quad G_3 \overset{def}{=}-\mathbb{P}(\tv^h \cdot \n_h \de),\\
 G_4 \overset{def}{=}-\mathbb{P}(\de \cdot \nabla \We),\quad G_5 \overset{def}{=}-\mathbb{P}(\We \cdot \n \de),\quad G_6 \overset{def}{=}-\mathbb{P}(\tv^h \cdot \n_h \We),\\
 G_7 \overset{def}{=}-\mathbb{P}(\We \cdot \n \tv^h,0,0),\quad G_8 \overset{def}{=}-\mathbb{P}(\We \cdot \n \We),\\
 G_9 \overset{def}{=}-\mathbb{P}(0,0,0,\de^3 \cdot \d_3 \tThee),\quad G_{10} \overset{def}{=}-\mathbb{P}(0,0,0,\We^3 \cdot \d_3 \tThee).
\end{cases}
 \label{systdenunu}
\end{equation}
Following the same steps as in the general case, we obtain that for all $t\leq T_{\ee,2}$ (where $T_{\ee,2}$ is the same as in \eqref{DefTe}),
\begin{multline}
 \|\de(t)\|_{H^s}^2+\frac{\nu}2 \int_0^t \|\n \de(t')\|_{H^s}^2 dt'\\ \leq \Bigg[\|\UoeS^h-\tv^h\|_{L^2 \cap \dot{H}^s}^2 +\Do \left(\|\n \We\|_{L_t^2 L^3}^2 +\|\We\|_{L_t^4 L^6}^2 +\|\We\|_{L_t^2 L^8}^2 +\|\We\|_{L_t^8 L_{v,h}^{\infty, 2}}^2\right)\Bigg]\\
 \times \exp \Bigg\{\Do \Big(1 +\|\n \We\|_{L_t^2 L^3}^2 +\|\We\|_{L_t^4 L^6}^4 +\|\We\|_{L_t^\frac2{1-s} L^6}^\frac2{1-s} +\|\n \We\|_{L_t^2 L^\frac83}^2\Big)\Bigg\},
 \label{EstimEnergieHsnunu}
\end{multline}

\subsection{Strichartz estimates}
We will prove in this section the following result:

\begin{prop}
 \sl{There exists a constant $\Do=\Do(\nu,\Co,\delta,\eta)>0$ such that for any $t\geq 0$,
 \begin{equation}
  \begin{cases}
  \vspace{0.1cm}
   \|\n \We\|_{L_t^2 L^3} +\|\We\|_{L_t^4 L^6} \leq \Do \ee^{\frac{\delta}2} \big(\|\Uoeosc\|_{\dot{H}^{\frac12-c\delta} \cap \dot{H}^{\frac12+\delta}}+1\big),\\
   \vspace{0.1cm}
   \|\We\|_{L_t^{\frac2{1-s}} L^6} \leq \Do \ee^{(1-\eta)\frac{\delta}2} \big(\|\Uoeosc\|_{\dot{H}^{\frac12-c\delta} \cap \dot{H}^{\frac12+\delta}}+1\big),\\
   \ee^{-\frac18}\|\We\|_{L_t^2 L^8} +\|\n \We\|_{L_t^2 L^\frac83} + \ee^{-\frac1{16}}\|\We\|_{L_t^8 L_{v,h}^{\infty, 2}} \leq \Do \ee^{\frac1{16}} \big(\|\Uoeosc\|_{\dot{H}^{\frac12-c\delta} \cap \dot{H}^{\frac12+\delta}}+1\big).\\
  \end{cases}
 \end{equation}
\label{Estimdispgen2}
 }
\end{prop}
\textbf{Proof:} using Proposition \ref{EstimStri1aniso} with $(d,p,r,q, \theta)=(1,2,3,2, \frac{\delta}6)$, we obtain (with the same arguments as in the general case) that there exists $C=C(\nu,\delta)>0$ such that for any $t\geq 0$,
\begin{multline}
 \|\n \We\|_{L_t^2 L^3} \leq \|\n \We\|_{L_t^2 \dot{B}_{3,2}^0} \leq \|\n \We\|_{\tilde{L}_t^2 \dot{B}_{3,2}^0} \leq C \ee^{\frac{\delta}2} \left(\|\Uoeosc\|_{\dot{H}^{\frac12+\delta}} +\|\tG\|_{L^1 \dot{H}^{\frac12+\delta}}\right)\\
 \leq C \ee^{\frac{\delta}2} \left(\|\Uoeosc\|_{\dot{H}^{\frac12-c\delta} \cap \dot{H}^{\frac12+\delta}} +C(\Co, \nu, \delta)\right) \leq \Do \ee^{\frac{\delta}2} \big(\|\Uoeosc\|_{\dot{H}^{\frac12-c\delta} \cap \dot{H}^{\frac12+\delta}}+1\big).
\end{multline}
This choice for $\theta$ requires that $\delta\leq \frac16$, and the condition $p\leq \frac4{\theta(1-\frac2{r})}$ is trivially satisfied. The second and third estimates are obtained similarly, applying the same proposition successively for $(d,p,r,q, \theta) =(0,4,6,2,3\delta)$ and $(0,\frac2{1-s}, 6,2,3(1-\eta)\delta)$ (and does not require any additionnal assumption as we already have $\delta \leq \frac16$).

With $(d,p,r,q, \theta)=(0,2,8,2,1)$ (we took $\theta=1$ as the Sobolev index is $\sigma_1=\frac{1+3\theta}8$) we obtain that
$$
\|\We\|_{L_t^2 L^8} \leq \Do \ee^{\frac3{16}}\big(\|\Uoeosc\|_{\dot{H}^{\frac12-c\delta} \cap \dot{H}^{\frac12+\delta}}+1\big),
$$
which gives the fourth estimates. The Fifth estimates is obtained choosing $(d,p,r,q, \theta)=(1,2,\frac83,2,1)$. As it involves anisotropic spaces, we use Proposition \ref{EstimStri2aniso} with $(d,p,m, \theta)=(0,8,\infty,1)$, Theorem \ref{ThSNS}, and obtain, combining the arguments from \eqref{estimBsHs2} with interpolation, that:
\begin{multline}
 \|\We\|_{L_t^8 L_{v,h}^{\infty, 2}} \leq C_{\nu} \ee^{\frac18} \left(\|\Uoeosc\|_{\dot{B}_{2,1}^\frac12} +\|\tG\|_{L^1 \dot{B}_{2,1}^\frac12}\right)\\
\leq C_{\nu} \ee^{\frac18} \left(\|\Uoeosc\|_{\dot{H}^{\frac12-c\delta}}^\frac12 \|\Uoeosc\|_{\dot{H}^{\frac12+c\delta}}^\frac12 +\|\tG\|_{L^1 \dot{H}^{\frac12-c\delta}}^\frac12 \|\tG\|_{L^1 \dot{H}^{\frac12+c\delta}}^\frac12\right)\\
\leq C_{\nu} \ee^{\frac18} \left(\|\Uoeosc\|_{\dot{H}^{\frac12-c\delta} \cap \dot{H}^{\frac12+\delta}} +C(\Co, \nu, \delta)\right).
\end{multline}

\subsection{Results when $s=\frac12$}

When we only assume that there exists $c,m_0>0$ such that:
$$
\|\Uoeosc\|_{\dot{H}^{\frac12-c\delta} \cap \dot{H}^{\frac12+\delta}} \leq m_0 \ee^{-\frac{\delta}2},
$$
gathering the Strichartz estimates from the previous section into \eqref{EstimEnergieHsnunu} entails that for any $t\leq T_{\ee,2}$ (in the present case $\eta=0$),
\begin{multline}
\|\de(t)\|_{H^\frac12}^2+\frac{\nu}2 \int_0^t \|\n \de(t')\|_{H^\frac12}^2 dt' \leq \Do e^{\Do \left\{1+ (m_0+\ee^{\frac{\delta}2})^2 \Big(1+\ee^{\frac18-\delta}\Big) +2(m_0+\ee^{\frac{\delta}2})^4 \right\}}\\
\times \left[\ee^{2\aa_0} +(m_0+\ee^{\frac{\delta}2})^2\Big(1+\ee^{\frac3{16}-\delta} +\ee^{\frac14-\delta}\Big) \right].
\end{multline}
If we choose $\ee,m_0>0$ so small that (we recall that $\delta\leq \frac18$):
$$
\begin{cases}
\vspace{0.2cm}
 (m_0+\ee^{\frac{\delta}2})^2 \Big(1+\ee^{\frac18-\delta}\Big) +2(m_0+\ee^{\frac{\delta}2})^4 \leq 1,\\
 \Do e^{2\Do} \left[\ee^{2\aa_0} +(m_0+\ee^{\frac{\delta}2})^2\Big(1+\ee^{\frac3{16}-\delta} +\ee^{\frac14-\delta}\Big) \right] \leq \left(\frac{\nu}{8C} \right)^2,
\end{cases}
$$
then we prove as in the general case that $T_{\ee,2}=T_\ee^*=\infty$.
If in addition $m_0$ is replaced by some $m(\ee)\underset{\ee \rightarrow 0}{\longrightarrow} 0$, we obtain that when $\ee>0$ is small enough:
$$
\|\de\|_{E^{\frac12}} =\|\de\|_{\dot{E}^0 \cap \dot{E}^{\frac12}}  \leq \Do \max \big(\ee^{\aa_0}, \ee^{\frac{\delta}2}, m(\ee)\big).
$$

\subsection{Precise convergence rates}

With the following stronger assumption,
$$
\|\Uoeosc\|_{\dot{H}^{\frac12-c\delta} \cap \dot{H}^{\frac12+\delta}} \leq \Co \ee^{-\gamma},
$$
the Strichartz estimates from Proposition \ref{Estimdispgen2} now become when we introduce $\eta_0>0$ so that $\gamma=\frac{\delta}2 (1-2\eta_0)$ (we also recall that $\delta\leq \frac18$):
$$
\begin{cases}
  \vspace{0.2cm}
   \|\n \We\|_{L_t^2 L^3} +\|\We\|_{L_t^4 L^6} \leq \Do \ee^{\frac{\delta}2 -\gamma} =\Do \ee^{\eta_0 \delta},\\
   \vspace{0.1cm}
   \|\We\|_{L_t^{\frac2{1-s}} L^6} \leq \Do \ee^{(1-\eta)\frac{\delta}2 -\gamma} =\Do \ee^{(\eta_0-\frac{\eta}2) \delta},\\
   \ee^{-\frac18}\|\We\|_{L_t^2 L^8} +\|\n \We\|_{L_t^2 L^\frac83} + \ee^{-\frac1{16}}\|\We\|_{L_t^8 L_{v,h}^{\infty, 2}} \leq \Do \ee^{\frac1{16}-\gamma} =\Do \ee^{\frac12(\frac18-\delta) +\eta_0 \delta} \leq \Do \ee^{\eta_0 \delta}.\\
  \end{cases}
$$
Gathering these estimates in \eqref{EstimEnergieHsnunu} we have that for any $t\leq T_{\ee,2}$ (here $s=\frac12+\eta \delta$ and as in the general case $\frac2{1-s}\geq4$),
\begin{multline}
\|\de(t)\|_{H^s}^2+\frac{\nu}2 \int_0^t \|\n \de(t')\|_{H^s}^2 dt' \leq \Do e^{\Do \left\{1+ 2\ee^{2\eta_0 \delta} +\ee^{4\eta_0 \delta}  +\ee^{4(\eta_0-\frac{\eta}2)\delta} \right\}}\\
\times \left[\ee^{2\aa_0} +\ee^{2\eta_0 \delta}\Big(2+\ee^{\frac14}+\ee^{\frac18}\Big) \right].
\end{multline}
Now, as we need $\eta\leq \min(2\eta_0,\frac12)$ (with $ \eta_0\in ]0,\frac12$), we can simply choose $\eta= \eta_0$ (now $s=\frac12+\eta_0\delta=\frac12 +\frac{\delta}2-\gamma$) and as $\ee\in]0,1]$ then
$$
\|\de(t)\|_{H^s}^2+\frac{\nu}2 \int_0^t \|\n \de(t')\|_{H^s}^2 dt' \leq \Do e^{5\Do} \left[\ee^{2\aa_0} +4 \ee^{2\eta_0 \delta} \right] \leq \Do \ee^{2 \min(\aa_0, \eta_0 \delta)}.
$$
Once more this allows us to prove that $T_{\ee,2}=T_\ee^*=\infty$ and as the previous estimates is now valid for any $t\geq 0$ we obtain that:
$$
\|\de\|_{E^s} =\|\de\|_{E^{\frac12+\eta_0 \delta}} =\|\de\|_{\dot{E}^0 \cap \dot{E}^{\frac12+\eta_0 \delta}}  \leq \Do \ee^{\min(\aa_0, \eta_0 \delta)}.
$$
As in the general case, using Proposition \ref{injectionLr} and Lemma \ref{estimBsHs}, we have (we recall that $\gamma=\frac{\delta}2 (1-2\eta_0)$ with $\eta_0\in]0,\frac12[$)
\begin{multline}
\|\de\|_{L_t^2 L^\infty} \leq C \|\de\|_{L_t^2 \dot{B}_{2,1}^\frac32} \leq C\|\n \de\|_{L_t^2 \dot{H}^{\frac12-\eta_0\delta}}^\frac12 \|\n\de\|_{L_t^2 \dot{H}^{\frac12+\eta_0\delta}}^\frac12\\
\leq C \|\n \de\|_{L_t^2 H^{\frac12+\eta_0\delta}} \leq \Do \ee^{\min(\aa_0, \eta_0 \delta)} =\Do \ee^{\min(\aa_0, \frac{\delta}2-\gamma)}.
\label{Estimde}
\end{multline}
All that remains is then to use Proposition \ref{EstimStri2} with $(d,p,r,q)=(0,2,\infty,1)$ and obtain that for any $\theta\in[0,1]$, and $t\geq 0$
\begin{equation}
  \|\We\|_{L_t^2 L^\infty} \leq \|\We\|_{L_t^2 \dot{B}_{\infty,1}^0}  \leq \|\We\|_{\tilde{L}_t^2 \dot{B}_{\infty,1}^0} \leq C_{\theta,\nu} \ee^{\frac{\theta}4} \left(\|\Uoeosc\|_{\dot{B}_{2,1}^{\frac12+\frac{\theta}2}} +\|\tG\|_{L^1 \dot{B}_{2,1}^{\frac12+\frac{\theta}2}}\right).
  \label{EstimWe}
\end{equation}
As in \cite{FCcompl} and \cite{FCRF}, applying Lemma \ref{estimBsHs} with $(\aa,\bb)=(\frac{\theta}2(1-a), \frac{\theta}2(1+b))$ for $a,b>0$ (and $b$ small as we will see in what follows), there exists a constant $C=C(a,b,\theta)>0$ such that for any function we have:
$$
\|f\|_{\dot{B}_{2,1}^{\frac12+\frac{\theta}2}} \leq C\|f\|_{\dot{H}^{\frac12+\frac{\theta}2(1-a)}}^{\frac{b}{a+b}} \|f\|_{\dot{H}^{\frac12+\frac{\theta}2(1+b)}}^{\frac{a}{a+b}}.
$$
Trying to use the assumptions we will choose $a,b$ so that
$$
(\frac{\theta}2(1-a), \frac{\theta}2(1+b))=(-c\delta, \delta),
$$
which is realized when $\theta=\frac{2\delta}{1+b}$ (this is possible as we already ask $\delta \leq \frac18$) and $a=1+c(1+b)$ so that we obtain:
\begin{multline}
 \|\Uoeosc\|_{\dot{B}_{2,1}^{\frac12+\frac{\theta}2}} \leq C_{b,c,\delta} \|\Uoeosc\|_{\dot{H}^{\frac12-c\delta}}^{\frac{b}{(1+b)(1+c)}} \|\Uoeosc\|_{\dot{H}^{\frac12+\delta}}^{\frac{1+c(1+b)}{(1+b)(1+c)}}\\
 \leq C_{b,c,\delta} \|\Uoeosc\|_{H^{\frac12+\delta}} \leq C_{b,c,\delta,\Co} \ee^{-\gamma}.
\end{multline}
Similarly, we obtain that
$$
\|\tG\|_{L^1 \dot{B}_{2,1}^{\frac12+\frac{\theta}2}} \leq \|\tG\|_{L^1 H^{\frac12+\delta}} \leq C_{\delta,\nu,\Co}.
$$
Gathering the previous estimates into \eqref{EstimWe},
$$
\|\We\|_{L_t^2 L^\infty} \leq C_{\nu,\Co,\delta,b} \ee^{\frac{\delta}{2(1+b)}-\gamma} =C_{\nu,\Co,\delta,b} \ee^{\frac{\delta}2 (\frac1{1+b}-1+2\eta_0)}=C_{\nu,\Co,\delta,b} \ee^{\frac{\delta}2 (2\eta_0-\frac{b}{1+b})}.
$$
When some $k\in]0,1[$ is given (as close to $1$ as we wish), choosing $b=\frac{2\eta_0(1-k)}{1-2\eta_0(1-k)}$ we finally get that:
$$
\|\We\|_{L_t^2 L^\infty} \leq C_{\nu,\Co,\delta,\eta_0,k} \ee^{k\eta_0 \delta}.
$$
Combining this with \eqref{Estimde}, we finally obtain that
\begin{multline}
 \|\Ue-(\tv,0, \tThee)\|_{L_t^2 L^\infty} =\|\De\|_{L_t^2 L^\infty} =\|\de -\Wet\|_{L_t^2 L^\infty}\\
 \leq \Do \ee^{\min(\aa_0, \eta_0 \delta)} +C_{\nu,\Co,\delta,\eta_0,k} \ee^{k\eta_0 \delta} \leq \Do \ee^{\min(\aa_0, k\eta_0 \delta)},
\end{multline}
which concludes the proof of the theorem. $\blacksquare$

\section{Appendix}

\subsection{Notations, Sobolev spaces and Littlewood-Paley decomposition}

As in \cite{FCPAA, FCcompl, FCStratif1}, this section roughly presents the spaces and norms that we will use. For a complete presentation of the Sobolev spaces and the Littlewood-Paley decomposition, we refer to \cite{Dbook}. Let us just recall that if $\phi:\R_+\rightarrow \R$ a smooth function supported in the ball $[0,\frac43]$, equal to 1 in a neighborhood of $[0,\frac43]$ and nonincreasing over $\R_+$. If we set $\varphi(r)=\phi(r/2)-\phi(r)$, then $\varphi$ is compactly supported in the set $\cC=[\frac34, \frac83]$ and we define the homogeneous dyadic blocks: for all $j\in \Z$,
$$
\ddj u:= \varphi(2^{-j}|D|)u =2^{jd} h(2^j.)* u, \quad \mbox{with } h(x)=\cF^{-1} \big(\varphi(|\xi|)\big).
$$
We recall that $\hat{k(D)u}(\xi)=k(\xi) \hat{u} (\xi)$ and we can define the homogeneous Besov norms and spaces:
\begin{defi}
\sl{For $s\in\R$ and $1\leq p,r\leq\infty,$ we set
$$
\|u\|_{\dot B^s_{p,r}}:=\bigg(\sum_{l\in \Z} 2^{rls}
\|\ddl u\|^r_{L^p}\bigg)^{\frac{1}{r}}\ \text{ if }\ r<\infty
\quad\text{and}\quad
\|u\|_{\dot B^s_{p,\infty}}:=\sup_{l} 2^{ls}
\|\ddl u\|_{L^p}.
$$
The homogeneous Besov space $\dot B^s_{p,r}$ is the subset of tempered distributions such that $\lim_{j \rightarrow -\infty} \|\dot{S}_j u\|_{L^\infty}=0$ and $\|u\|_{\dot B^s_{p,r}}$ is finite (where $\dot{S}_j u=\Sum_{l\leq j-1} \ddl u=\phi(2^{-j}|D|)u$).
}
\end{defi}
Let us first mention the following lemma:
\begin{prop}
 \sl{(\cite{Dbook} Chapter 2) The following continuous injections hold:
$$
 \begin{cases}
\mbox{For any } p\geq 1, & \dot{B}_{p,1}^0 \hookrightarrow L^p,\\
\mbox{For any } p\in[2,\infty[, & \dot{B}_{p,2}^0 \hookrightarrow L^p,\\
\mbox{For any } p\in[1,2], & \dot{B}_{p,p}^0 \hookrightarrow L^p.
\end{cases}
$$
}
 \label{injectionLr}
\end{prop}
Sometimes it is more convenient to work in a slight modification of the classical $L_t^p \dot{B}_{q,r}^s$ Spaces: the Chemin-Lerner time-space Besov spaces. As explained in the following definition, the integration in time is performed before the summation with respect to the frequency decomposition index:
\begin{defi} \cite{Dbook}
 \sl{For $s,t\in \R$ and $a,b,c\in[1,\infty]$, we define the following norm
 $$
 \|u\|_{\tilde{L}_t^a \dot{B}_{b,c}^s}= \Big\| \left(2^{js}\|\ddj u\|_{L_t^a L^b}\right)_{j\in \Z}\Big\|_{l^c(\Z)}.
 $$
 The space $\tilde{L}_t^a \dot{B}_{b,c}^s$ is defined as the set of tempered distributions $u$ such that $\lim_{j \rightarrow -\infty} S_j u=0$ in $L^a([0,t],L^\infty(\R^d))$ and $\|u\|_{\tilde{L}_t^a \dot{B}_{b,c}^s} <\infty$.
 }
 \label{deftilde}
\end{defi}
We refer once more to \cite{Dbook} (Section 2.6.3) for more details and will only recall the following proposition:
\begin{prop}
\sl{
For all $a,b,c\in [1,\infty]$ and $s\in \R$:
     $$
     \begin{cases}
    \mbox{if } a\leq c,& \forall u\in L_t^a \dot{B}_{b,c}^s, \quad \|u\|_{\tilde{L}_t^a \dot{B}_{b,c}^s} \leq \|u\|_{L_t^a \dot{B}_{b,c}^s}\\
    \mbox{if } a\geq c,& \forall u\in\tilde{L}_t^a \dot{B}_{b,c}^s, \quad \|u\|_{\tilde{L}_t^a \dot{B}_{b,c}^s} \geq \|u\|_{L_t^a \dot{B}_{b,c}^s}.
     \end{cases}
     $$
     \label{Propermut}
     }
\end{prop}
Let us end with the following lemma whose proof is close to Lemma $5$ from \cite{FCestimLp}  (see also Section 2.11 in \cite{Dbook}):
\begin{lem}
 \sl{For any $\aa, \beta>0$ there exists a constant $C_{\aa, \beta}>0$ such that for any $u\in \dot{H}^{s-\aa} \cap \dot{H}^{s+\beta}$, then $u\in\dot{B}_{2,1}^s$ and:
\begin{equation}
 \|u\|_{\dot{B}_{2,1}^s} \leq C_{\aa, \beta} \|u\|_{\dot{H}^{s-\aa}}^{\frac{\beta}{\aa + \beta}} \|u\|_{\dot{H}^{s+\beta}}^{\frac{\aa}{\aa + \beta}}.
\end{equation}
 }
\label{estimBsHs}
 \end{lem}

\subsection{Truncations}
\label{Troncatures}
In this section we define a particular truncation operator introduced in \cite{FCStratif1} that we will also abundantly use in the present article: let $\chi\in \mathcal{C}_0^\infty (\R, \R)$ taking values into $[0,1]$ and such that:
$$
\begin{cases}
 \mbox{supp }\chi \subset [-1,1],\\
 \chi \equiv 1 \mbox{ near } [-\frac12,\frac12].
\end{cases}
$$
Given $0<r<R$ we will denote by $\cC_{r,R}$ the following set (where $\xi=(\xi_h,\xi_3)$ and $\xi_h=(\xi_1,\xi_2)$):
\begin{equation}
\cC_{r,R} =\{\xi \in \R^3, \quad |\xi|\leq R \mbox{ and } |\xi_h|\geq r\}.
 \label{CrR}
\end{equation}
Defining $f_{r,R}(\xi)=\chi (\frac{|\xi|}{R})\big(1-\chi (\frac{|\xi_h|}{2r})\big)$, we have:
\begin{equation}
 \begin{cases}
 \mbox{supp }f_{r,R} \subset \cC_{r,R},\\
 f_{r,R}\equiv 1 \mbox{ on } \cC_{2r,\frac{R}2}.
\end{cases}
\label{frR}
\end{equation}
Let us define the following frequency truncation operator on $\cC_{r,R}$ ($\mathcal{F}^{-1}$ denotes the inverse Fourier transform and $|D|^s$ the classical derivation (non-local pseudo differential) operator: $|D|^s f =\mathcal{F}^{-1} (|\xi|^s \hat{f}(\xi))$.):
\begin{multline}
 \cP_{r,R} u= f_{r,R}(D) u =\chi (\frac{|D|}{R})\big(1-\chi (\frac{|D_h|}{2r})\big) u \\
 =\mathcal{F}^{-1} \Big(f_{r,R}(\xi) \hat{u}(\xi)\Big) = \mathcal{F}^{-1} \Big(\chi (\frac{|\xi|}{R})\big(1-\chi (\frac{|\xi_h|}{2r})\big) \hat{u}(\xi)\Big),
\label{PrR}
\end{multline}
Thanks to \eqref{frR}, we have:
\begin{equation}
 f_{\frac{r}2,2R}(D) f_{r,R}(D) u =f_{r,R}(D) u.
\end{equation}
In what follows we will use these objects, as in \cite{FC5, FCPAA, FCStratif1}, choosing in particular $r_\ee=\ee^m$ and $R_\ee =\ee^{-M}$, where $m$ and $M$ are precised in the proofs of the main results. Let us first recall the following anisotropic Bernstein-type result (more details in \cite{Dragos1, FC1, FCPAA, FCStratif1}):
\begin{lem}
\sl{There exists a constant $C>0$ such that for all function $f$, $\aa>0$, $1\leq q \leq p \leq \infty$ and all $0<r<R$, we have
\begin{equation}
\begin{cases}
 \vspace{0.2cm}
 \|\chi (\frac{|D|}{R}) \chi (\frac{|D_h|}{r}) f\|_{L^p} \leq C(R r^2)^{\frac{1}{q}-\frac{1}{p}} \|\chi (\frac{|D|}{R}) \chi (\frac{|D_h|}{r}) f\|_{L^q} \leq C(R r^2)^{\frac{1}{q}-\frac{1}{p}} \|f\|_{L^q}\\
 ||D|^\aa \cP_{r,R} f\|_{L^p} \leq C R^\aa \|\cP_{r,R} f\|_{L^p}.
\end{cases}
\end{equation}
}
\label{lemaniso}
\end{lem}
Let us end this section with the following proposition which adapts Lemma 2.3 from \cite{Dbook}. We refer to the last section of \cite{FCStratif1} for the proof.
\begin{prop} \cite{FCStratif1}
 \sl{Let $0<r<R$ be fixed. There exists a constant $C$ such that for any $p\in [1,\infty]$, $t\geq 0$ and any function $u$ we have:
 $$
 \mbox{Supp }\hat{u} \subset \cC_{r,R} \Rightarrow \|e^{t\D} u\|_{L^p} \leq C\frac{R^3}{r^4} e^{-\frac{t}2 r^2} \|u\|_{L^p}.
 $$
 }
 \label{Heatflow}
\end{prop}
\subsection{Eigenelements of the linearized system}

The linearized system of \eqref{Stratif} is written as follows (with $f_0, \Fe$ being divergence-free, the second form is obtained using the Leray orthogonal projector $\bP$ on divergence-free vectorfields):
\begin{equation}
\begin{cases}
\d_t f-(L-\frac{1}{\ee} \cB) f=\Fe,\\
\div f=0,\\
f_{|t=0}=f_0.
\end{cases}
\Longleftrightarrow
\begin{cases}
\d_t f-(L-\frac{1}{\ee} \bP \cB) f=\Fe,\\
f_{|t=0}=f_0.
\end{cases}
\label{systdisp}
\end{equation}
Applying the Fourier transform turns the equation into (as in \cite{FC1, Scro3}):
$$
\d_t \hat{f}(\xi)- \mathbb{B}(\xi, \ee)\hat{f}(\xi)=\hat{\Fe}(t,\xi),
$$
where
$$\mathbb{B}(\xi, \ee)= \hat{L-\frac{1}{\ee} \bP \cB} =\left(
\begin{array}{cccc}
-\nu(\xi_2^2+\xi_3^2) & \nu \xi_1 \xi_2 & \nu \xi_1 \xi_3 & \displaystyle{\frac{\xi_1\xi_3}{\ee |\xi|^2}}\\
\nu \xi_1 \xi_2 & -\nu(\xi_1^1+\xi_3^2) & \nu \xi_2 \xi_3 & \displaystyle{\frac{\xi_2\xi_3}{\ee |\xi|^2}}\\
\nu \xi_1 \xi_3 & \nu \xi_2 \xi_3 & -\nu(\xi_1^1+\xi_2^2) & \displaystyle{-\frac{\xi_1^2+\xi_2^2}{\ee |\xi|^2}}\\
0 & 0 & \displaystyle{\frac{1}{\ee}} & \displaystyle{-\nu'|\xi|^2}
\end{array}
\right).
$$

We refer to \cite{FCStratif1} for details about the following proposition gathering the properties needed to obtain the Strichartz estimates.
\begin{prop}
\label{estimvp}
\sl{
If $\nu\neq \nu'$, for all $m,M>0$ with $3M+m<1$, for all $\ee< \ee_1= \left(\frac{\sqrt{2}}{|\nu-\nu'|}\right)^\frac1{1-(3M+m)}$, if $\re=\ee^m$ and $\Re =\ee^{-M}$ (that is such that $|\nu-\nu'|\ee \Re^2\leq \re \sqrt{2}$), then for all $\xi \in \mathcal{C}_{r_{\ee}, R_{\ee}}$, the matrix $\mathbb{B}(\xi, \ee) = \widehat{L-\frac{1}{\ee} \bP \cB}$ is diagonalizable and its eigenvalues satisfy:
\begin{equation}
\label{vp}
\begin{cases}
\vspace{0.2cm} \lambda_1\exi=0,\\
\vspace{0.2cm} \lambda_2\exi=-\nu |\xi|^2,\\
\vspace{0.2cm} \lambda_3\exi = -\frac{\nu+\nu'}2 |\xi|^2+i\frac{|\xi_h|}{\ee|\xi|} -i \ee D\exi,\\
\vspace{0.2cm} \lambda_4\exi =\overline{\lambda_3\exi},
\end{cases}
\end{equation}
where $D\exi$ satisfies for all $\xi \in\cC_{\re, \Re}$ (with $k\in\{1,2\}$):
$$
\begin{cases}
\vspace{0.2cm}
|D\exi| \leq (\nu-\nu')^2\frac1{4\sqrt{2}} \frac{|\xi|^5}{|\xi_h|} \leq C_0 (\nu-\nu')^2 \frac{\Re^5}{\re}=C_0 (\nu-\nu')^2 \ee^{-(5M+m)},\\
\vspace{0.2cm}
|\d_{\xi_k}D\exi| \leq (\nu-\nu')^2\frac{9}{2\sqrt{2}} \frac{|\xi|^5}{|\xi_h|^2} \leq C_0 (\nu-\nu')^2 \frac{\Re^5}{\re^2}=C_0 (\nu-\nu')^2 \ee^{-(5M+2m)},\\
|\d_{\xi_3}D\exi| \leq (\nu-\nu')^2\frac{15}{4\sqrt{2}} \frac{|\xi|^4}{|\xi_h|} \leq C_0 (\nu-\nu')^2 \frac{\Re^4}{\re}=C_0 (\nu-\nu')^2 \ee^{-(4M+m)},
\end{cases}
$$
Moreover, if we denote by $\cP_k(\ee,\xi)$ the projectors onto the k-th eigenvector (we refer to the appendix from \cite{FCStratif1} for details), and by $\bP_k f=\bP_k(\ee,D) f=\cF^{-1}\left(\cP_k(\ee,\xi) \big(\hat{f}(\xi)\big)\right)$, then $\bP_2=\cQ$ (defined in Proposition \ref{PropSosc}) and for any divergence-free $\R^4$-valued vectorfield $f$, we have:
\begin{equation}
\begin{cases}
\vspace{0.2cm}
 \bP_2 f =(\n_h^\perp \D_h^{-1} \omega(f),0,0), \quad \mbox{with }\omega(f)=\d_1 f^2-\d_2 f^1,\\
 \|\bP_2 f\|_{\Hs} \leq \|(f^1,f^2)\|_{\Hs} \leq \|f\|_{\Hs}, \quad \mbox{for any }s\in \R.
\end{cases}
\label{estimvp2}
\end{equation}
and
\begin{equation}
\begin{cases}
\vspace{0.2cm}
 (I_d-\bP_2) f =(\n_h \D_h^{-1} \div_h f^h,f^3,f^3), \quad \mbox{with }\div_h f^h=\d_1 f^1+\d_2 f^2,\\
 \|(I_d-\bP_2) f\|_{\Hs} \leq \|f\|_{\Hs}, \quad \mbox{for any }s\in \R.
\end{cases}
\end{equation}
Finally for $k=3,4$,
\begin{equation}
\label{estimvp34}
\|\bP_k \cP_{\re, \Re} f\|_{\Hs }\leq \sqrt{2} \frac{\Re}{\re} \|\cP_{\re, \Re} f\|_{\Hs } =\sqrt{2} \ee^{-(m+M)} \|\cP_{\re, \Re} f\|_{\Hs }.
\end{equation}
If $\nu=\nu'$, there is no need anymore of a frequency truncation or an expansion for the last two eigenvalues (no $\ee_1$ either is necessary), and the $\bP_k$ ($k\in\{2,3,4\}$) are orthogonal so for any divergence-free $\R^4$-valued vectorfield $f$, we have:
$$
\|\bP_k f\|_{\Hs} \leq \|f\|_{\Hs}, \quad \mbox{for any }s\in \R.
$$}
\end{prop}

\subsection{Isotropic Strichartz estimates}

\label{Section:iso}
We list in this section the following Strichartz estimates that we proved in \cite{FCStratif1}: namely Propositions 5.4 and 5.6. We state the first one a little differently compared to \cite{FCStratif1} because we wish to estimate $f$ in the case when its $\bP_2$-part is zero (that is $f=\bP_3 f +\bP_4 f$). We recall that the operator $\cPrR$ is defined in Section \ref{Troncatures} and that we chose $(\re,\Re)=(\ee^m, \ee^{-M})$.

\begin{prop} (\cite{FCStratif1}, $\nu\neq \nu'$)
 \sl{For any $d\in \R$, $r\geq2$, $q\geq 1$ and $p\in[1,\frac{8}{1-\frac2{r}}]$, there exists a constant $C_{p,r}>0$ such that for any $\ee\in]0,\ee_1]$ (where $\ee_1= \left(\sqrt{2}/|\nu-\nu'|\right)^\frac1{1-(3M+m)}$) and any $f$ solving \eqref{systdisp} with initial data $f_0$ and external force $\Fe$ such that $\div f_0=\div \Fe=0$ and $\omega(f_0)=\omega(\Fe)=0$, then,
 \begin{multline}
  \||D|^d \cPrR f\|_{\tilde{L}_t^p\dot{B}_{r, q}^0}\\
  \leq \frac{C_{p,r}}{(\nu+\nu')^{\frac1{p}-\frac18(1-\frac2{r})}}  \frac{\Re^{7 -\frac9{r}}}{\re^{\frac{13}2+\frac2{p}-\frac7{r}}} \ee^{\frac18(1-\frac{2}{r})}\left( \|\cPrR f_0\|_{\dot{B}_{2, q}^d}+ \|\cPrR \Fe\|_{L^1 \dot{B}_{2, q}^d}\right).
 \end{multline}
\label{EstimStri1}
}
\end{prop}
When $\nu=\nu'$, usual simplifications allow better results: we have $L=\nu\Delta$ and System \eqref{systdisp} becomes:
\begin{equation}
\begin{cases}
\d_t f-(\nu \D-\frac{1}{\ee} \bP \cB) f=\Fe,\\
f_{|t=0}=f_0.
\end{cases}
\label{systdispb}
\end{equation}
\begin{prop} (\cite{FCStratif1}, $\nu= \nu'$)
 \sl{For any $d\in \R$, $r\geq 2$, $q\geq 1$, $\theta\in[0,1]$ and $p\in[1, \frac{4}{\theta (1-\frac{2}{r})}]$, there exists a constant $C=C_{p,r,\theta}>0$ such that for any $f$ solving \eqref{systdispb} for initial data $f_0$ and external force $\Fe$ such that $\div f_0=\div \Fe=0$ and $\omega(f_0)=\omega(\Fe)=0$, then,
 \begin{equation}
  \||D|^d f\|_{\tilde{L}_t^p\dot{B}_{r, q}^0} \leq \frac{C_{p,r,\theta}}{\nu^{\frac{1}{p}-\frac{\theta}{4}(1-\frac{2}{r})}} \ee^{\frac{\theta}{4}(1-\frac{2}{r})} \left( \|f_0\|_{\dot{B}_{2, q}^{\sigma_1}} +\|\Fe\|_{\tilde{L}_t^1 \dot{B}_{2, q}^{\sigma_1}} \right),
 \end{equation}
 where $\sigma_1= d+\frac32-\frac{3}{r}-\frac{2}{p}+\frac{\theta}{2} (1-\frac{2}{r})$.
\label{EstimStri2}
 }
\end{prop}

\subsection{Anisotropic Strichartz estimates}

\label{Section:aniso}
As observed in \cite{CDGG}, dealing with functions only depending on $x_3$ requires special versions of the Strichartz estimates: the space in $x$ now becomes of the form $L_{v,h}^{p, q}$ (anisotropic integrability in $x$) as introduced in \eqref{Anisospace}. We emphasize that, as described in Remark \ref{Anisorem}, the vertical/horizontal integrations are swapped compared to \cite{CDGG}.

The aim of this section is to state and prove the following anisotropic results:
\begin{prop} ($\nu \neq \nu'$)
 \sl{For any $m\geq2$, $p\in[1,\frac{8}{1-\frac2{m}}]$, there exists a constant $C_{p,m}>0$ such that for any $\ee\in]0,\ee_1]$ (where $\ee_1= \left(\sqrt{2}/|\nu-\nu'|\right)^\frac1{1-(3M+m)}$) and any $f$ solving \eqref{systdisp} with initial data $f_0$ and external force $\Fe$ such that $\div f_0=\div \Fe=0$ and $\omega(f_0)=\omega(\Fe)=0$, then
 \begin{multline}
  \|\cPrR f\|_{L_t^p L_{v,h}^{m,2}}\\
  \leq \frac{C_{p,m}}{(\nu+\nu')^{\frac1{p}-\frac18(1-\frac2{m})}} \frac{\Re^{6 -\frac7{m}}}{\re^{\frac{13}2+\frac2{p}-\frac7{m}}} \ee^{\frac18(1-\frac{2}{m})}\Big( \|\cPrR f_0\|_{L^2}+ \|\cPrR \Fe\|_{L^1 L^2}\Big),
 \end{multline}
\label{EstimStri1aniso}
}
\end{prop}
As usual, when $\nu=\nu'$ we can improve the previous result:
\begin{prop}($\nu= \nu'$)
 \sl{For any $d\in \R$, $m> 2$, $\theta\in[0,1]$ and $p\in[1, \frac{8}{\theta (1-\frac{2}{m})}]$, there exists a constant $C_{p,m,\theta}$ such that for any $f$ solving \eqref{systdispb} for initial data  $f_0$ and external force $\Fe$ such that $\div f_0=\div \Fe=0$ and $\omega(f_0)=\omega(\Fe)=0$, then
 \begin{equation}
  \||D|^d f\|_{L_t^p L_{v,h}^{m,2}} \leq \frac{C_{p,m,\theta}}{\nu^{\frac{1}{p}-\frac{\theta}{8}(1-\frac{2}{m})}} \ee^{\frac{\theta}{8}(1-\frac{2}{m})} \left( \|f_0\|_{\dot{B}_{2, q}^{\sigma_2}} +\|\Fe\|_{\tilde{L}_t^1 \dot{B}_{2, q}^{\sigma_2}} \right),
 \end{equation}
 where $\sigma_2= d+\frac12-\frac{1}{m}-\frac{2}{p}+\frac{\theta}{4} (1-\frac{2}{m})$.
\label{EstimStri2aniso}
 }
\end{prop}

\subsubsection{Proof of the anisotropic Strichartz estimates when $\nu\neq\nu'$}
The proof of Proposition \ref{EstimStri1aniso} is inspired by the one from \cite{CDGG} but, as in \cite{FCStratif1}, will require important adaptations. As usual we first assume $\Fe=0$ (and the inhomogeneous case is obtained reproducing the arguments on the Duhamel term). Starting close to what we did in \cite{FCStratif1}, we will skip details and point out what is new. Let $\cA$ be the following set:
$$
\cA\overset{def}{=}\{\psi \in \cC_0^\infty (\R_+\times \R^3, \R), \quad \|\psi\|_{L^{\bar{p}}(\R_+, L_{v,h}^{\bar{m},2}(\R^3))}\leq 1\}.
$$
As $\div f_0=\div \Fe=0$ and $\omega(f_0)=\omega(\Fe)=0$, we have $f=\bP_3 f +\bP_4 f$ so we can reduce to study $\bP_3 f$ (having in mind the norm of projectors $\bP_{3,4}$ given in Proposition \ref{estimvp}). Thanks to Plancherel and \eqref{PrR}, using the arguments from \cite{FCStratif1} (section 5.3.2)
\begin{multline}
 \|\bP_3 \cPrR f\|_{L_t^p L_{v,h}^{m,2}}= \sup_{\psi \in \cA} \int_0^\infty \int_{\R^3} \bP_3 \cPrR f (t,x) \psi(t,x) dx dt\\
 =C \sup_{\psi \in \cA} \int_0^\infty \int_{\R^3} e^{-\frac{\nu+\nu'}2 t|\xi|^2+i\frac{t}{\ee}\frac{|\xi_h|}{|\xi|}-i t \ee D\exi} \cF\big(\bP_3 \cPrR f_0\big)(\xi) f_{\frac{\re}2, 2\Re} (\xi)\hat{\psi}(t,\xi) d\xi dt\\
 \leq C \|\bP_3 \cPrR f_0\|_{L^2} \\
 \times \sup_{\psi \in \cA} \left[\int_0^\infty \int_0^\infty \|L_{\ee,t,t'} \psi(t,.)\|_{L_{v,h}^{m,2}} \|e^{\frac{\nu+\nu'}4 (t+t')\D} \overline{\cPrRb \psi (t',.)}\|_{L_{v,h}^{\bar{m},2}} dtdt'\right]^{\frac12},
 \label{estimTT2}
\end{multline}
where for some $g$:
\begin{equation}
 \left(L_{\ee,t,t'} g\right)(x)= \int_{\R^3} e^{ix\cdot \xi} e^{-\frac{\nu+\nu'}4 (t+t')|\xi|^2+i\frac{t-t'}{\ee}\frac{|\xi_h|}{|\xi|}-i(t-t') \ee D\exi} \chi (\frac{|\xi|}{2\Re})\big(1-\chi (\frac{|\xi_h|}{\re}) \hat{g}(\xi) d\xi.
\end{equation}
As in \cite{FCStratif1}, it is not possible to directly use the smoothing effect of the heat flow from Lemma 2.3 in \cite{Dbook} (Section 2.1.2), and we use Proposition \ref{Heatflow} which is an adaptation for the set $\cC_{r,R}$ (defined in \eqref{CrR}). The fact that in the present article, the spaces are anisotropic does not change the result as the bounds are obtained through convolution estimates, so that we obtain:
\begin{equation}
 \|e^{\frac{\nu+\nu'}4 (t+t')\D} \overline{\cPrRb \psi (t',.)}\|_{L_{v,h}^{\bar{m},2}} \leq C\frac{\Re^3}{\re^4} e^{-\frac{\nu+\nu'}{32} (t+t')\re^2} \|\psi(t',.)\|_{L_{v,h}^{\bar{m},2}}.
 \label{estimheat}
\end{equation}
The other term will require the Riesz-Thorin theorem, and thanks to \cite{FCStratif1} we already have:
\begin{equation}
 \|L_{\ee,t,t'}\|_{L_{v,h}^{2,2} \rightarrow L_{v,h}^{2,2}}\leq C_0 e^{-\frac{\nu+\nu'}{16} (t+t')\re^2},
\label{RieszTho}
\end{equation}
Obtaining a bound for $\|L_{\ee,t,t'}\|_{L_{v,h}^{1,2} \rightarrow L_{v,h}^{\infty,2}}$ will require us (as in \cite{CDGG}) to rewrite this operator. Let us first introduce the horizontal and vertical Fourier transforms: for a function $g$ depending on $x=(x_h,x_3)\in \R^3$,
$$
\cF_h g(\xi_h,x_3) \overset{def}{=} \int_{\R^2} e^{-i x_h\cdot \xi_h} g(x_h,x_3)dx_h, \quad \mbox{and}\quad \cF_v g(x_h,\xi_3) \overset{def}{=} \int_{\R} e^{-i x_3\cdot \xi_3} g(x_h,x_3)dx_3.
$$
Of course, $\cF =\cF_h \circ \cF_v =\cF_v \circ \cF_h$ and we easily obtain that, if we introduce:
\begin{equation}
 I_{\ee,t,t'} (\xi_h,x_3)=(2\pi)^{-1}\int_{\R} e^{ix_3\cdot \xi_3} e^{-\frac{\nu+\nu'}4 (t+t')|\xi|^2+i\frac{t-t'}{\ee}\frac{|\xi_h|}{|\xi|}-i(t-t') \ee D\exi} \chi (\frac{|\xi|}{2\Re})\big(1-\chi (\frac{|\xi_h|}{\re}) d\xi_3,
 \label{Def:I}
\end{equation}
then (also denoting as $\cF_v$ the vertical Fourier transform of a function depending on $(\xi_h,x_3)$):
\begin{multline}
\left(L_{\ee,t,t'} g\right)(x)\\
=\int_{\R^2} e^{ix_h\cdot \xi_h} \left(\int_{\R} e^{ix_3\cdot \xi_3} e^{-\frac{\nu+\nu'}4 (t+t')|\xi|^2+i\frac{t-t'}{\ee}\frac{|\xi_h|}{|\xi|}-i(t-t') \ee D\exi} \chi (\frac{|\xi|}{2\Re})\big(1-\chi (\frac{|\xi_h|}{\re}) \hat{g}(\xi) d\xi_3 \right) d\xi_h\\
=\int_{\R^2} e^{ix_h\cdot \xi_h} \left(\int_{\R} e^{ix_3\cdot \xi_3} \cF_v(I_{\ee,t,t'})(\xi_h,\xi_3)\cdot \cF_v \cF_h{g}(\xi_h,\xi_3) d\xi_3 \right) d\xi_h \\
= C \int_{\R^2} e^{ix_h\cdot \xi_h} (\cF_v)^{-1} \Big(\cF_v(I_{\ee,t,t'})(\xi_h,\xi_3)\cdot \cF_v \cF_h{g}(\xi_h,\xi_3)\Big) d\xi_h\\
= C \int_{\R^2} e^{ix_h\cdot \xi_h} \Big(I_{\ee,t,t'}(\xi_h,x_3) \ast_{x_3} (\cF_h g)(\xi_h,x_3) \Big)d\xi_h\\
=C \cF_h^{-1}\Big(I_{\ee,t,t'}(\xi_h,x_3) \ast_{x_3} (\cF_h g)(\xi_h,x_3)\Big).
\end{multline}
Thanks to \eqref{Anisospace}, the Plancherel, Minkowski and Young estimates, and to Remark 1.1 from \cite{Dragos1},
\begin{multline}
 \|L_{\ee,t,t'} g\|_{L_{v,h}^{\infty,2}} = C \|\cF_h^{-1}\Big(I_{\ee,t,t'}(\xi_h,x_3) \ast_{x_3} (\cF_h g)(\xi_h,x_3)\Big)\|_{L_{v,h}^{\infty,2}}\\
= C \sup_{x_3\in \R} \left(\int_{\R^2} |\cF_h^{-1}\Big(I_{\ee,t,t'}(\xi_h,x_3) \ast_{x_3} (\cF_h g)(\xi_h,x_3)\Big)|^2 dx_h\right)^\frac12\\
= C \sup_{x_3\in \R} \left(\int_{\R^2} |I_{\ee,t,t'}(\xi_h,x_3) \ast_{x_3} (\cF_h g)(\xi_h,x_3)|^2 d\xi_h\right)^\frac12\\
\leq C \left(\int_{\R^2} \Big\|I_{\ee,t,t'}(\xi_h,x_3) \ast_{x_3} (\cF_h g)(\xi_h,x_3)\Big\|_{L^\infty(\R_{x_3})}^2 d\xi\right)^\frac12\\
\leq C \Big\| \|I_{\ee,t,t'}(\xi_h,\cdot)\|_{L^\infty(\R_{x_3})} \|\cF_h g(\xi_h,.)\|_{L^1(\R_{x_3})}\Big\|_{L^2(\R_{\xi_h}^2)}\\
\leq  C \|I_{\ee,t,t'}(\xi_h,\cdot)\|_{L^\infty \big(\R_{\xi_h}^2, L^\infty(\R_{x_3})\big)} \|\cF_h g(\xi_h,.)\|_{L^2\big(\R_{\xi_h}^2,  L^1(\R_{x_3})\big)}\\
\leq C \|I_{\ee,t,t'}\|_{L^\infty (\R_{\xi_h}^2 \times \R_{x_3})} \|\cF_h g(\xi_h,.)\|_{L^1\big(\R_{x_3},L^2(\R_{\xi_h}^2)\big)} \leq C \|I_{\ee,t,t'}\|_{L^\infty} \|g(x_h,.)\|_{L^1\big(\R_{x_3},L^2(\R_{x_h}^2)\big)}\\
\leq C \|I_{\ee,t,t'}\|_{L^\infty} \|g\|_{L_{v,h}^{1,2}},
\end{multline}
which implies that
\begin{equation}
 \|L_{\ee,t,t'}\|_{L_{v,h}^{1,2} \rightarrow L_{v,h}^{\infty,2}}\leq C\|I_{\ee,t,t'}\|_{L^\infty}.
\label{RieszTho2}
\end{equation}
Thanks to \eqref{Def:I}, we immediately see that:
\begin{equation}
 \|I_{\ee,t,t'}\|_{L^\infty} \leq C_0 \Re e^{-\frac{\nu+\nu'}{16} (t+t')\re^2}.
 \label{normeLinfiniI1}
\end{equation}
In order to obtain a finer estimate, we will adapt the proof of Proposition $5.4$ from \cite{FCStratif1}: as $I_{\ee,t,t'}(\xi_h,-x_3)=I_{\ee,t,t'}(\xi_h,x_3)$ we can assume that $x_3\geq 0$. Moreover for any $t,t',\ee$,
$$
\|I_{\ee,t,t'}\|_{L^\infty(\R^3)}= \sup_{(\xi_h,x_3)\in \R^3} \|I_{\ee,t,t'}(\xi_h,\frac{t-t'}{\ee} x_3)\|,
$$
so that we will bound:
\begin{equation}
 I_{\ee,t,t'} (\xi_h,\frac{t-t'}{\ee} x_3)=(2\pi)^{-1}\int_{\R} e^{-\frac{\nu+\nu'}4 (t+t')|\xi|^2+i\frac{t-t'}{\ee} a(\xi)-i(t-t') \ee D\exi} \chi (\frac{|\xi|}{2\Re})\big(1-\chi (\frac{|\xi_h|}{\re}) d\xi_3,
 \label{Def:I2}
\end{equation}
where function $a$ is the same as in \cite{FCStratif1}:
$$
a(\xi) \overset{def}{=} x_3\cdot \xi_3+\frac{|\xi_h|}{|\xi|}.
$$
If we also introduce the same operator $\mathcal{L}$:
\begin{equation}
 \mathcal{L}f=
 \begin{cases}
 \vspace{0.2cm}
  \displaystyle{\frac1{1+\frac{t-t'}{\ee}\aa(\xi)^2} \left(f(\xi)+i\aa(\xi) \d_{\xi_3} f(\xi)\right)} & \mbox{if } t>t',\\
  \displaystyle{\frac1{1+\frac{t'-t}{\ee}\aa(\xi)^2} \left(f(\xi)-i\aa(\xi) \d_{\xi_3} f(\xi)\right)} & \mbox{else },
 \end{cases}
\end{equation}
with
$$
\aa(\xi)=-\d_{\xi_3} a(\xi)= -(x_3-\frac{\xi_3 |\xi_h|}{|\xi|^3}),
$$
then, performing an integration by parts, we obtain
\begin{multline}
 I_{\ee,t,t'}(\xi_h,\frac{t-t'}{\ee} x_3)\\
 =\int_{\R} e^{i\frac{t-t'}{\ee}a(\xi)}\big(1-\chi(\frac{|\xi_h|}{\re})\big) ^t\mathcal{L} \left(e^{-\frac{\nu+\nu'}4 (t+t')|\xi|^2-i(t-t') \ee D\exi} \chi(\frac{|\xi|}{2\Re})\right) d\xi_3.
\end{multline}
As the computation is the same as in \cite{FCStratif1}, we do not give details and only jump to the following bound:
\begin{equation}
 | ^t\mathcal{L} \left(e^{-\frac{\nu+\nu'}4 (t+t')|\xi|^2-i(t-t') \ee D\exi} \chi(\frac{|\xi|}{2\Re})\right)| \leq C_0 \frac{e^{-\frac{\nu+\nu'}{32} (t+t')\re^2}}{1+\frac{t-t'}{\ee} \aa^2} \left(\frac1{\re^2}+ \frac{|\aa|}{\re}\right),
\end{equation}
and
\begin{multline}
 |K_{\ee,t,t'}(\xi_h,\frac{t-t'}{\ee} x_3)|\\
 \leq C_0 \big|1-\chi(\frac{|\xi_h|}{\re})\big| e^{-\frac{\nu+\nu'}{32} (t+t')\re^2} \int_{-\sqrt{(2\Re)^2-|\xi_h|^2}}^{\sqrt{(2\Re)^2-|\xi_h|^2}} \frac1{1+\frac{t-t'}{\ee} \aa^2} \left(\frac1{\re^2}+ \frac{|\aa|}{\re}\right) d\xi_3.
 \label{normeLinfiniI2}
\end{multline}
We bounded a similar term in \cite{FCStratif1} (in the present article there is no horizontal integration) but will give a few details. It is easy to bound the second term using that $|\aa| =\left(\frac{t-t'}{\ee}\right)^{-\frac12} \left(\frac{t-t'}{\ee}\right)^{\frac12} |\aa|\leq  \frac12\left(\frac{t-t'}{\ee}\right)^{-\frac12}(1+\frac{t-t'}{\ee} \aa^2)$:
$$
 \int_{-\sqrt{(2\Re)^2-|\xi_h|^2}}^{\sqrt{(2\Re)^2-|\xi_h|^2}} \frac{|\aa|}{1+\frac{t-t'}{\ee} \aa^2} d\xi_3 \leq \frac12\left(\frac{t-t'}{\ee}\right)^{-\frac12} 4\Re.
$$
The first integral is split into two halves and the first half is easily bounded using the change of variable $z=\left(\frac{t-t'}{\ee}\right)^\frac12 \frac{\re}{16 \Re^3} \xi_3$:
$$
 \int_{-\sqrt{(2\Re)^2-|\xi_h|^2}}^0 \frac1{1+\frac{t-t'}{\ee} \aa^2} d\xi_3 \leq \int_{-\sqrt{(2\Re)^2-|\xi_h|^2}}^0 \frac1{1+\frac{t-t'}{\ee} \frac{\xi_3^2 \re^2}{16^2 \Re^6}} d\xi_3
 \leq C_0 \left(\frac{t-t'}{\ee}\right)^{-\frac12} \frac{\Re^3}{\re}.
$$
The most difficult part is to correctly bound the second half of the integral. In \cite{FCStratif1} we did it thanks to the following technical result:
\begin{prop}(\cite{FCStratif1} Proposition 6.1)
 \sl{There exists a constant $C_0>0$ such that for any $\aa>0$, $R\geq \frac2{\sqrt{3}}\aa$ and all $\bb\geq 0$,
 \begin{equation}
 I_{\aa, \bb}^R(\sigma) \overset{def}{=} \int_0^{\sqrt{R^2-\aa^2}} \frac{dx}{1+\sigma(f_\aa(x)-\bb)^2} \leq C_0 \frac{R^7}{\aa^\frac{11}2} \min(1,\sigma^{-\frac14}).
 \end{equation}
Moreover, the exponent $-\frac14$ is optimal in the sense that there exist $c_0,\sigma_0>0$ such that for any $R\geq \frac{\sqrt{3}}{\sqrt{2}} \aa$ and $\sigma\geq\sigma_0$,
 $$
 \sup_{\bb \in \R_+} I_{\aa, \bb}^R(\sigma) \geq c_0 \sigma^{-\frac14} \aa^\frac32.
 $$
 }
 \label{Proptech}
\end{prop}
This implies there exists a constant $C_0>0$ such that
$$
 \int_0^{\sqrt{(2\Re)^2-|\xi_h|^2}} \frac1{1+\frac{t-t'}{\ee} \aa^2} d\xi_3 =I_{|\xi_h|,x_3}^{2\Re}(\frac{t-t'}{\ee}) \leq C_0 \frac{\Re^7}{\re^\frac{11}2} \min\left(1, \big(\frac{t-t'}{\ee}\big)^{-\frac14}\right),
$$
which finally leads to:
$$
\|I_{\ee,t,t'}\|_{L^\infty}\leq C_0 \frac{\Re^7}{\re^\frac{15}2} \min\left(1,\frac{\ee^\frac14}{|t-t'|^\frac14}\right) e^{-\frac{\nu+\nu'}{32} (t+t')\re^2} \leq C_0 \frac{\Re^7}{\re^\frac{15}2} \frac{\ee^\frac14}{|t-t'|^\frac14} e^{-\frac{\nu+\nu'}{32} (t+t')\re^2}.
$$
Using this together with \eqref{RieszTho} and \eqref{RieszTho2}, we obtain thanks to the Riesz-Thorin theorem that for any $r\in[2, \infty]$:
$$
 \|L_{\ee,t,t'} g\|_{L_{v,h}^{m,2}} \leq C_0 \left(\frac{\Re^7}{\re^\frac{15}2} \frac{\ee^\frac14}{|t-t'|^\frac14}\right)^{1-\frac2{m}} e^{-\frac{\nu+\nu'}{32} (t+t')\re^2} \|g\|_{L_{v,h}^{\bar{m},2}}.
$$
Gathering this estimates together with \eqref{estimheat}, and thanks to \eqref{estimvp34}, we can properly bound \eqref{estimTT2} and obtain that:
\begin{multline}
  \|\bP_3 \cPrR f\|_{L_t^p L_{v,h}^{m,2}}\\
  \leq C_0 \|\cPrR f_0\|_{L^2} \sup_{\psi \in \cA} \frac{\Re^{1+\frac32 +\frac72 (1-\frac2{m})}}{\re^{1+2+\frac{15}4 (1-\frac2{m})}} \ee^{\frac18 (1-\frac2{m})} \left[\int_0^\infty \int_0^\infty \frac{h(t)h(t')}{|t-t'|^{\frac14(1-\frac2{m})}} dtdt'\right]^{\frac12},
\end{multline}
where $h(t)=e^{-\frac{\nu+\nu'}{16} t\re^2} \|\psi(t,.)\|_{L_{v,h}^{\bar{m},2}}$. The rest of the proof is identical to \cite{FCStratif1} (end of Section 5.3.2) so we directly write the bound:
$$
 \|\bP_3 \cPrR f\|_{L_t^p L_{v,h}^{m,2}} \leq \frac{C_{p,m}}{(\nu+\nu')^{\frac1{p}-\frac18(1-\frac2{m})}} \frac{\Re^{6-\frac7{m}}}{\re^{\frac{13}2+\frac2{p}-\frac7{m}}} \ee^{\frac18 (1-\frac2{m})} \|\cPrR f_0\|_{L^2},
$$
where $C_{p,m}=C_0 \big[16(\frac1{p}-\frac18(1-\frac2{m}))\big]^{\frac1{p}-\frac18(1-\frac2{m})}$ which concludes the proof. $\blacksquare$

\subsubsection{Proof of the anisotropic Strichartz estimates when $\nu=\nu'$}
As in the previous section, we are reduced to study $\bP_3 f$ in the case  $\Fe=0$, but when $\nu=\nu'$ additionnal simplifications arise (described in Proposition \ref{estimvp}):
\begin{itemize}
 \item The projectors $\bP_{3,4}$ become mutually orthogonal (we recall that in the general case they are orthogonal to $\bP_2$) and their norms become 1,
 \item Frequency truncations are not needed anymore for the eigenvalues (and projectors) in the case $k\in\{3,4\}$, and we can consider $\bP_k f$ (instead of $\bP_k \cPrR f$ in the previous part).
\end{itemize}
Nevertheless, to prove Proposition \ref{EstimStri2aniso}, we will begin as in \cite{FCPAA, FCRF} by frequency localization (we refer to Section A2 from \cite{FCPAA} for the notations related to the Besov spaces, and more generally to \cite{Dbook} for a complete presentation of the Littlewood-Paley theory). Introducing the complete truncation operator $\ddj u =\varphi(2^{-j} |D|)u$ and its horizontal counterpart $\ddk u=\varphi(2^{-k} |D_h|)u$:
$$
\|\bP_3 \ddj f\|_{L_t^p L_{v,h}^{m,2}} \leq \Sum_{k=-\infty}^{j+1} \|\bP_3 \ddk \ddj f\|_{L_t^p L_{v,h}^{m,2}},
$$
and
\begin{multline}
 \|\bP_3 \ddk \ddj f\|_{L_t^p L_{v,h}^{m,2}} = \sup_{\psi \in \cA} \int_0^\infty \int_{\R^3} \bP_3 \ddk \ddj f (t,x) \psi(t,x) dx dt\\
 \leq C \|\bP_3 \ddk \ddj f_0\|_{L^2} \\
 \times \sup_{\psi \in \cA} \left[\int_0^\infty \int_0^\infty \|L_{\ee,t,t'}^{j,k} \psi(t,.)\|_{L_{v,h}^{m,2}} \|e^{\frac{\nu}2 (t+t')\D} \varphi_1(2^{-j} |D|) \varphi_1(2^{-k} |D_h|) \overline{\psi (t',.)}\|_{L_{v,h}^{\bar{m},2}} dtdt'\right]^{\frac12},
 \label{estimTT2nunu}
\end{multline}
where $\varphi_1$ is a function (with values in $[0,1]$) supported in the set $\cC'=[c_0, C_0]$ (say $(c_0, C_0)=(\frac35,3)$) and equal to 1 close to $\cC=[\frac34, \frac83]$ (introduced in the first section of the appendix), and for some $g$ we define the analoguous of the operator $L_{\ee,t,t'}$ from the previous section:
\begin{equation}
 \left(L_{\ee,t,t'}^{j,k} g\right)(x)= \int_{\R^3} e^{ix\cdot \xi} e^{-\frac{\nu}2 (t+t')|\xi|^2+i\frac{t-t'}{\ee}\frac{|\xi_h|}{|\xi|}} \varphi_1(2^{-j} |\xi|) \varphi_1(2^{-k} |\xi_h|) \hat{g}(\xi) d\xi.
\end{equation}
The heat term is estimated without resorting to Proposition \ref{Heatflow} thanks to the following fact: introducing $h_1(x_h)=\cF_h^{-1} \varphi_1(|\xi_h|)$, for any $p,q\in [1,\infty]$ and any function $g$ we have:
\begin{multline}
 \|\ddk \ddj f\|_{L_{v,h}^{p,q}}^p\\
 =\int_\R \left(\int_{R^2} |\ddk \ddj f(x_h,x_3)|^q dx_h \right)^{\frac{p}{q}} dx_3 =\int_R \|2^{2k} h_1(2^k |.|) \star_{x_h}\ddj f(.,x_3)\|_{L^q(\R_h^2)}^p dx_3\\
 \leq \int_R \left(\|h_1\|_{L^1(\R_h^2)} \|\ddj f(.,x_3)\|_{L^q(\R_h^2)}\right)^p dx_3 =\|h_1\|_{L^1(\R_h^2)}^p \|\ddj f\|_{L_{v,h}^{p,q}}^p,
\end{multline}
so that (as explained in the previous section, obtaining an anisotropic version of Lemma 2.3 from \cite{Dbook} is easy as the proof involves convolutions) there exists a constant $C>0$ such that ($c_0=3/5$ as recalled above):
\begin{multline}
 \|e^{\frac{\nu}2 (t+t')\D} \varphi_1(2^{-j} |D|) \varphi_1(2^{-k} |D_h|) \overline{\psi (t',.)}\|_{L_{v,h}^{\bar{m},2}} \leq \|e^{\frac{\nu}2 (t+t')\D} \varphi_1(2^{-j} |D|) \overline{\psi (t',.)}\|_{L_{v,h}^{\bar{m},2}}\\
 \leq C e^{-\frac{\nu}2 (t+t')c_0^2 2^{2j}} \|\psi (t',.)\|_{L_{v,h}^{\bar{m},2}}.
\end{multline}
With a view to use the Riesz-Thorin theorem, similarly as in the previous section we have:
\begin{equation}
 \|L_{\ee,t,t'}^{j,k}\|_{L_{v,h}^{2,2} \rightarrow L_{v,h}^{2,2}}\leq C_0 e^{-\frac{\nu}2 (t+t')c_0^2 2^{2j}}.
\label{RieszThonunu}
\end{equation}
Introducing
\begin{equation}
 I_{\ee,t,t'}^{j,k} (\xi_h,x_3)=(2\pi)^{-1}\int_{\R} e^{ix_3\cdot \xi_3} e^{-\frac{\nu}2 (t+t')|\xi|^2+i\frac{t-t'}{\ee}\frac{|\xi_h|}{|\xi|}} \varphi_1(2^{-j} |\xi|) \varphi_1(2^{-k} |\xi_h|) d\xi_3,
 \label{Def:Inunu}
\end{equation}
and reproducing the arguments from the previous section leads to
$$
\left(L_{\ee,t,t'}^{j,k} g\right)(x) =C \cF_h^{-1}\Big(I_{\ee,t,t'}^{j,k}(\xi_h,x_3) \ast_{x_3} (\cF_h g)(\xi_h,x_3)\Big),
$$
and
$$
\|L_{\ee,t,t'}^{j,k} g\|_{L_{v,h}^{\infty,2}} \leq \|I_{\ee,t,t'}^{j,k}\|_{L^\infty} \|g\|_{L_{v,h}^{1,2}},
$$
so that (thanks to \eqref{Def:Inunu}):
\begin{equation}
 \|L_{\ee,t,t'}^{j,k}\|_{L_{v,h}^{1,2} \rightarrow L_{v,h}^{\infty,2}} \leq \|I_{\ee,t,t'}^{j,k}\|_{L^\infty} \leq C_0 2^j e^{-\frac{\nu}2 (t+t') c_0^2 2^{2j}}.
 \label{RieszTho2nunu}
\end{equation}
Next, performing the change of variable $\xi_3=2^j \eta_3$, we can write that
$$
I_{\ee,t,t'}^{j,k} (2^j \eta_h,x_3)=2^j \overline{I_{\ee,t,t'}^{j,k}} (\eta_h,2^j x_3),
$$
where
\begin{equation}
 \overline{I_{\ee,t,t'}^{j,k}} (\eta_h,x_3)=(2\pi)^{-1}\int_{\R} e^{ix_3\cdot \eta_3} e^{-\frac{\nu}2 2^{2j}(t+t')|\eta|^2+i\frac{t-t'}{\ee}\frac{|\eta_h|}{|\eta|}} \varphi_1(|\eta|) \varphi_1(2^{j-k} |\eta_h|) d\eta_3,
 \label{Def:Inunub}
\end{equation}
which entails that
\begin{equation}
 \|I_{\ee,t,t'}^{j,k}\|_{L^\infty} =2^j \|\overline{I_{\ee,t,t'}^{j,k}}\|_{L^\infty}.
\end{equation}
In $\overline{I_{\ee,t,t'}^{j,k}}$, the frequencies are now truncated as follows: $c_0\leq |\eta| \leq C_0$ and $c_0 2^{k-j}\leq |\eta_h| \leq C_0 2^{k-j}$, so that we can reproduce the arguments from the previous section (see also \cite{FCStratif1}) to the vertical rescale $\overline{I_{\ee,t,t'}^{j,k}}(\eta_h,\frac{t-t'}{\ee} x_3)$ with $(\re,\Re)$ replaced by $(C_0,c_0 2^{k-j})$ and obtain that (also using \eqref{RieszTho2nunu}) for all $|\eta_h|\geq c_0 2^{k-j}$ and $x_3$,
$$
|\overline{I_{\ee,t,t'}^{j,k}}(\eta_h,\frac{t-t'}{\ee} x_3)| \leq C_0  \int_{-\sqrt{C_0^2-|\eta_h|^2}}^{\sqrt{C_0^2-|\eta_h|^2}} \frac{e^{-\frac{\nu}4 (t+t') c_0^2 2^{2j}}}{1+\frac{|t-t'|}{\ee} \aa(\eta)^2}\left((1+\frac4{|\eta|^2}) +|\aa(\eta)|(\frac1{|\eta|}+1)\right) d\eta_3,
$$
so that (we recall that in the present case $|\eta|\geq c_0$, which is better than in the previous section) with the same steps as in the previous part (see also Section 5.3.2 and Proposition 6.1 from \cite{FCStratif1}, which is recalled in the present article as Proposition \ref{Proptech}), for every $\theta\in[0,1]$
\begin{multline}
  \|I_{\ee,t,t'}^{j,k}\|_{L^\infty} \leq C_0 2^j \frac{C_0^7}{(c_0 2^{k-j})^\frac{11}2} e^{-\frac{\nu}4 (t+t') c_0^2 2^{2j}} \min \left(1, \frac{\ee^\frac14}{|t-t'|^\frac14}\right)\\
  \leq C_0 2^j \frac{C_0^7}{(c_0 2^{k-j})^\frac{11}2} e^{-\frac{\nu}4 (t+t') c_0^2 2^{2j}} \frac{\ee^\frac{\theta}4}{|t-t'|^\frac{\theta}4}.
\end{multline}
\begin{rem}
 \sl{As in the previous section the most difficult is to correctly bound the following integral, which is done using Proposition \ref{Proptech}:
 $$
 \int_0^{\sqrt{C_0^2-|\eta_h|^2}} \frac1{1+\frac{|t-t'|}{\ee} \aa(\eta)^2} d\xi_3 =I_{|\xi_h|, x_3}^{C_0} (\frac{t-t'}{\ee}) \leq C_1 \frac{C_0^7}{(c_0 2^{k-j})^\frac{11}2} \min\left(1, \big(\frac{t-t'}{\ee}\big)^{-\frac14}\right).
 $$
 }
\end{rem}
Gathering the previous bound for $\|I_{\ee,t,t'}^{j,k}\|_{L^\infty}$ with \eqref{RieszThonunu}, thanks to the Riesz-Thorin theorem we finally obtain that with $m\geq 2$ and for any $\theta\in[0,1]$,
$$
 \|L_{\ee,t,t'}^{j,k} \psi(t,\cdot)\|_{L_{v,h}^{m,2}} \leq C_0 e^{-\frac{\nu}4 (t+t') c_0^2 2^{2j}} \left(2^{j+\frac{11}2(j-k)} \frac{\ee^\frac{\theta}4}{|t-t'|^\frac{\theta}4}\right)^{1-\frac2{m}} \|\psi(t)\|_{L_{v,h}^{\bar{m},2}}.
$$
Plugging this into \eqref{estimTT2nunu}, we obtain
\begin{multline}
 \|\bP_3 \ddk \ddj f\|_{L_t^p L_{v,h}^{m,2}}\\ \leq C_0 \|\bP_3 \ddj f_0\|_{L^2} 2^{\big(j+\frac{11}2(j-k)\big)(\frac12-\frac1{m})} \ee^{\frac{\theta}8 (1-\frac2{m})} \sup_{\psi \in \cA} \left[\int_0^\infty \int_0^\infty \frac{h(t)h(t')}{|t-t'|^{\frac{\theta}4(1-\frac2{m})}} dtdt'\right]^{\frac12},
\end{multline}
with $h(t)=e^{-\frac{3\nu}4 c_0^2 2^{2j}t} \|\psi(t,.)\|_{L_{v,h}^{\bar{m},2}}$. Using once more the Hardy-Littlewood theorem, and introducing $k_1,\bb\geq 1$ defined as follows (the condition on $p$ comes from here)
$$
\frac1{k_1}=1-\frac{\theta}8(1-\frac2{m}), \quad \mbox{and} \quad \frac1{\bb}=\frac1{p}-\frac{\theta}8(1-\frac2{m}),
$$
we obtain that:
$$
\left[\int_0^\infty \int_0^\infty \frac{h(t)h(t')}{|t-t'|^{\frac{\theta}4(1-\frac2{m})}} dtdt'\right]^{\frac12} \leq C \|h\|_{L^{k_1}(\R)} \leq C \left(\int_0^\infty e^{-\frac{3\nu}4 c_0^2 2^{2j} \bb t}\right)^\frac1{\bb} \|\psi\|_{L^{\bar{p}} L_{v,h}^{\bar{m},2}},
$$
and
\begin{multline}
 \|\bP_3 \ddk \ddj f\|_{L_t^p L_{v,h}^{m,2}} \leq \frac{C_{p,m,\theta}}{\nu^{\frac1{p}-\frac{\theta}8(1-\frac2{m})}} \|\bP_3 \ddj f_0\|_{L^2} \ee^{\frac{\theta}8 (1-\frac2{m})} 2^{\big(j+\frac{11}2(j-k)\big)(\frac12-\frac1{m})-\frac2{\bb}}\\
 \leq \frac{C_{p,m,\theta}}{\nu^{\frac1{p}-\frac{\theta}8(1-\frac2{m})}} \|\bP_3 \ddj f_0\|_{L^2} \ee^{\frac{\theta}8 (1-\frac2{m})} 2^{\frac{11}2(j-k)(\frac12-\frac1{m})} 2^{j\big(\frac12-\frac1{m}-\frac2{p}+\frac{\theta}4 (1-\frac2{m})\big)}.
\end{multline}
Now, all that remains is to sum for $k\leq j+1$, which is possible if and only if $m>2$, and:
$$
 \|\bP_3 \ddj f\|_{L_t^p L_{v,h}^{m,2}} \leq \frac{C_{p,m,\theta}}{\nu^{\frac1{p}-\frac{\theta}8(1-\frac2{m})}} \|\bP_3 \ddj f_0\|_{L^2} \ee^{\frac{\theta}8 (1-\frac2{m})} 2^{j\big(\frac12-\frac1{m}-\frac2{p}+\frac{\theta}4 (1-\frac2{m})\big)}.
$$
Multiplying by $2^{jd}$ and summing over $j\in \Z$, we conclude the proof. $\blacksquare$

\subsection{Proof of theorem \ref{TH0FK}}
\label{PreuveFK}
Let us fix some $\ee>0$. We will use the Friedrich's scheme introduced in \cite{FCStratif1}. If $\Den=(\Ven,\Hen)$, projecting over divergence-free vectorfields with the Leray projector $\bP$, this system is written as follows (for $n\in\N$, $J_n$ is the Fourier truncation operator on the ball centered at zero and with radius $n$):
\begin{equation}
\begin{cases}
 \d_t \Den-L \Den +\frac1{\ee} \bP \cB J_n\Den = -J_n \bP \left[\Den\cdot \n \Den
 +\left(\begin{array}{c}
  \Den \cdot \n \tv^h \\0\\ D_{\ee}^{n,3}\cdot \d_3 \tThee \end{array}\right) +\tv^h\cdot \n_h \Den \right] +J_n \tG,\\
 D_{\ee|t=0}^n=J_n\left(\Uoeosc +(\UoeS-(\tvo^h,0,0))\right) =J_n\left(\Uoeosc +(\UoeS^h-\tvo^h,0,0)\right).
\end{cases}
\label{SchemasF}
\end{equation}
In order to neutralize the constant term $\tG$, let us introduce the following Stokes-type system:
\begin{equation}
\begin{cases}
 \d_t \Ee-L \Ee +\frac1{\ee} \cB \Ee = -\left(\begin{array}{c}\n Q_\ee \\0\end{array}\right) +\tG,\\
 E_{\ee|t=0}=\Uoeosc +(\UoeS^h-\tvo^h,0,0).
\end{cases}
\label{SystE}
\end{equation}
It is easy to prove that if $\Ee(0) \in \dot{H}^s$ for some $s\in[0,\frac12+\delta]$, there exists a unique global solution satisfying for all $t\geq 0$:
\begin{equation}
 \|\Ee(t)\|_{\dot{H}^s}^2 +\nu_0 \int_0^t \|\n \Ee(\tau)\|_{\dot{H}^s}^2 d\tau \leq \left(\|\Ee(0)\|_{\dot{H}^s}^2 +\int_0^t \|\tG(\tau)\|_{\dot{H}^s} d\tau \right) e^{\int_0^t \|\tG(\tau)\|_{\dot{H}^s} d\tau}.
 \label{EstimEe}
\end{equation}
Now, we introduce $\Een=J_n \Ee$ and $\Fen=\Den-\Een$ which satisfy $F_{\ee|t=0}^n=0$ and:
\begin{multline}
 \d_t \Fen-L \Fen +\frac1{\ee} \bP \cB J_n\Fen = -J_n \bP \Bigg[
 \left(\begin{array}{c}
  (\Een+\Fen) \cdot \n \tv^h \\0\\ (E_{\ee}^{n,3}+F_{\ee}^{n,3})\cdot \d_3 \tThee \end{array}\right)  +\tv^h\cdot \n_h (\Een +\Fen)\\
  +\Fen\cdot \n \Fen +\Fen\cdot \n \Een +\Een\cdot \n \Fen +\Een\cdot \n \Een\Bigg],\\
\end{multline}
Reproducing the arguments from Section \ref{Apriorigen} we obtain that
\begin{multline}
\frac{d}{dt} \|\Fen\|_{L^2}^2 +\nu_0 \|\n \Fen\|_{L^2}^2 \leq \frac{C}{\nu_0} \|\Fen\|_{L^2}^2 \Bigg(\|\n \Ee\|_{\dot{H}^\frac12}^2 +\|\n \tv^h\|_{\dot{H}^\frac12}^2 +\nu_0^\frac23 \|\tThee\|_{\dot{H}^1(\R)}^\frac43 \Bigg)\\
 +C \Bigg[\frac1{\nu_0} \|\Ee\|_{\dot{H}^\frac12}^2 (\|\n \Ee\|_{L^2}^2 +\|\n \tv^h\|_{L^2}^2) +(1+\frac1{\nu_0}\|\tv^h\|_{\dot{H}^\frac12}^2)\|\n \Ee\|_{L^2}^2 +\|\tThee\|_{\dot{H}^1(\R)}^\frac43 \|\Ee\|_{L^2}^2\Bigg]
\end{multline}
and if $\frac12\leq s < \frac12+\frac{\delta}2$ (for some $\bb>0$)
\begin{multline}
\frac{d}{dt} \|\Fen\|_{\dot{H}^s}^2 +\nu_0 \|\Fen\|_{\dot{H}^s}^2 \leq 2C \|\Fen\|_{\dot{H}^\frac12} \|\n\Fen\|_{\dot{H}^s}^2 +\frac{C}{\nu_0} \|\Fen\|_{H^{\frac12-\bb}}^2 \|\tThee\|_{\dot{H}^{s+\bb}(\R)}^2  +\frac{\nu_0}4 \|\n \Fen\|_{H^s}^2\\
 +\frac{C}{\nu_0} \|\Fen\|_{\dot{H}^s}^2 \Bigg( (1+\frac1{\nu_0^2} \|\tv^h\|_{\dot{H}^\frac12}^2) \|\n \tv^h\|_{\dot{H}^\frac12}^2 +(1+\frac1{\nu_0^2} \|\Ee\|_{\dot{H}^\frac12}^2) \|\n \Ee\|_{\dot{H}^\frac12}^2\Bigg)\\
 + \frac{C}{\nu_0} \Bigg[\|\Ee\|_{\dot{H}^s}^2(\|\n \Ee\|_{\dot{H}^\frac12}^2 +\|\n \tv^h\|_{\dot{H}^\frac12}^2) +\|\tv^h\|_{\dot{H}^s}^2 \|\n \Ee\|_{\dot{H}^\frac12}^2 +\|\Een\|_{H^{\frac12-\bb}}^2 \|\tThee\|_{\dot{H}^{s+\beta}(\R)}^2\Bigg].
\end{multline}
If we define
$$
T_\ee^{n,1}=\sup \{t>0,\; \forall t'\in[0,t], \; \|\Fen(t')\|_{\dot{H}^\frac12} \leq \frac{\nu_0}{4C}\},
$$
then we have $T_\ee^{n,1}>0$ (we recall that $\Fen(0)=0$) and thanks to the previous estimates with $s=\frac12$, Theorem \ref{ThSNS} and \eqref{EstimEe}, there exists a constant $\mathbb{D}=\mathbb{D}(\Co, \delta, \nu,\nu', \|\Uoeosc\|_{H^\frac12})$ for all $t\leq T_\ee^{n,1}$,
\begin{multline}
 \|\Fen(t)\|_{H^\frac12}^2 +\frac{\nu_0}2 \int_0^t \|\n \Fen(\tau)\|_{H^\frac12}^2 d\tau\\
 \leq \Do \left[\|\n\Ee\|_{L_t^2 L^2}^2 +\|\n\tv^h\|_{L_t^2 L^2}^2 +\|\n\Ee\|_{L_t^2 \dot{H}^\frac12}^2 +\|\n\tv^h\|_{L_t^2 \dot{H}^\frac12}^2 +\|\tThee\|_{L_t^\frac43 \dot{H}^1}^\frac43 +\|\tThee\|_{L_t^2 \dot{H}^{\frac12+\bb}}^2 \right].
\end{multline}
All these quantities are independant of $n$ and converge to zero as $t$ goes to zero, so if we define $T_\ee^2>0$ such that the right-hand side is bounded by $\left(\frac{\nu_0}{8C}\right)^2$ for all $t\leq T_\ee^2$, then with the same arguments as in the previous bootstrap, we obtain that $T_\ee^{n,1} \geq T_\ee^2>0$ for all $n$,  and for all $t\leq T_\ee^2$,
$$
\|\Fen(t)\|_{H^\frac12}^2 +\frac{\nu_0}2 \int_0^t \|\n \Fen(\tau)\|_{H^\frac12}^2 d\tau \leq \mathbb{D}'(\Co, \delta, \nu,\nu', \|\Uoeosc\|_{H^\frac12}),
$$
which allows to prove (with classical arguments) existence of a strong solution as described in Theorem \ref{TH0FK}.

The propagation of the regularity and the blow-up criterion are proved through classical ideas thanks to the following estimates (which are proved with the very same arguments): for all $s\in[0,\frac12+\delta]$ and $t\in[0, T_\ee^*[$,
\begin{multline}
 \|\De(t)\|_{H^s}^2 +\frac{\nu_0}2 \int_0^t \|\n \De(\tau)\|_{H^s}^2 d\tau \leq \left[\|\Uoeosc\|_{H^s}^2 +\|\UoeS^h-\tvo^h\|_{H^s}^2 +\|\tG\|_{L_t^1 H^s}\right]\\
 \times e^{\frac{C}{\nu_0}\left[\|\n \De\|_{L_t^2 \dot{H}^\frac12}^2 +\|\n \tv^h\|_{L_t^2 \dot{H}^\frac12}^2 (1+\frac1{\nu_0^2} \|\tv^h\|_{L_t^\infty \dot{H}^\frac12}^2) +\nu_0^\frac23 \|\tThee\|_{L_t^\frac43 \dot{H}^1}^\frac43 +\|\tThee\|_{L_t^2 \dot{H}^{\frac12+\bb}}^2\right]},
\end{multline}
which ends the proof of Theorem \ref{TH0FK}. $\blacksquare$

\textbf{Aknowledgement :} The author wishes to thank the anonymous referees for useful suggestions.

%\textbf{Aknowledgements :}

\end{document}